\pdfoutput=1
\documentclass{article}

\usepackage{graphicx}
\usepackage{amssymb}
\usepackage{amsmath}
\usepackage[all]{xypic}
\usepackage{rotating}   
\usepackage{hyperref}
\xyoption{poly}
\entrymodifiers={+!!<0pt,\fontdimen22\textfont2>}    

\numberwithin{equation}{section}
\numberwithin{figure}{section}
\numberwithin{table}{section}

\newcommand{\abs}[1]{\left|#1\right|}
\newcommand{\vnorm}[1]{\left|\left|#1\right|\right|}
\newcommand{\im}[1]{\mbox{im}#1}
\newcommand{\Z}{\mathbb{Z}}
\newcommand{\R}{\mathbb{R}}
\newcommand{\C}{\mathbb{R}}
\long\def\symbolfootnote[#1]#2{\begingroup%
\def\thefootnote{\fnsymbol{footnote}}\footnote[#1]{#2}\endgroup}

\newcommand{\foursqwidth}{0.45} 
\newcommand{\foursqwidtha}{.6} 
\newcommand{\threewidth}{1}   
\newcommand{\onewidth}{0.9}     
\newcommand{\tabonewd}{13cm}    
\newcommand{\tabtwowd}{20cm}    
\newcommand{\tabthreewd}{13cm}  
\newcommand{\tabthreewda}{11cm} 
\newcommand{\tabfourwd}{13cm}   
\newcommand{\tabfourwda}{14cm}  
\newcommand{\tabfivewd}{13cm}   
\newcommand{\tabfivewda}{11.3cm}
\newcommand{\tabsixwd}{13cm}    
\newcommand{\tabsevenwd}{9cm}   
\newcommand{\tabeightwd}{8cm}   
\newcommand{\tabninewd}{10cm}   
\newcommand{\tabtenwd}{7cm}     

\def\alphanum{\ifcase\xypolynode\or \makebox[.8cm][r]{\Huge -$v$\:}\or \makebox[.8cm][r]{\Huge -$v$\:}\or \makebox[.8cm][r]{\Huge $u$\:}\or \makebox[.8cm][r]{\Huge $u$\:}\or \makebox[.8cm][r]{\Huge $v$\:}\or \makebox[.8cm][r]{\Huge $v$\:}\or \makebox[.8cm][r]{\Huge -$u$\:}\or \makebox[.8cm][r]{\Huge -$u$\:}\fi}
\def\betanum{\ifcase\xypolynode\or \makebox[.8cm][r]{\Huge -$u$\:}\or \makebox[.8cm][r]{\Huge -$u$\:}\or \makebox[.8cm][r]{\Huge $v$\:}\or \makebox[.8cm][r]{\Huge $v$\:}\or \makebox[.8cm][r]{\Huge $u$\:}\or \makebox[.8cm][r]{\Huge $u$\:}\or \makebox[.8cm][r]{\Huge -$v$\:}\or \makebox[.8cm][r]{\Huge -$v$\:}\fi}
\def\gammanum{\ifcase\xypolynode\or \makebox[.8cm][r]{\Huge -$u$\:}\or \makebox[.8cm][r]{\Huge $u$\:}\or \makebox[.8cm][r]{\Huge $u$\:}\or \makebox[.8cm][r]{\Huge -$u$\:}\or \makebox[.8cm][r]{\Huge -$u$\:}\or \makebox[.8cm][r]{\Huge $u$\:}\or \makebox[.8cm][r]{\Huge $u$\:}\or \makebox[.8cm][r]{\Huge -$u$\:}\fi}
\def\deltanum{\ifcase\xypolynode\or \makebox[.8cm][r]{\Huge -$v$\:}\or \makebox[.8cm][r]{\Huge -$v$\:}\or \makebox[.8cm][r]{\Huge $v$\:}\or \makebox[.8cm][r]{\Huge $v$\:}\or \makebox[.8cm][r]{\Huge -$v$\:}\or \makebox[.8cm][r]{\Huge -$v$\:}\or \makebox[.8cm][r]{\Huge $v$\:}\or \makebox[.8cm][r]{\Huge $v$\:}\fi}

\def\blfootnote{\xdef\@thefnmark{}\@footnotetext}

\begin{document}

\title{Outer Approximation of the Spectrum of a Fractal Laplacian}

\author{Tyrus Berry \and Steven M. Heilman \and Robert S. Strichartz$^{1}$}

\maketitle

\begin{abstract}
We present a new method to approximate the Neumann spectrum of a Laplacian on a fractal K in the plane as a renormalized limit of the Neumann spectra of the standard Laplacian on a sequence of domains that approximate K from the outside.  The method allows a numerical approximation of eigenvalues and eigenfunctions for lower portions of the spectrum.  We present experimental evidence that the method works by looking at examples where the spectrum of the fractal Laplacian is known (the unit interval and the Sierpinski Gasket (SG)).  We also present a speculative description of the spectrum on the standard Sierpinski carpet (SC), where existence of a self-similar Laplacian is known, and also on nonsymmetric and random carpets and the octagasket, where existence of a self-similar Laplacian is not known.  At present we have no explanation as to why the method should work.  Nevertheless, we are able to prove some new results about the structure of the spectrum involving ``miniaturization'' of eigenfunctions that we discovered by examining the experimental results obtained using our method.
\end{abstract}

\footnotetext[1]{The first and second authors were supported by the National Science Foundation through the Research experiences for Undergraduates (REU) Program at Cornell University.  The third author was supported in part by the National Science Foundation, Grant DMS 0652440.}

\setlength{\parindent}{6mm}

\section{Introduction}
\label{secintro}

Laplacians arise in many different mathematical contexts; three in particular that will interest us: manifolds, graphs and fractals.  There are connections relating these different types of Laplacians.  Manifold Laplacians may be obtained as limits of graph Laplacians for graphs arising from triangulations of the manifold (\cite{verdiere98,dodziuk76}).  Kigami's approach of construction Laplacians on certain fractals, such as the Sierpinski gasket (SG), also involves taking limits of graph Laplacians for graphs that approximate the fractal (\cite{kigami01,strichartz99,strichartz06}).  In this paper we present another connection, where we approximate the fractal from without by planar domains, and attempt to capture spectral information about the fractal Laplacian from spectral information about the standard Laplacian on the domains.  Thus we add an arrow to the diagram:
$$
\xymatrix
{
 & \mbox{graphs} \ar[dl] \ar[dr] &\\
\mbox{manifolds} \ar@{.>}[rr] & &\mbox{fractals}
}
$$
We should point out that the probabilistic approach to constructing Laplacians on fractals also involves approximating from without, but in that case it is the stochastic process generated by the Laplacian that is approximated, so it is not clear how to obtain spectral information.

We may describe our method succinctly as follows.  Suppose we have a self-similar fractal $K$ in the plane, determined by the identity
\begin{equation} \label{one1}
K=\bigcup F_{i}K
\end{equation}
where $\{F_{i}\}$ is a finite set of contractive similarities (called an \emph{iterated function system,} IFS).  Choose a bounded open set $\Omega$ whose closure contains $K$, and form the sequence of domains
\begin{equation} \label{one2}
\begin{split}
\Omega_{0} & =\Omega\\
\Omega_{m} & =\bigcup F_{i}\Omega_{m-1}\qquad\mbox{for }m\geq1
\end{split}
\end{equation}
Consider the standard Laplacian $\Delta$ on $\Omega_{m}$ with Neumann boundary conditions (recall that such conditions make sense even for domains with rough boundary).  Let $\{\lambda_{n}^{(m)}\}$ denote the eigenvalues in increasing order (repeated in case of nontrivial multiplicity) with eigenfunctions $\{u_{n}^{(m)}\}$ ($L^{2}$ normalized).  So
\begin{equation} \label{one3}
-\Delta u_{n}^{(m)}=\lambda_{n}^{(m)}u_{n}^{(m)}
\end{equation}
Of course $\lambda_{0}^{m}=0$ with $u_{0}^{m}$ constant.  We then hope to find a renormalization factor $r$ such that
\begin{equation} \label{one4}
\lim_{m\to\infty}r^{m}\lambda_{n}^{(m)}=\lambda_{n}
\end{equation}
exists and
\begin{equation} \label{one5}
\lim_{m\to\infty}u_{n}^{(m)}|_{K}=u_{n}
\end{equation}
exists.  (We have to be careful in cases of nontrivial multiplicity, and we may have to adjust $u_{n}^{(m)}$ by a minus sign in general).  If this is the case then we may simply define a self-adjoint operator $\Delta$ on $K$ by
\begin{equation} \label{one6}
-\Delta u_{n}=\lambda_{n}u_{n}
\end{equation}
Of course we would also like to identify $\Delta$ with a previously defined Laplacian, if such is possible, or at least show that $\Delta$ is a local operator satisfying some sort of self-similarity.

This may seem like wishful thinking, but it is not implausible.  After all, many other types of structures on fractals can be obtained as limits of structures on $\Omega_{m}$, so why not a Laplacian?  After reading this paper, we hope the reader will agree that there is a lot of evidence that this method should work in many cases. We leave to the future the challenge of describing exactly when it works, and why.

We note one great advantage of our method: it not only approximates the Laplacian, but it gives information about the spectrum.  Other methods of constructing Laplacians on fractals do not yield spectral information directly.  Of course, not all spectral information is immediately available.  In particular, asymptotic information must be lost, since we know from Weyl's law that $\lambda_{n}^{(m)}=O(n)$ for each fixed $m$, but for fractals Laplacians this is not the case.  This means, in particular, that the limit (\ref{one4}) is not uniform in $n$.  To get information about $\lambda_{n}$ for large $n$ requires taking a large value for $m$.  In practice, our numerical calculations get stuck around $m=4$.  So we only see an approximation to a segment at the bottom of the spectrum. But this is already enough to reveal aspects of the spectrum that are provable.  Briefly, if the fractal has a nontrivial finite group of symmetries, then every Neumann eigenfunction can be miniaturized, and so there is an eigenvalue renormalization factor $R$ such that if $\lambda$ is an eigenvalue then so is $R\lambda$.  The argument for this works for the approximating domains and also for a self-similar Laplacian on the fractal.  (In fact the argument could be presented on the fractal alone, so its validity is independent of the validity of the outer approximation method, but in fact it was discovered by examining the experimental data!)

So what is the evidence for the validity of the outer approximation method?  First we show that it works for the case when $K$ is the unit interval (embedded in the x-axis in the plane).  In this case we can take $F_{0}(x,y)=\left(\frac{1}{2}x,\frac{1}{2}y\right)$ and $F_{1}(x,y)=\left(\frac{1}{2}x+\frac{1}{2},\frac{1}{2}y\right)$.  If we take $\Omega$ to be the unit square, then we can compute the spectra of $\Omega_{m}$ (rectangles) and verify everything by hand ($r=1$ in this case).  We do this in section \ref{secunit}, where we also look at different choices of $\Omega$, producing sawtooth shaped domains, whose spectra are computed numerically.

In section \ref{secgasket} we look at the case of SG, where the spectrum is known exactly.  Here we see numerically how the spectra of the approximating domains approaches the known spectra.  This computation shows that the accuracy falls off rapidly as $n$ increases.  We are also able to compare the eigenfunctions of the approximating domains with the known eigenfunctions on SG.  In this case it is natural to take $\Omega$ to be a triangle containing SG in its interior since this yields connected domains $\Omega_{m}$.  We examine how the size of the overlap influences the spectra.  After the work reported in section \ref{secgasket} was completed, a different approach to outer approximation on SG was studied in \cite{blasiak08}.  In particular, different methods for choosing approximating domains are used, and a whole family of different Laplacians are studied.

In section \ref{secnonpcf} we examine numerical data for some fractals for which very little had been known about the spectrum of the Laplacian, and in some cases where even the existence of a Laplacian is unknown.  These examples fall outside of the postcritically finite (PCF) category defined in \cite{kigami01}.  The first example is the standard Sierpinski carpet SC (cut out the middle square in tic-tac-toe and iterate).  Here it is known that a self-similar Laplacian exists \cite{barlow95}, but the construction is indirect, and uniqueness is not known.  (After this work was completed, uniqueness was established in \cite{barlow08}.)  But we also examine some nonsymmetric variants of SC for which the existence of a Laplacian is unknown.  We also examine a symmetric fractal, the octagasket, where existence of a Laplacian is unknown.  In all cases the spectra of the approximating regions appear to converge when appropriately renormalized.  We can identify features of the spectrum, such as multiple eigenvalues, and eigenvalue renormalization factors $R$, and we produce rough graphs of eigenfunctions on the fractal.  In particular, there is no discernible difference between the behavior in the case of the standard SC and the other examples.

In section \ref{secmini} we describe the miniaturization process that produces the eigenvalue renormalization factor.  For this to work we need a dihedral group of symmetries of the fractal.  We only deal with the examples at hand, but it is clear that it works quite generally (we also explain how it works on the square).  For the approximating regions, this shows how $R'\lambda_{n}^{(m)}$ shows up in the spectrum on $\Omega_{m+1}$ (the factor $R'$ is not the same as $R$).

In section \ref{secranc} we examine numerical data of randomly constructed variants of SC, where the existence of Laplacians is unknown .  To make these carpets, we modify the construction of SC.  We fix the number of squares cut out at each recursive step, but we randomly determine which squares are removed.  Then, we achieve connected domains $\Omega_{m}$ with a suitable change to the above algorithm and properly chosen parameters.  Here we again see convergence of normalized eigenvalues.  These random carpets are related to the Mandelbrot percolation process.  See \cite{chayes88} and \cite{broman08}, for example.

How do we compute the spectrum of the Laplacian on the approximating domain?  We use a finite element method solver in Matlab, Matlab's own \verb!pdeeig! function.  To do this we only need to describe the geometry of the polygonal domain $\Omega_{m}$.  Then we either choose a triangulation (exclusive to Section \ref{secranc}) or let Matlab's triangulation functions \verb!decsg! and \verb!initmesh! produce a triangulation and then use piecewise linear splines in the finite element method.  Note that it would be preferable to use higher-order splines, at least piecewise cubic, since these increase accuracy dramatically for a fixed memory space and running time.  As a concession, all of our triangulations may be further refined with the \verb!refinemesh! function.  The advantage of automating the triangulation is that it saves a tremendous amount of work; in particular it chooses nonregular triangulations that increase accuracy.  The disadvantage is that the program usually does not pick a triangulation with the same symmetry as the domain.  This means that the eigenspaces that have nontrivial multiplicity in the domain end up being split into clusters of eigenspaces with eigenvalues close but not quite equal.  Since a lot of the structure of the spectrum we are trying to observe has to do with multiplicities, this forces us to make ad hoc judgements as to when we have close but unequal eigenvalues, versus multiple eigenvalues.

Why do we deal exclusively with Neumann spectra?  The main reason is that Neumann boundary conditions on the approximating domains appears to lead to Neumann boundary conditions for the Laplacian on the fractal in the case of the interval and SG, while at the same time Dirichlet boundary conditions on the approximating domains do not lead to Dirichlet boundary conditions for the Laplacian on the fractal.  For example, in the case of the interval you would need to use a mix of Dirichlet and Neumann boundary conditions on different portions of the boundary.  It is not at all clear what to do for other fractals.  Indeed for SC it is not even clear what to choose for the boundary.  The advantage of Neumann boundary conditions is that one can dispense with all notions of boundary, and define eigenfunctions simply as stationary points of the Rayleigh quotient with no boundary restrictions.  All our programs, as well as further numerical data may be found on the websites \verb!www.math.cornell.edu/~thb9d! [and \verb!www.math.cornell.edu/~smh82!].

Finally, we note that \cite{kuchment01} have studied similar outer approximations in the context of quantum graphs.

\clearpage
\section{The Unit Interval}
\label{secunit}

For the unit interval $I$ with the second derivative as Laplacian, the Neumann eigenfunctions are $\cos n\pi x$ with eigenvalues $(\pi n)^{2}$.  If we take $\Omega$ to be the unit square, then $\Omega_{m}$ is the rectangle $[0,1]\times[0,2^{-m}]$, with Neumann eigenfunctions $\cos n\pi x\cos 2^{m}k\pi y$ and eigenvalues $(\pi n)^{2}+(\pi 2^{m}k)^{2}$.  If we restrict attention to a fixed bottom segment of the spectrum, we will see eigenvalues with $k=0$ just for $m$ large enough (specifically, eigenvalues up to $L$ provided $L\leq(\pi 2^{m})^{2}$).  So $\lambda_{n}^{(m)}=\lambda_{n}$ exactly for large enough $m$.  Of course the corresponding eigenfunctions restricted to the interval give the exact eigenfunctions of the Laplacian on the interval.  Note that for each $m$ there are many other eigenfunctions on $\Omega_{m}$ (those with $k\neq0$), but they are ``blown away'' in the limit.  A similar analysis holds if we start with $\Omega$ equal to any rectangle with sides parallel to the axes.  Note that we do not have to renormalize the spectrum, or equivalently, we can take $r=1$ in (\ref{one4}).

We also note how other structures on $I$ may be approximated from corresponding structures on $\Omega_{m}$.  For example, Lebesgue measure on $I$ is the limit of Lebesgue measure on $\Omega_{m}$ suitably renormalized in the sense that
\begin{equation}\label{two1}
\lim_{m\to\infty}2^{m}\iint_{\Omega_{m}}u(x,y)dxdy=\int_{0}^{1}u(x,0)dx
\end{equation}
if $u(x,y)$ is continuous on $\Omega$ (the result is independent of the continuous extension $u(x,y)$ to $\Omega$ of $u(x,0)$ on $I$).  A similar result holds for energy, provided we use the minimum energy extension.  In other words, given $f\in H^{1}(I)$, let $u$ be the minimum energy function with $u(x,0)=f(x)$.  Then
\begin{equation}\label{two2}
\lim_{m\to\infty}2^{m}\iint_{\Omega_{m}}\abs{\nabla u(x,y)}^{2}dxdy=\int_{0}^{1}\abs{f'(x)}^{2}dx
\end{equation}
In order to see this we expand $f$ in a Fourier cosine series
\begin{equation}\label{two3}
f(x)=\sum_{k=0}^{\infty}a_{k}\cos \pi nx
\end{equation}
for which we have
\begin{equation}\label{two4}
\int_{0}^{1}\abs{f'(x)}^{2}dx=\frac{1}{2}\sum_{k=1}^{\infty}(\pi k)^{2}\abs{a_{k}}^{2}
\end{equation}
The minimum energy extension to $\Omega_{m}$ is easily seen to be
\begin{equation}\label{two5}
u(x,y)=a_{0}+\sum_{k=1}^{\infty}a_{k}\cos \pi kx\frac{\cosh 2\pi k(2^{-m}-y)}{\cosh \pi k2^{-m}}
\end{equation}
with
\begin{equation}\label{two6}
\int_{\Omega_{m}}\abs{\nabla u(x,y)}^{2}dxdy
=\sum_{k=1}^{\infty}\abs{a_{k}}^{2}\pi k\left(\frac{\sinh 2\pi k2^{-m}}{4\cosh^{2} \pi k2^{-m}}\right)
\end{equation}
Then (\ref{two2}) follows from (\ref{two4}) and (\ref{two6}).  Note that we obtain the same result if we use the simpler extension $u(x,y)=f(x)$, although this extension does not minimize energy.  (The energy minimizing extension must be harmonic on the interior and satisfy Neumann boundary conditions on the portion of the boundary of $\Omega_{m}$ disjoint from $I$, and this explains (\ref{two5})).  We also have a bilinear version: let
\begin{equation}\label{two7}
\mathcal{E}_{I}(f,g)=\int_{0}^{1}f'(x)g'(x)dx
\end{equation}
and
\begin{equation}\label{two8}
\mathcal{E}_{m}(u,v)=\int_{\Omega_{m}}(\nabla u\cdot\nabla v)dxdy
\end{equation}
If $u_{m}$ and $v_{m}$ denote the minimum energy extensions of $f$ and $g$ to $\Omega_{m}$, then
\begin{equation}\label{two9}
\lim_{m\to\infty}2^{m}\mathcal{E}_{m}(u_{m},v_{m})=\mathcal{E}_{I}(f,g)
\end{equation}
We can use this to ``define'' a Laplacian on $I$ via the weak formulation
\begin{equation}\label{two10}
\mathcal{E}_{I}(f,g)=-\int_{0}^{1}f''(x)g(x)dx
\end{equation}
if $g$ vanishes at $0$ and $1$.  By the usual Gauss-Green formula
\begin{equation}\label{two11}
\mathcal{E}_{m}(u_{m},v_{m})=\int_{\partial\Omega_{m}}(\partial_{n}u_{m})v_{m},
\end{equation}
and $\partial_{n}u_{m}=0$ on all of $\partial\Omega_{m}$ except $I$, where $\partial_{n}u_{m}=-\frac{\partial}{\partial y}u_{m}$, so
\begin{equation}\label{two12}
\mathcal{E}_{m}(u_{m},v_{m})=-\int_{0}^{1}\left(\frac{\partial u_{m}}{\partial y}\right)gdx.
\end{equation}
Combining (\ref{two9}), (\ref{two10}) and (\ref{two12}) yields at least formally
\begin{equation}\label{two13}
f''(x)=\lim_{m\to\infty}2^{m}\frac{\partial u_{m}}{\partial y}(x,0)
\end{equation}
We can verify this by differentiating (\ref{two5}) directly (assuming $f$ is smooth enough) to obtain
\begin{equation}\label{two14}
\frac{\partial u_{m}}{\partial y}(x,0)
=-\sum_{k=1}^{\infty}(\pi k)^{2}a_{k}\cos \pi kx\left(\frac{\pi k\sinh 2\pi k2^{-m}}{\cosh \pi k2^{-m}}\right)
\end{equation}
and taking the limit to obtain
\begin{equation}\label{two15}
\lim_{m\to\infty}2^{m}\frac{\partial u_{m}}{\partial y}(x,0)=-\sum_{k=1}^{\infty}(\pi k)^{2}a_{k}\cos \pi kx,
\end{equation}
and this is the same value for $f''(x)$ that we obtain by differentiating (\ref{two3}) directly.

For a less trivial example we need only to take a geometrically more interesting $\Omega$.  In particular, let $\Omega$ be a triangle with vertices $(-\epsilon,0),$ $(1+\epsilon,0)$ and $\left(\frac{1}{2},h\right)$ for some choice of positive parameters $\epsilon$ and $h$.  Then $\Omega_{m}$ is a sawtooth region with $2^{m}$ teeth, maximum height $2^{-m}h$ and overlaps of length $2^{-m}\epsilon$.  It is not feasible to compute the Neumann spectrum of the $\Omega_{m}$ exactly, so we use numerical methods.  In Tables \ref{tabletwo1} and \ref{tabletwo3} we present the eigenvalues for several choices of parameters and level $m=2,3,4$ (we also vary the number of refinements used in the FEM approximation).  Actually the computations are done for a similar image of $\Omega_{m}$ so that the base is exactly $I$, but this makes no difference in the limit.  The evidence suggests that we get $c(\epsilon,h)\frac{d^{2}}{dx^{2}}$ in the limit for some constant that depends on the parameters.

In Tables \ref{tabletwo2} and \ref{tabletwo4} we present the same data, but we normalize by dividing $\lambda_{n}^{(m)}$ by $\lambda_{1}^{(m)}$.  This enables us to compare the normalized eigenvalues with the expected values $n^{2}$.  Note that with level $m=5$ we see about a $1\%$ deviation already at $n=6$.

\begin{figure}[htbp!]
  \includegraphics[width=\threewidth\textwidth]{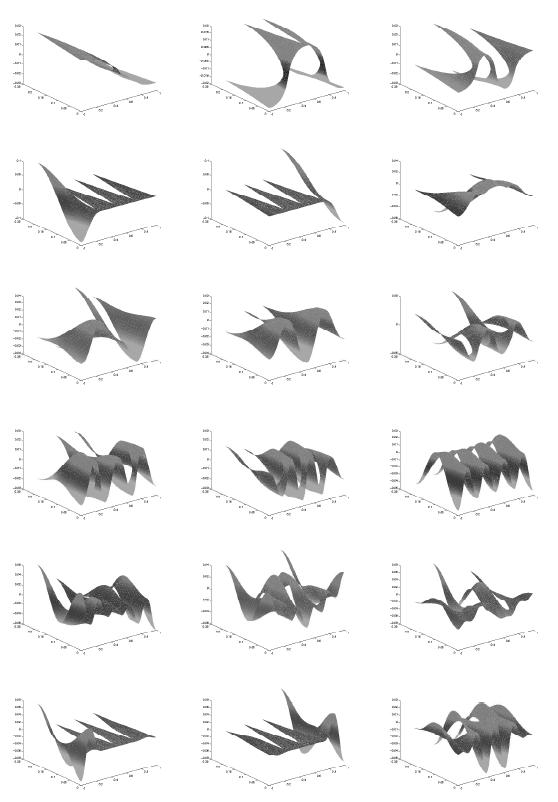}
\caption{Sawtooth Eigenfunctions, m=2}
\label{figtwo1}
\end{figure}

\begin{figure}
\begin{center}
\includegraphics[width=\onewidth\textwidth]{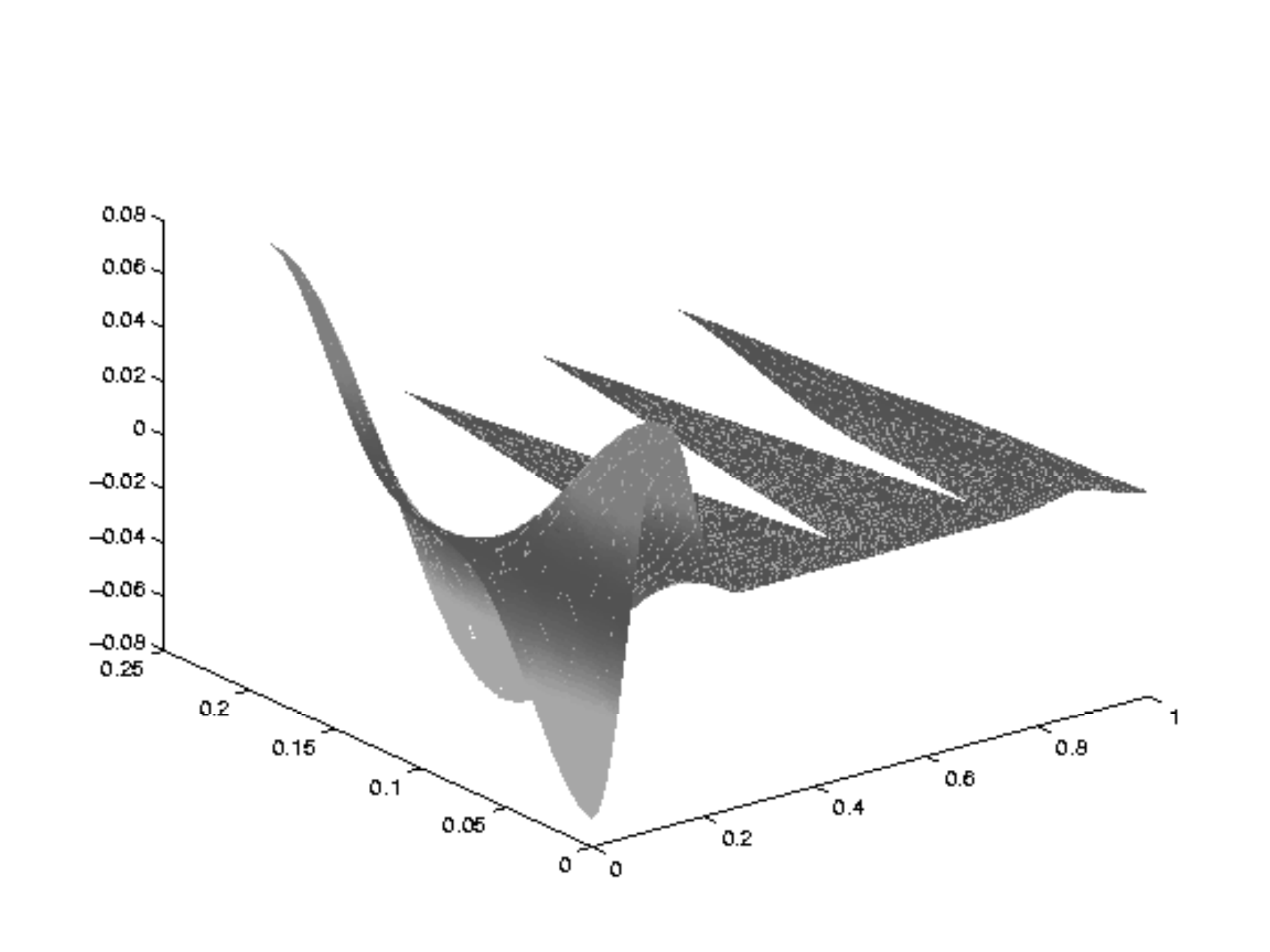}
\end{center}
\caption{Almost Localized Sawtooth Eigenfunction, m=2}
\label{figtwo2}
\end{figure}

In Figure \ref{figtwo1} we show some graphs of eigenfunctions on $\Omega_{m}$, that approximate eigenfunctions on $I$.  In Figure \ref{figtwo2} we show the graph of an eigenfunction on $\Omega_{2}$ that does not approximate an eigenfunction on $I$.  Indeed, this eigenfunction appears to be almost localized to one of the teeth.  This phenomenon is discussed in \cite{heilman10}.  Unfortunately, we do not know if we can define energy on $I$ via (\ref{two2}) for a sawtooth region approximation.  Indeed, we have no idea what the minimum energy extension looks like.

\begin{table}[htbp!]
  \resizebox{\tabonewd}{!}{
  \begin{tabular}{r | r*{21}{r}r}
    \multicolumn{13}{p{21cm}}{Equilateral Triangles Sawtooth Region (height determined by requirement that triangles are equilateral, overlaps set to $(2^{-m})/10)$}\\
    \hline
    Level: & 2 & 2 & 2 & 2 & 3 & 3 & 3 & 3 & 4 & 4 & 4 & 4 \\
    Refinement: & 1 & 2 & 3 & 4 & 1 & 2 & 3 & 4 & 1 & 2 & 3 & 4 \\
    \hline
    n\\	
    1 & 4.905	  & 4.823  &  4.790 & 4.777  & 4.868 & 4.789   & 4.756  & 4.743  & 4.828  & 4.749  & 4.717 & 4.703\\
    2 & 17.980  & 17.662 & 17.535 & 17.483 & 19.097 & 18.782 & 18.655	& 18.602 & 19.218 & 18.903 & 18.776	& 18.724\\
    3 & 33.418  & 32.790 & 32.53  & 32.436 & 41.513 & 40.808 & 40.525	& 40.408 & 42.900 & 42.191 & 41.906	& 41.787\\
    4 & 246.809 & 243.909	& 243.176 & 242.991	& 69.950 & 68.724 & 68.232 & 68.029 & 75.398 & 74.140 & 73.635	& 73.425\\
    5 & 246.809 & 243.910	& 243.176 & 242.991	& 100.984 & 99.165 & 98.436	& 98.134 & 116.012	& 114.066	& 113.283	 & 112.958\\
    6 & 248.850 & 246.991	& 246.524 & 246.407	& 129.818 & 127.424 & 126.463 & 126.064 & 163.836 & 161.053 & 159.935 & 159.471\\
    7 & 250.833 & 248.743	& 248.218 & 248.087	& 150.737 & 147.863	& 146.713 & 146.238 & 217.674	& 213.915	& 212.408	& 211.783\\
    8 & 253.564 & 251.508	& 250.992 & 250.863	& 959.592 & 952.139	& 950.177 & 949.677	& 276.002	& 271.161	& 269.220	& 268.417\\
    9 & 337.235 & 332.179	& 330.654 & 330.157	& 959.592 & 952.139	& 950.177 & 949.677	& 336.966	& 330.967	& 328.563	& 327.569\\
    10 & 389.324 & 382.371	& 380.228 & 379.513	& 970.250 & 963.501	& 961.797 & 961.369	& 398.337	& 391.162	& 388.285	& 387.094\\
    11 & 449.038 & 440.249	& 437.449 & 436.483	& 971.305 & 964.578	& 962.895 & 962.474	& 457.640	& 449.268	& 445.913	& 444.524\\
    12 & 782.622 & 760.213	& 754.319 & 752.824	& 973.587 & 966.696	& 964.959 & 964.524	& 512.105	& 502.535	& 498.708	& 497.126\\
    13 & 817.310 & 787.079	& 779.191 & 777.161	& 977.148 & 970.044	& 968.246 & 967.794	& 558.576	& 548.015	& 543.791	& 542.045\\
    14 & 884.992 & 851.631	& 843.276 & 841.121	& 982.108 & 974.799	& 972.955 & 972.492	& 594.214	& 582.929	& 578.409	& 576.538\\
    15 & 931.268 & 900.120	& 892.271 & 890.247	& 987.564 & 980.415	& 978.615 & 978.163	& 616.898	& 605.004	& 600.254	& 598.291\\
    16 & 1022.576 & 984.632 & 974.972 & 972.543 & 991.919 & 985.129	& 983.418 & 982.989	\\			
    17 & 1022.620 & 984.636 & 974.973 & 972.544 & 1253.731 & 1237.690	& 1232.817 & 1231.219\\			
    18 & 1023.447 & 995.768 & 988.646	& 986.841 & 1315.232 & 1297.342	& 1291.832 & 1290.003\\			
    19 & 1042.893 & 1004.686 & 995.410	& 993.103 & 1407.013 & 1386.000	& 1379.452 & 1377.254\\			
    20 & 1050.656 & 1012.066 & 1002.876 & 1000.611	& 1518.109 & 1493.410 & 1485.623 & 1482.974\\		
    21 & 1269.022 & 1219.852 & 1206.751 & 1203.162	& 1636.329 & 1608.281 & 1599.289 & 1596.175\\		
    22 & 1288.306 & 1235.661 & 1221.617 & 1217.694	& 1747.569 & 1716.437 & 1706.278 & 1702.698\\		
    23 & 1304.295 & 1255.860 & 1242.652 & 1238.793	& 1833.745 & 1798.417 & 1786.867 & 1782.782\\		
    24 & 1855.667 & 1749.382 & 1712.980 & 1703.675\\				
    25 & 1878.164 & 1749.398 & 1712.983 & 1703.676\\							
    26 & 1883.733 & 1766.473 & 1743.511 & 1737.727\\							
    27 & 1883.965 & 1785.448 & 1761.326 & 1755.233\\							
    28 & 1909.401 & 1821.564 & 1798.312 & 1792.391\\							
	29 &          & 1948.299 & 1915.997 & 1907.606					
    \end{tabular}
    }
    \caption{Sawtooth Unnormalized Eigenvalues, built with Equilateral Triangles}
    \label{tabletwo1}
\end{table}

\begin{table}[htbp!]
  \resizebox{\tabonewd}{!}{
  \begin{tabular}{r | r*{12}{r}r}
    \multicolumn{13}{p{21cm}}{Equilateral Triangles Sawtooth Region (height determined by requirement that triangles are equilateral, overlaps set to $(2^{-m})/10)$}\\
    \hline
    Level: & 2 & 2 & 2 & 2 & 3 & 3 & 3 & 3 & 4 & 4 & 4 & 4 \\
    Refinement: & 1 & 2 & 3 & 4 & 1 & 2 & 3 & 4 & 1 & 2 & 3 & 4 \\
    \hline
    n\\	
    1 & 1.000 & 1.000 & 1.000 & 1.000 & 1.000 & 1.000 & 1.000 & 1.000 & 1.000 & 1.000 & 1.000 & 1.000\\
    2 & 3.665 & 3.662 & 3.661 & 3.660 & 3.923 & 3.922 & 3.922 & 3.922 & 3.981 & 3.981 & 3.981 & 3.981\\
    3 & 6.813 & 6.798 & 6.793 & 6.790 & 8.527 & 8.522 & 8.520 & 8.519 & 8.886 & 8.885 & 8.884 & 8.884\\
    4 & 50.316 & 50.570 & 50.764 & 50.870 & 14.368 & 14.352 & 14.345 & 14.343 & 15.618 & 15.613 & 15.612	& 15.611\\
    5 & 50.316 & 50.570 & 50.764 & 50.870 & 20.742 & 20.709 & 20.695 & 20.690 & 24.030 & 24.021 & 24.017	& 24.016\\
    6 & 50.732 & 51.209 & 51.463 & 51.585 & 26.665 & 26.610 & 26.588 & 26.578 & 33.936 & 33.916 & 33.908	& 33.905\\
    7 & 51.136 & 51.572 & 51.816 & 51.936 & 30.962 & 30.878 & 30.845 & 30.832 & 45.088 & 45.048 & 45.033	& 45.027\\
    8 & 51.693 & 52.145 & 52.395 & 52.518 & 197.104 & 198.834 & 199.766 & 200.222 & 57.170 & 57.104 & 57.078	& 57.068\\
    9 & 68.750 & 68.871 & 69.025 & 69.118 & 197.104 & 198.834 & 199.766 & 200.222 & 69.797 & 69.698 & 69.659	& 69.644\\
    10 & 79.370 & 79.277 & 79.374 & 79.450 & 199.293 & 201.207 & 202.209 & 202.687 & 82.510 & 82.375 & 82.321	& 82.299\\
    11 & 91.543 & 91.277 & 91.319 & 91.377 & 199.510 & 201.432 & 202.440 & 202.920 & 94.793 & 94.611 & 94.539	& 94.509\\
    12 & 159.549 & 157.615 & 157.466 & 157.602 & 199.978 & 201.874 & 202.874 & 203.352 & 106.075 & 105.829	& 105.732 & 105.693\\
    13 & 166.621 & 163.185 & 162.658 & 162.697 & 200.710 & 202.573 & 203.564 & 204.042 & 115.701 & 115.406	& 115.290 & 115.243\\
    14 & 180.419 & 176.568 & 176.036 & 176.087 & 201.729 & 203.566 & 204.554 & 205.032 & 123.083 & 122.759	& 122.630 & 122.576\\
    15 & 189.853 & 186.622 & 186.264 & 186.372 & 202.849 & 204.739 & 205.745 & 206.228 & 127.781 & 127.408	& 127.261 & 127.201\\
    16 & 208.468 & 204.144 & 203.528 & 203.600 & 203.744 & 205.723 & 206.754 & 207.245\\
    17 & 208.477 & 204.144 & 203.528 & 203.600 & 257.521 & 258.466 & 259.188 & 259.580\\
    18 & 208.645 & 206.452 & 206.383 & 206.593 & 270.153 & 270.923 & 271.595 & 271.974\\
    19 & 212.610 & 208.301 & 207.795 & 207.904 & 289.006 & 289.437 & 290.017 & 290.369\\
    20 & 214.192 & 209.831 & 209.353 & 209.476 & 311.825 & 311.867 & 312.338 & 312.658\\
    21 & 258.709 & 252.912 & 251.913 & 251.880 & 336.108 & 335.856 & 336.235 & 336.525\\
    22 & 262.641 & 256.189 & 255.016 & 254.922 & 358.957 & 358.442 & 358.729 & 358.983\\
    23 & 265.900 & 260.377 & 259.407 & 259.339 & 376.658 & 375.562 & 375.672 & 375.867\\
    24 & 378.306 & 362.699 & 357.590 & 356.661\\
    25 & 382.892 & 362.702 & 357.590 & 356.661\\
    26 & 384.028 & 366.242 & 363.963 & 363.790\\
    27 & 384.075 & 370.176 & 367.682 & 367.455\\
    28 & 389.260 & 377.664 & 375.403 & 375.234\\
	29 &         & 403.940 & 399.970 & 399.354\\
    \end{tabular}
    }
    \caption{Sawtooth Normalized Eigenvalues, built with Equilateral Triangles}
    \label{tabletwo2}
\end{table}

\begin{sidewaystable}[htbp!]
\centering
\resizebox{\tabtwowd}{!}{
  \begin{tabular}{r | *{20}{r}}
    \multicolumn{11}{p{10cm}}{Eigenvalue Data for Sawtooth Regions with different Parameters}\\
    \hline
    Level: & 2 & 2 & 2 & 2 & 3 & 3 & 3 & 3 & 4 & 4 & 4 & 4 & 5 & 5 & 5 & 4 & 4 & 4 & 4 & 5 \\
   Height: & 0.100 &   0.010 &   0.001 &   0.0005 &   0.100 &   0.010 &   0.001 &   0.0005 &   0.100 &   0.010
            & 0.001 &   0.0005 &   0.100 &   0.010 &   0.001 &   0.100 &   0.010 &   0.001 &   0.0005 &   0.001\\
    Refinement: & 2 & 2 & 2 & 2 & 2 & 2 & 2 & 2 & 2 & 2 & 2 & 2 & 2 & 2 & 2 & 3 & 3 & 3 & 3 & 3 \\
    \hline
    n\\	
  1 &   6.256 &   7.155 &   7.550 &   8.017 &   5.093 &   6.996 &   7.229 &   7.326 &   3.836 &   6.782 &   7.060 &   7.155 &   2.867 &   6.341 &   6.981 &   3.795 &   6.768 &   7.011 &   7.061 &   6.931\\
  2 &  23.467 &  26.916 &  29.096 &  31.996 &  20.023 &  27.616 &  28.606 &  29.123 &  15.210 &  27.041 &  28.154 &  28.543 &  11.388 &  25.345 &  27.902 &  15.046 &  26.987 &  27.957 &  28.160 &  27.704\\
  3 &  45.256 &  52.154 &  56.017 &  60.827 &  43.675 &  60.656 &  63.040 &  64.508 &  33.701 &  60.510 &  63.022 &  63.921 &  25.325 &  56.948 &  62.698 &  33.341 &  60.387 &  62.567 &  63.035 &  62.251\\
  4 & 294.835 & 327.674 & 331.873 & 333.511 &  73.960 & 103.746 & 108.178 & 111.071 &  58.625 & 106.725 & 111.207 & 112.857 &  44.286 & 101.047 & 111.258 &  58.004 & 106.506 & 110.372 & 111.227 & 110.461\\
  5 & 319.781 & 356.983 & 368.346 & 377.480 & 107.475 & 152.707 & 159.642 & 163.925 &  89.037 & 164.994 & 172.017 & 174.690 &  67.757 & 157.493 & 173.423 &  88.094 & 164.649 & 170.667 & 172.045 & 172.175\\
 6 & 388.559 & 435.688 & 465.781 & 504.639 & 139.164 & 200.594 & 210.056 & 214.981 & 123.691 & 234.348 & 244.478 & 248.465 &  95.054 & 226.090 & 248.988 & 122.397 & 233.848 & 242.467 & 244.511 & 247.185\\
 7 & 478.443 & 535.783 & 582.431 & 651.190 & 162.497 & 237.152 & 248.566 & 253.227 & 161.172 & 313.478 & 327.260 & 332.860 & 125.431 & 306.595 & 337.691 & 159.508 & 312.794 & 324.432 & 327.290 & 335.229\\
 8 & 696.381 & 851.228 & 1007.893 & 1106.376 & 1017.767 & 1277.229 & 1293.087 & 1295.630 & 199.921 & 400.647 & 418.579 & 426.082 & 158.038 & 398.709 & 439.215 & 197.867 & 399.751 & 414.783 & 418.596 & 435.987\\
 9 & 766.306 & 930.457 & 1254.199 & 1439.666 & 1027.972 & 1304.205 & 1326.104 & 1331.648 & 238.235 & 493.580 & 516.082 & 525.745 & 192.104 & 502.074 & 553.181 & 235.820 & 492.443 & 511.178 & 516.074 & 549.076\\
 10 & 890.920 & 1092.418 & 1677.158 &  & 1068.070 & 1393.716 & 1422.250 & 1432.226 & 274.559 & 589.334 & 616.709 & 628.727 & 226.653 & 616.269 & 679.138 & 271.805 & 587.934 & 610.581 & 616.660 & 674.045\\
 11 & 1003.561 & 1266.292 & 1819.137 &  & 1121.065 & 1527.962 & 1568.308 & 1587.383 & 307.455 & 684.176 & 716.553 & 731.030 & 260.952 & 740.802 & 816.560 & 304.396 & 682.499 & 709.130 & 716.448 & 810.365\\
 12 & 1557.709 & 1812.385 &  &  & 1175.215 & 1694.892 & 1751.443 & 1784.010 & 335.656 & 773.516 & 810.781 & 827.688 & 294.369 & 875.099 & 964.834 & 332.378 & 771.560 & 802.054 & 810.605 & 957.419\\
 13 & 1561.255 & 1931.647 &  &  & 1220.465 & 1882.831 & 1958.367 &  & 358.389 & 851.985 & 893.700 & 912.842 & 326.301 & 1018.500 & 1123.252 & 354.918 & 849.767 & 883.758 & 893.449 & 1114.496\\
 14 & 1608.660 &  &  &  & 1252.193 &  &  &  & 375.031 & 913.808 & 959.146 & 980.122 & 356.303 & 1170.244 & 1290.993 & 371.402 & 911.372 & 948.192 & 958.826 & 1280.776\\
 15 & 1616.961 &  &  &  & 1270.474 &  &  &  & 385.025 & 953.571 & 1001.298 & 1023.492 & 384.206 & 1329.457 & 1467.118 & 381.364 & 950.989 & 989.668 & 1000.930 & 1455.315\\
 16 & 1618.303 &  &  &  & 1277.740 &  &  &  & 1427.143 &  &  &  & 409.710 & 1495.142 & 1650.547 & 1425.511 &  &  &  & 1637.034\\
 17 & 1865.746 &  &  &  & 1592.815 &  &  &  & 1429.230 &  &  &  & 432.971 & 1666.162 & 1840.045 & 1427.493 &  &  &  & 1824.701\\
 18 & 1996.687 &  &  &  & 1610.143 &  &  &  & 1437.825 &  &  &  & 453.826 & 1841.225 &  & 1436.042 &  &  &  & \\
 19 &  &  &  &  & 1641.624 &  &  &  & 1452.086 &  &  &  & 472.532 &  &  & 1450.237 &  &  &  & \\
 20 &  &  &  &  & 1689.362 &  &  &  & 1472.750 &  &  &  & 489.103 &  &  & 1470.701 &  &  &  & \\
 21 &  &  &  &  & 1751.504 &  &  &  & 1500.282 &  &  &  & 503.710 &  &  & 1497.967 &  &  &  & \\
 22 &  &  &  &  & 1818.223 &  &  &  & 1535.318 &  &  &  & 516.571 &  &  & 1532.667 &  &  &  & \\
 23 &  &  &  &  & 1871.686 &  &  &  & 1578.626 &  &  &  & 527.762 &  &  & 1575.480 &  &  &  & \\
 24 &  &  &  &  &  &  &  &  & 1630.678 &  &  &  & 537.315 &  &  & 1626.922 &  &  &  & \\
 25 &  &  &  &  &  &  &  &  & 1691.729 &  &  &  & 545.583 &  &  & 1687.201 &  &  &  & \\
 26 &  &  &  &  &  &  &  &  & 1761.543 &  &  &  & 552.515 &  &  & 1755.965 &  &  &  & \\
 27 &  &  &  &  &  &  &  &  & 1838.629 &  &  &  & 558.237 &  &  & 1831.784 &  &  &  & \\
 28 &  &  &  &  &  &  &  &  & 1919.939 &  &  &  & 562.760 &  &  & 1911.655 &  &  &  & \\
 29 &  &  &  &  &  &  &  &  &  &  &  &  & 566.286 &  &  & 1990.465 &  &  &  & \\
 30 &  &  &  &  &  &  &  &  &  &  &  &  & 568.765 &  &  &  &  &  &  & \\
 31 &  &  &  &  &  &  &  &  &  &  &  &  & 570.211 &  &  &  &  &  &  & \\
 32 &  &  &  &  &  &  &  &  &  &  &  &  & 1466.956 &  &  &  &  &  &  & \\
 33 &  &  &  &  &  &  &  &  &  &  &  &  & 1469.178 &  &  &  &  &  &  & \\
     \end{tabular}
     }
    \caption{Sawtooth Unnormalized Eigenvalues}
    \label{tabletwo3}
\end{sidewaystable}

\begin{sidewaystable}[htbp!]
\centering
\resizebox{\tabtwowd}{!}{
  \begin{tabular}{r | *{20}{r}}
    \multicolumn{11}{p{10cm}}{Eigenvalue Data for Sawtooth Regions with different Parameters}\\
    \hline
    Level: & 2 & 2 & 2 & 2 & 3 & 3 & 3 & 3 & 4 & 4 & 4 & 4 & 5 & 5 & 5 & 4 & 4 & 4 & 4 & 5 \\
    Height: & 0.100 &   0.010 &   0.001 &   0.0005 &   0.100 &   0.010 &   0.001 &   0.0005 &   0.100 &   0.010
            & 0.001 &   0.0005 &   0.100 &   0.010 &   0.001 &   0.100 &   0.010 &   0.001 &   0.0005 &   0.001\\
    Refinement: & 2 & 2 & 2 & 2 & 2 & 2 & 2 & 2 & 2 & 2 & 2 & 2 & 2 & 2 & 2 & 3 & 3 & 3 & 3 & 3 \\
    \hline
    n\\	
  1 &   1.000 &   1.000 &   1.000 &   1.000 &   1.000 &   1.000 &   1.000 &   1.000 &   1.000 &   1.000 &   1.000 &   1.000 &   1.000 &   1.000 &   1.000 &   1.000 &   1.000 &   1.000 &   1.000 &   1.000\\
  2 &   3.751 &   3.762 &   3.854 &   3.991 &   3.931 &   3.947 &   3.957 &   3.975 &   3.965 &   3.987 &   3.988 &   3.989 &   3.973 &   3.997 &   3.997 &   3.965 &   3.987 &   3.987 &   3.988 &   3.997\\
  3 &   7.234 &   7.289 &   7.419 &   7.587 &   8.575 &   8.670 &   8.721 &   8.805 &   8.786 &   8.922 &   8.927 &   8.933 &   8.834 &   8.980 &   8.981 &   8.786 &   8.922 &   8.924 &   8.927 &   8.981\\
  4 &  47.131 &  45.799 &  43.955 &  41.602 &  14.521 &  14.829 &  14.965 &  15.160 &  15.283 &  15.737 &  15.752 &  15.773 &  15.448 &  15.935 &  15.937 &  15.285 &  15.736 &  15.742 &  15.752 &  15.936\\
  5 &  51.119 &  49.895 &  48.786 &  47.086 &  21.102 &  21.827 &  22.084 &  22.375 &  23.211 &  24.328 &  24.366 &  24.414 &  23.636 &  24.836 &  24.842 &  23.214 &  24.326 &  24.342 &  24.365 &  24.840\\
 6 &  62.113 &  60.896 &  61.690 &  62.948 &  27.323 &  28.671 &  29.059 &  29.343 &  32.245 &  34.554 &  34.630 &  34.725 &  33.158 &  35.654 &  35.666 &  32.253 &  34.550 &  34.583 &  34.627 &  35.661\\
 7 &  76.481 &  74.886 &  77.140 &  81.228 &  31.905 &  33.897 &  34.386 &  34.564 &  42.015 &  46.222 &  46.356 &  46.520 &  43.754 &  48.349 &  48.373 &  42.032 &  46.214 &  46.274 &  46.350 &  48.364\\
 8 & 111.320 & 118.976 & 133.491 & 138.007 & 199.828 & 182.557 & 178.882 & 176.844 &  52.117 &  59.075 &  59.291 &  59.548 &  55.128 &  62.875 &  62.915 &  52.140 &  59.061 &  59.161 &  59.281 &  62.900\\
 9 & 122.498 & 130.049 & 166.113 & 179.581 & 201.832 & 186.413 & 183.449 & 181.760 &  62.105 &  72.778 &  73.102 &  73.477 &  67.011 &  79.175 &  79.240 &  62.140 &  72.756 &  72.909 &  73.085 &  79.215\\
 10 & 142.418 & 152.687 & 222.132 &  & 209.704 & 199.207 & 196.750 & 195.488 &  71.574 &  86.896 &  87.356 &  87.870 &  79.063 &  97.184 &  97.283 &  71.623 &  86.864 &  87.087 &  87.330 &  97.245\\
 11 & 160.424 & 176.989 & 240.936 &  & 220.109 & 218.395 & 216.955 & 216.666 &  80.149 & 100.881 & 101.499 & 102.167 &  91.027 & 116.822 & 116.968 &  80.211 & 100.836 & 101.143 & 101.461 & 116.912\\
 12 & 249.007 & 253.316 &  &  & 230.741 & 242.255 & 242.290 & 243.504 &  87.501 & 114.054 & 114.846 & 115.676 & 102.684 & 138.000 & 138.208 &  87.585 & 113.994 & 114.397 & 114.796 & 138.127\\
 13 & 249.574 & 269.985 &  &  & 239.626 & 269.118 & 270.915 &  &  93.427 & 125.624 & 126.591 & 127.577 & 113.823 & 160.614 & 160.900 &  93.524 & 125.548 & 126.050 & 126.528 & 160.789\\
 14 & 257.152 &  &  &  & 245.855 &  &  &  &  97.765 & 134.740 & 135.862 & 136.980 & 124.288 & 184.544 & 184.929 &  97.868 & 134.650 & 135.241 & 135.786 & 184.778\\
 15 & 258.479 &  &  &  & 249.444 &  &  &  & 100.371 & 140.603 & 141.833 & 143.041 & 134.022 & 209.651 & 210.158 & 100.493 & 140.503 & 141.156 & 141.749 & 209.959\\
 16 & 258.694 &  &  &  & 250.871 &  &  &  & 372.037 &  &  &  & 142.918 & 235.779 & 236.433 & 375.635 &  &  &  & 236.175\\
 17 & 298.248 &  &  &  & 312.733 &  &  &  & 372.581 &  &  &  & 151.032 & 262.748 & 263.578 & 376.157 &  &  &  & 263.250\\
 18 & 319.180 &  &  &  & 316.135 &  &  &  & 374.821 &  &  &  & 158.307 & 290.355 &  & 378.409 &  &  &  & \\
 19 &  &  &  &  & 322.316 &  &  &  & 378.539 &  &  &  & 164.832 &  &  & 382.150 &  &  &  & \\
 20 &  &  &  &  & 331.689 &  &  &  & 383.926 &  &  &  & 170.613 &  &  & 387.542 &  &  &  & \\
 21 &  &  &  &  & 343.890 &  &  &  & 391.103 &  &  &  & 175.708 &  &  & 394.727 &  &  &  & \\
 22 &  &  &  &  & 356.989 &  &  &  & 400.236 &  &  &  & 180.194 &  &  & 403.871 &  &  &  & \\
 23 &  &  &  &  & 367.486 &  &  &  & 411.526 &  &  &  & 184.098 &  &  & 415.153 &  &  &  & \\
 24 &  &  &  &  &  &  &  &  & 425.095 &  &  &  & 187.430 &  &  & 428.708 &  &  &  & \\
 25 &  &  &  &  &  &  &  &  & 441.010 &  &  &  & 190.314 &  &  & 444.592 &  &  &  & \\
 26 &  &  &  &  &  &  &  &  & 459.210 &  &  &  & 192.732 &  &  & 462.712 &  &  &  & \\
 27 &  &  &  &  &  &  &  &  & 479.305 &  &  &  & 194.728 &  &  & 482.691 &  &  &  & \\
 28 &  &  &  &  &  &  &  &  & 500.502 &  &  &  & 196.306 &  &  & 503.738 &  &  &  & \\
 29 &  &  &  &  &  &  &  &  &  &  &  &  & 197.536 &  &  & 524.505 &  &  &  & \\
 30 &  &  &  &  &  &  &  &  &  &  &  &  & 198.401 &  &  &  &  &  &  & \\
 31 &  &  &  &  &  &  &  &  &  &  &  &  & 198.905 &  &  &  &  &  &  & \\
 32 &  &  &  &  &  &  &  &  &  &  &  &  & 511.714 &  &  &  &  &  &  & \\
 33 &  &  &  &  &  &  &  &  &  &  &  &  & 512.490 &  &  &  &  &  &  & \\
     \end{tabular}
     }
    \caption{Sawtooth Normalized Eigenvalues}
    \label{tabletwo4}
\end{sidewaystable}

There is yet another outer approximation approach to $I$, in which we regard it as the bottom line in SG.  So we take $\Omega=\mbox{SG}$ and $\Omega_{m+1}=F_{1}\Omega_{m}\bigcup F_{2}\Omega_{m}$.  Then $\Omega_{m}$ is a fractafold in the sense of \cite{strichartz03} consisting of $2^{m}$ cells of level $m$ along the bottom of SG.  The bottom $2^{m}$ Neumann eigenfunctions of the fractal Laplacian on $\Omega_{m}$ are obtained by the method of spectral decimation as follows.  Fix a parameter $j$ satisfying $0\leq j<2^{m}$.  Let $x_{k}=\frac{k}{2^{m}}$ for $0\leq k\leq 2^{m}$ denote the points along $I$ where the cells of $\Omega_{m}$ intersect, and let $y_{k}$ for $1\leq k\leq 2^{m}$ denote the top vertices of the cells (so cell number $k$ has vertices $x_{k-1},x_{k},y_{k}$). Then $u_{j}$ restricted to these points is defined by
\begin{equation}\label{two16}
\left\{
\begin{aligned}
u_{j}(x_{k})& =\frac{1}{2}(\cos \pi jx_{k}+\cos \pi jx_{k+1})\\
u_{j}(y_{k})& =\cos\pi jx_{k}
\end{aligned}
\right.
\end{equation}
One can check that for a graph Laplacian $\Delta_{m}$ on the graph $\{x_{k},y_{k}\}$ we have
\begin{equation}\label{two17}
-\Delta_{m}u_{y}=\left(2-2\cos\frac{\pi j}{2^{m}}\right)u_{j}
\end{equation}
with the appropriate Neumann conditions at the boundary points $x_{0},x_{2m}$.  Let
\begin{equation}\label{two18}
\phi_{-}(t)=\frac{5-\sqrt{25-4t}}{2}
\end{equation}
and
\begin{equation}\label{two19}
\Phi(t)=\lim_{n\to\infty}5^{n}\phi_{-}^{(n)}(t),
\end{equation}
where $\phi_{-}^{(n)}(t)$ denotes the n-fold composition.  In particular, $\Phi$ is a smooth function with $\Phi(0)=0$ and $\Phi'(0)=1$.  Then the method of spectral decimation (See \cite{strichartz06} for a detailed explanation) says that $u_{j}$ may be extended to eigenfunctions of the fractal Laplacian on $\Omega_{m}$ with eigenvalue
\begin{equation}\label{two20}
\lambda_{j}^{(m)}=\frac{3}{2}\lim_{n\to\infty}5^{m+n}\phi_{-}^{(n)}\left(2-2\cos\frac{\pi j}{2^{m}}\right)
=\frac{3}{2}5^{m}\Phi\left(2-2\cos\frac{\pi j}{2^{m}}\right)
\end{equation}
Now observe that $2-2\cos\frac{\pi j}{2^{m}}\approx(\frac{\pi j}{2^{m}})^{2}$ for large $m$ so
\begin{equation}\label{two21}
\lim_{m\to\infty}\left(\frac{4}{5}\right)^{m}\lambda_{j}^{(m)}=\frac{3}{2}(\pi j)^{2}.
\end{equation}
Of course $(\pi j)^{2}$ is the correct eigenvalue for the eigenfunction $\cos \pi jx$ on $I$, which is clearly the limit of $u_{j}$ as $u\to\infty$.

\clearpage
\section{The Sierpinski Gasket}
\label{secgasket}

Let $\{q_{0},q_{1},q_{2}\}$ denote the vertices of a unit length equilateral triangle in the plane, and let $F_{i}x=\frac{1}{2}(x+q_{i})$ for $i=0,1,2$.  Then SG is the invariant set for this IFS.  We take $\Omega$ to be the equilateral triangle dilated by a factor $1+\epsilon$.  Then $\Omega_{m}$ is a union of $3^{m}$ triangles of size $2^{-m}$ that overlap in triangles of size $(1+\epsilon)2^{-m}$.

\begin{figure}
\includegraphics[width=\threewidth\textwidth]{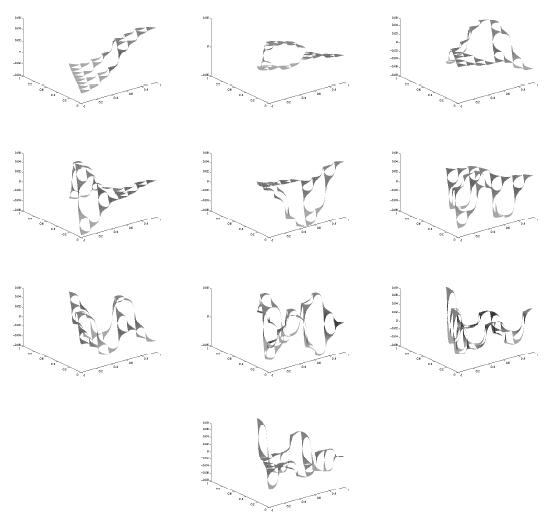}
\caption{Sierpinski Gasket (SG) Eigenfunctions, Level 3}
\label{figthree1}
\end{figure}

In Tables \ref{tablethree1} and \ref{tablethree2} we present the same data as in Tables \ref{tabletwo1} through \ref{tabletwo4} for this example.  The multiplicities and normalized eigenvalues agree with the known values for the Neumann spectrum of the standard Laplacian on SG \cite{strichartz03}.  For example, the first six distinct normalized eigenvalues on SG are $1,5,8.103,10.305,25,31.784$.  So the numerical accuracy improves as we decrease $\epsilon$ but the error remains significant.  (Much better accuracy is achieved in \cite{blasiak08}).  Nevertheless, the qualitative features of the spectrum, including high multiplicities and large gaps, are already apparent.  In Figure \ref{figthree1} we show some graphs of eigenfunctions.  Actual graphs of Dirichlet eigenfunctions on SG may be found in \cite{dalrymple99}.

In this case we know the eigenfunction renormalization factor $R=5$, so we expect $r=1.25$ in (\ref{one4}).  The data is not inconsistent with this expectation, but it is impossible to deduce these values from the data alone.

\begin{table}[htbp!]
  \resizebox{\tabthreewd}{!}{
  \begin{tabular}{r | r*{11}{r}r}
    \multicolumn{7}{p{9cm}}{Sierpinski Gasket Eigenvalue Data}\\
    \hline
    Level: & 2 & 2 & 2 & 2 & 2 & 3 & 3 & 3 & 3 & 4 & 5 \\
    Refinement: & 0 & 1 & 2 & 3 & 4 & 0 & 1 & 2 & 3 & 0 & 0 \\
    \hline
    n\\
1 & 5.0727 & 4.8920 & 4.8223 & 4.7946 & 4.7832 & 4.1689 & 4.0255 & 3.9697 & 3.9473 & 3.3372 & 2.6327\\
2 & 5.0729 & 4.8924 & 4.8226 & 4.7948 & 4.7833 & 4.1690 & 4.0255 & 3.9697 & 3.9473 & 3.3376 & 2.6333\\
3 & 20.6394 & 19.9346 & 19.6622 & 19.5531 & 19.5080 & 18.2283 & 17.6218 & 17.3890 & 17.2965 & 15.0372 & 12.2130\\
4 & 20.6560 & 19.9498 & 19.6783 & 19.5698 & 19.5251 & 18.2452 & 17.6387 & 17.4059 & 17.3134 & 15.0492 & 12.2253\\
5 & 20.6657 & 19.9529 & 19.6796 & 19.5704 & 19.5254 & 18.2457 & 17.6389 & 17.4060 & 17.3135 & 15.0518 & 12.2256\\
6 & 35.4198 & 34.0098 & 33.4700 & 33.2558 & 33.1678 & 32.1806 & 31.0512 & 30.6141 & 30.4394 & 26.2223 & 20.8931\\
7 & 35.4331 & 34.0165 & 33.4733 & 33.2574 & 33.1685 & 32.1839 & 31.0522 & 30.6144 & 30.4394 & 26.2245 & 20.8956\\
8 & 43.3830 & 41.5793 & 40.8896 & 40.6160 & 40.5037 & 41.3292 & 39.8556 & 39.2856 & 39.0578 & 33.7524 & 26.8513\\
9 & 271.4544 & 266.9576 & 265.7778 & 265.4780 & 265.4017 & 83.0086 & 80.1965 & 79.1253 & 78.7016 & 71.6959	 & 58.5692\\
10 & 271.5749 & 266.9740 & 265.7838 & 265.4848 & 265.4100 & 83.0336 & 80.2024 & 79.1260 & 78.7019 & 71.6980 & 58.5772\\
11 & 272.0985 & 267.1555 & 265.8297 & 265.4951 & 265.4123 & 83.0497 & 80.2096 & 79.1293 & 78.7031 & 71.7027 & 58.5823\\
12 & 272.4539 & 268.5083 & 267.5057 & 267.2524 & 267.1889 & 83.2255 & 80.4029 & 79.3294 & 78.9054 & 71.9215 & 58.7675\\
13 & 272.6653 & 268.5884 & 267.5274 & 267.2580 & 267.1903 & 83.2406 & 80.4093 & 79.3318 & 78.9064 & 71.9232 & 58.7823\\
14 & 273.1340 & 269.0642 & 268.1061 & 267.8663 & 267.8062 & 83.3102 & 80.4900 & 79.4147 & 78.9892 & 72.0256 & 58.8670\\
15 & 299.2086 & 293.1155 & 291.4105 & 290.9176 & 290.7663 & 110.1633 & 106.0977 & 104.5471 & 103.9332 & 96.1550 & 78.0084\\
16 & 316.7469 & 309.9446 & 307.9827	& 307.3882 & 307.1950 & 119.9714 & 115.5265 & 113.8296 & 113.1570 & 106.2227	 & 86.4177\\
17 & 316.9947 & 310.0058 & 307.9991 & 307.3930 & 307.1966 & 120.0147 & 115.5398 & 113.8336 & 113.1582 & 106.2293	 &  86.4276\\
18 & 344.9089 & 336.4442 & 333.9342 & 333.1492 & 332.8849 & 130.5284 & 125.6595 & 123.8065 & 123.0738 & 118.5657	 &  97.4488\\
19 & 345.3793 & 337.1256 & 334.6768 & 333.9079 & 333.6483 & 130.5647 & 125.7034 & 123.8531 & 123.1212 & 118.6246	 &  97.5097\\
20 & 345.4605 & 337.1269 & 334.6778 & 333.9091 & 333.6491 & 130.5856 & 125.7091 & 123.8549 & 123.1219 & 118.6412	 &  97.5150\\
21 & 427.1256 & 414.3260 & 410.2927 & 408.9367 & 408.4465 & 158.4240 & 152.2677 & 149.9313 & 149.0087 & 150.7050	 &  124.3511\\
22 & 427.2329 & 414.3685 & 410.3104 & 408.9440 & 408.4495 & 158.4260 & 152.2679 & 149.9313 & 149.0087 & 150.7129	 &  124.3577\\
23 & 427.7069 & 414.5074 & 410.3490 & 408.9552 & 408.4531 & 158.4307 & 152.2696 & 149.9319 & 149.0089 & 150.7148	 &  124.3641\\
24 & 437.5025 & 423.7340 & 419.2998 & 417.7744 & 417.2121 & 179.2826 & 171.6045 & 168.6461 & 167.4636 & 159.7685	 &  127.3438\\
25 & 437.6399 & 423.7789 & 419.3077 & 417.7749 & 417.2131 & 179.3092 & 171.6125 & 168.6488 & 167.4646 & 159.7975	 &  127.3581\\
26 & 462.1666 & 446.9485 & 441.9660 & 440.2315 & 439.5866 & 184.3542 & 176.4395 & 173.3946 & 172.1788 & 166.8725	 &  133.4426\\
27 & 851.9381 & 817.0677 & 807.8753 & 805.5351 & 804.9468 & 1071.5212 & 1057.8350 & 1053.1274 & 1051.9382 & 335.710 & 288.0482\\
28 & 888.8229 & 848.6904 & 837.8213 & 834.9749 & 834.2263 & 1071.7714 & 1057.8753 & 1053.1396 & 1051.9413 & 335.737 & 288.0994\\
29 & 891.2254 & 849.3639 & 837.9969 & 835.0193 & 834.2374 & 1071.8135 & 1058.3988 & 1053.2748 & 1051.975 &  335.770 & 288.1694\\
30 & 992.8633 & 939.1885 & 923.5286 & 919.3524 & 918.2220 & 1071.8613 & 1059.0201 & 1055.7633 & 1054.939 &  335.892 & 288.1905\\
31 & 993.9410 & 941.2736 & 927.4429 & 923.7731 & 922.7715 & 1072.3634 & 1059.1477 & 1055.7990 & 1054.949 &  335.896 & 288.2102\\
32 & 999.8911 & 941.9922 & 927.5956 & 923.8111 & 922.7812 & 1072.6331 & 1059.1996 & 1055.8151 & 1054.953 &  335.921 & 288.2306\\
33 & 1052.0007 & 999.5010 & 985.4672 & 981.7155 & 980.6827 & 1073.2824 & 1059.2529 & 1055.8720 & 1054.965 & 335.937 & 288.2510\\
34 & 1054.0103 & 999.9568 & 985.5786 & 981.7422 & 980.6887 & 1073.5674 & 1059.3885 & 1055.8825 & 1054.968 & 335.939 & 288.2655\\
35 & 1088.8728 & 1032.6001 & 1017.3610 & 1013.3103 & 1012.2129 & 1073.9888 & 1059.4116 & 1055.9128 & 1054.975 & 335.953 & 288.2700\\
36 & 1144.7398 & 1084.5089 & 1067.2785 & 1062.9436 & 1061.8501 & 1076.4724 & 1059.4875 & 1056.0661 & 1055.239 & 335.960 & 288.2956\\
37 & 1150.1233 & 1084.9562 & 1067.3690 & 1062.9740 & 1061.8796 & 1076.7105 & 1059.5993 & 1056.0783 & 1055.244 & 336.002 & 288.3556\\
38 & 1154.9383 & 1087.7050 & 1068.1575 & 1063.1712 & 1061.9241 & 1078.5174 & 1059.6378 & 1056.0826 & 1055.246 & 336.036 & 288.3577\\
39 & 1156.8863 & 1089.5572 & 1073.8220 & 1069.7891 & 1068.7754 & 1089.0059 & 1075.0220 & 1071.3853 & 1070.461 & 338.630 & 291.1745\\
40 & 1160.0776 & 1090.2141 & 1074.2545 & 1069.9123 & 1068.8069 & 1089.4789 & 1075.1464 & 1071.4168 & 1070.469 & 338.695 & 291.3120\\
41 & 1170.7057 & 1092.5885 & 1075.6761 & 1071.8657 & 1070.9267 & 1090.4701 & 1076.1214 & 1072.4046 & 1071.462 & 338.855 & 291.4274\\
42 & 1291.1778 & 1215.5339 & 1192.9766 & 1186.5037 & 1184.6079 & 1174.0043 & 1154.0648 & 1148.4476 & 1146.813 & 425.376 & 371.3145\\
43 & 1293.3068 & 1217.2739 & 1196.5209 & 1190.5539 & 1188.8129 & 1174.1797 & 1154.1186 & 1148.4619 & 1146.817 & 425.398 & 371.3541\\
44 & 1293.9845 & 1217.7969 & 1196.5696 & 1190.5670 & 1188.8201 & 1174.6064 & 1154.2334 & 1148.4912 & 1146.825 & 425.429 & 371.3888\\
45 & 1297.6383 & 1219.9161 & 1199.2968 & 1193.7045 & 1192.0944 & 1181.4177 & 1162.3500 & 1156.8737 & 1155.237 & 437.903 & 379.8021\\
46 & 1303.1037 & 1221.8710 & 1199.8517 & 1193.9074 & 1192.2241 & 1197.4002 & 1177.4067 & 1171.6223 & 1169.878 & 447.009 & 388.8530\\
47 & 1306.6105 & 1222.1916 & 1200.0086 & 1193.9475 & 1192.2329 & 1197.4897 & 1177.4369 & 1171.6306 & 1169.880 & 447.057 & 388.9220\\
48 & 1340.1246 & 1241.9893 & 1216.9260 & 1210.4150 & 1208.6378 & 1244.6256 & 1220.8687 & 1213.9726 & 1211.877 & 465.355 & 410.6224\\
49 & 1343.7226 & 1242.7779 & 1217.0955 & 1210.4525 & 1208.6464 & 1246.8784 & 1223.3028 & 1216.4472 & 1214.360 & 465.569 & 410.9257\\
50 & 1358.2827 & 1245.7145 & 1217.7373 & 1210.6107 & 1208.6875 & 1247.1785 & 1223.3827 & 1216.4669 & 1214.365 & 465.642 & 410.9518\\
51 & 1383.4011 & 1296.6487 & 1273.2079 & 1266.5890 & 1264.6036 & 1292.5390 & 1266.7158 & 1259.0430 & 1256.642 & 493.288 & 436.0154\\
52 & 1384.0920 & 1296.8827 & 1273.3351 & 1266.6775 & 1264.6515 & 1292.7459 & 1266.7771 & 1259.0591 & 1256.646 & 493.357 & 436.0628\\
53 & 1408.7496 & 1302.9443 & 1274.8446 & 1267.0475 & 1264.7924 & 1318.1064 & 1290.8384 & 1282.6980 & 1280.138 & 503.268 & 447.3366\\
54 &          & 1931.2295 & 1881.3515 & 1867.0606 & 1863.3648 & 1366.3125 & 1335.7326 & 1326.6378 & 1323.763 & 518.856 & 466.2302\\
55 &          & 1932.7859 & 1881.4113 & 1867.1235 & 1863.4643 & 1366.8361 & 1335.8772 & 1326.6733 & 1323.778 & 518.868 & 466.2960\\
56 &          & 1936.4743 & 1882.9905 & 1867.5235 & 1863.5513 & 1367.4013 & 1336.0234 & 1326.7149 & 1323.791 & 518.968 & 466.3212\\
57 &          & 1942.3175 & 1897.9042 & 1886.8888 & 1884.1212 & 1373.4237 & 1344.3897 & 1335.6043 & 1332.801 & 519.650 & 467.2563\\
58 &          & 1943.1835 & 1898.2372 & 1886.9745 & 1884.1428 & 1373.5194 & 1344.4135 & 1335.6090 & 1332.802 & 519.664 & 467.3472\\
59 &          & 1951.4159 & 1905.0061 & 1894.3829 & 1891.7581 & 1376.6480 & 1347.7240 & 1338.9737 & 1336.176 & 520.014 & 467.7506\\
60 &          & 1986.5116 & 1939.1087 & 1926.7835 & 1923.6106 & 1602.4894 & 1559.1123 & 1545.2809 & 1540.591 & 622.030 & 572.2704\\
    \end{tabular}
    }
    \caption{SG Unnormalized Eigenvalues}
    \label{tablethree1}
\end{table}

\begin{table}[htbp!]
  \resizebox{\tabfourwd}{!}{
  \begin{tabular}{r | r*{9}{r}r | r r}
    & \multicolumn{11}{p{9cm}}{Sierpinski Gasket Eigenvalue Data}
    & \multicolumn{2}{p{5cm}}{Spectral Decimation Eigenvalues}\\
    \hline
    Level: & 2 & 2 & 2 & 2 & 2 & 3 & 3 & 3 & 3 & 4 & 5 & \mbox{Actual} & \mbox{Actual}\\
    Refinement: & 0 & 1 & 2 & 3 & 4 & 0 & 1 & 2 & 3 & 0 & 0 & \mbox{Normalized} & \mbox{Unormalized}\\
    \hline
    n\\	
1 & 1.0000 & 1.0000 & 1.0000 & 1.0000 & 1.0000 & 1.0000 & 1.0000 & 1.0000 & 1.0000 & 1.0000	& 1.0000 & 1.0000 & 27.1144\\
2 & 1.0000 & 1.0001 & 1.0001 & 1.0000 & 1.0000 & 1.0000 & 1.0000 & 1.0000 & 1.0000 & 1.0001	& 1.0002 & 1.0000 &  27.1144\\
3 & 4.0687 & 4.0750 & 4.0773 & 4.0781 & 4.0784 & 4.3724 & 4.3775 & 4.3805 & 4.3819 & 4.5059	& 4.6390 & 5.0000 &  135.5721\\
4 & 4.0720 & 4.0781 & 4.0807 & 4.0816 & 4.0820 & 4.3764 & 4.3817 & 4.3847 & 4.3862 & 4.5095	& 4.6437 & 5.0000 &  135.5721\\
5 & 4.0739 & 4.0787 & 4.0809 & 4.0818 & 4.0820 & 4.3766 & 4.3818 & 4.3847 & 4.3862 & 4.5103	& 4.6437 & 5.0000 &  135.5721\\
6 & 6.9824 & 6.9522 & 6.9407 & 6.9361 & 6.9342 & 7.7191 & 7.7136 & 7.7120 & 7.7115 & 7.8575	& 7.9360 & 8.1039 &  219.7332\\
7 & 6.9850 & 6.9535 & 6.9414 & 6.9364 & 6.9343 & 7.7199 & 7.7139 & 7.7121 & 7.7115 & 7.8582	& 7.9370 & 8.1039 &  219.7332\\
8 & 8.5522 & 8.4995 & 8.4793 & 8.4712 & 8.4678 & 9.9136 & 9.9008 & 9.8964 & 9.8949 & 10.1139 & 10.1992 & 10.3056 & 279.4291\\
9 & 53.5124 & 54.5705 & 55.1142 & 55.3701 & 55.4858 & 19.9112 & 19.9221 & 19.9324 & 19.9381 & 21.4837 & 22.2469 & 25.0000 & 677.8606\\
10 & 53.5361 & 54.5739 & 55.1155 & 55.3715 & 55.4876 & 19.9172 & 19.9236 & 19.9326 & 19.9382 & 21.4843 & 22.2499 & 25.0000 & 677.8606\\
11 & 53.6394 & 54.6110 & 55.1250 & 55.3737 & 55.4880 & 19.9210 & 19.9254 & 19.9334 & 19.9385 & 21.4857 & 22.2519 & 25.0000 & 677.8606\\
12 & 53.7094 & 54.8875 & 55.4726 & 55.7402 & 55.8595 & 19.9632 & 19.9734 & 19.9838 & 19.9898 & 21.5513 & 22.3222 & 25.0000 & 677.8606\\
13 & 53.7511 & 54.9039 & 55.4771 & 55.7413 & 55.8598 & 19.9668 & 19.9750 & 19.9844 & 19.9900 & 21.5518 & 22.3278 & 25.0000 & 677.8606\\
14 & 53.8435 & 55.0012 & 55.5971 & 55.8682 & 55.9885 & 19.9835 & 19.9951 & 20.0053 & 20.0110 & 21.5825 & 22.3600 & 25.0000 & 677.8606\\
15 & 58.9836 & 59.9177 & 60.4297 & 60.6760 & 60.7887 & 26.4248 & 26.3564 & 26.3364 & 26.3303 & 28.8128 & 29.6307 & 51.5278 & 1397.1457\\
16 & 62.4410 & 63.3578 & 63.8663 & 64.1112 & 64.2233 & 28.7774 & 28.6987 & 28.6747 & 28.6670 & 31.8296 & 32.8248 & 35.1398 & 952.7966\\
17 & 62.4898 & 63.3703 & 63.8697 & 64.1122 & 64.2236 & 28.7878 & 28.7020 & 28.6757 & 28.6673 & 31.8316 & 32.8286 & 35.1398 & 952.7966\\
18 & 67.9926 & 68.7747 & 69.2478 & 69.4841 & 69.5941 & 31.3097 & 31.2159 & 31.1880 & 31.1794 & 35.5282 & 37.0149 & 40.5196 & 1098.6658\\
19 & 68.0854 & 68.9140 & 69.4018 & 69.6424 & 69.7537 & 31.3184 & 31.2268 & 31.1998 & 31.1913 & 35.5458 & 37.0380 & 40.5196 & 1098.6658\\
20 & 68.1014 & 68.9143 & 69.4020 & 69.6426 & 69.7539 & 31.3234 & 31.2282 & 31.2002 & 31.1915 & 35.5508 & 37.0400 & 40.5196 & 1098.6658\\
21 & 84.2002 & 84.6951 & 85.0823 & 85.2909 & 85.3913 & 38.0010 & 37.8258 & 37.7691 & 37.7496 & 45.1587 & 47.2335 &  8.1039 & 219.7331\\
22 & 84.2213 & 84.7038 & 85.0859 & 85.2925 & 85.3919 & 38.0015 & 37.8259 & 37.7691 & 37.7497 & 45.1611 & 47.2360 &  8.1039 & 219.7331\\
23 & 84.3148 & 84.7322 & 85.0939 & 85.2948 & 85.3927 & 38.0026 & 37.8263 & 37.7693 & 37.7497 & 45.1617 & 47.2384 & 31.7847 & 861.8226\\
24 & 86.2458 & 86.6182 & 86.9501 & 87.1342 & 87.2239 & 43.0044 & 42.6294 & 42.4836 & 42.4250 & 47.8746 & 48.3702 & 31.7847 & 861.8226\\
25 & 86.2729 & 86.6274 & 86.9517 & 87.1343 & 87.2241 & 43.0107 & 42.6314 & 42.4842 & 42.4253 & 47.8833 & 48.3756 & 31.7847 & 861.8226\\
26 & 91.1079 & 91.3637 & 91.6503 & 91.8180 & 91.9016 & 44.2209 & 43.8305 & 43.6797 & 43.6195 & 50.0033 & 50.6868 & 31.7847 & 861.8226\\
27 & 167.9444 & 167.0221 & 167.5288 & 168.0085 & 168.2851 & 257.0248 & 262.7838 & 265.2927 & 266.4966 & 100.5955 & 109.4121 & 202.5980 & 5493.3291\\
28 & 175.2155 & 173.4863 & 173.7387 & 174.1487 & 174.4064 & 257.0848 & 262.7938 & 265.2958 & 266.4974 & 100.6038 & 109.4315 & 202.5980 & 5493.3291\\
29 & 175.6892 & 173.6240 & 173.7751 & 174.1580 & 174.4087 & 257.0949 & 262.9239 & 265.3299 & 266.5060 & 100.6136 & 109.4581 & 202.5980 & 5493.3291\\
30 & 195.7253 & 191.9856 & 191.5118 & 191.7471 & 191.9668 & 257.1063 & 263.0782 & 265.9567 & 267.2570 & 100.6502 & 109.4661 & 202.5980 & 5493.3291\\
31 & 195.9377 & 192.4119 & 192.3235 & 192.6691 & 192.9179 & 257.2268 & 263.1099 & 265.9657 & 267.2593 & 100.6512 & 109.4736 & 202.5980 & 5493.3291\\
32 & 197.1107 & 192.5588 & 192.3552 & 192.6771 & 192.9200 & 257.2915 & 263.1228 & 265.9698 & 267.2604 & 100.6589 & 109.4813 & 202.5980 & 5493.3291\\
33 & 207.3831 & 204.3145 & 204.3560 & 204.7540 & 205.0251 & 257.4472 & 263.1360 & 265.9841 & 267.2634 & 100.6636 & 109.4891 & 202.5980 & 5493.3291\\
34 & 207.7793 & 204.4077 & 204.3791 & 204.7596 & 205.0263 & 257.5156 & 263.1697 & 265.9868 & 267.2643 & 100.6643 & 109.4946 & 202.5980 & 5493.3291\\
35 & 214.6518 & 211.0805 & 210.9698 & 211.3437 & 211.6169 & 257.6167 & 263.1754 & 265.9944 & 267.2660 & 100.6685 & 109.4963 & 202.5980 & 5493.3291\\
36 & 225.6650 & 221.6915 & 221.3211 & 221.6956 & 221.9942 & 258.2124 & 263.1943 & 266.0330 & 267.3330 & 100.6706 & 109.5060 & 202.5980 & 5493.3291\\
37 & 226.7263 & 221.7830 & 221.3399 & 221.7019 & 222.0004 & 258.2695 & 263.2221 & 266.0361 & 267.3342 & 100.6830 & 109.5288 & 202.5980 & 5493.3291\\
38 & 227.6754 & 222.3449 & 221.5034 & 221.7431 & 222.0097 & 258.7029 & 263.2316 & 266.0372 & 267.3347 & 100.6934 & 109.5296 & 202.5980 & 5493.3291\\
39 & 228.0595 & 222.7235 & 222.6781 & 223.1233 & 223.4420 & 261.2188 & 267.0533 & 269.8921 & 271.1893 & 101.4705 & 110.5995 & 202.5980 & 5493.3291\\
40 & 228.6886 & 222.8578 & 222.7678 & 223.1490 & 223.4486 & 261.3323 & 267.0843 & 269.9000 & 271.1913 & 101.4899 & 110.6518 & 202.5980 & 5493.3291\\
41 & 230.7837 & 223.3431 & 223.0626 & 223.5565 & 223.8918 & 261.5700 & 267.3264 & 270.1488 & 271.4429 & 101.5379 & 110.6956 & 202.5980 & 5493.3291\\
42 & 254.5326 & 248.4752 & 247.3871 & 247.4662 & 247.6584 & 281.6073 & 286.6889 & 289.3048 & 290.5322 & 127.4639 & 141.0399 & 158.9233 & 4309.1129\\
43 & 254.9523 & 248.8309 & 248.1221 & 248.3110 & 248.5375 & 281.6494 & 286.7022 & 289.3084 & 290.5332 & 127.4704 & 141.0549 & 158.9233 & 4309.1129\\
44 & 255.0859 & 248.9378 & 248.1322 & 248.3137 & 248.5390 & 281.7517 & 286.7308 & 289.3158 & 290.5351 & 127.4799 & 141.0681 & 158.9233 & 4309.1129\\
45 & 255.8062 & 249.3710 & 248.6978 & 248.9681 & 249.2236 & 283.3855 & 288.7471 & 291.4274 & 292.6663 & 131.2177 & 144.2638 & 158.9233 & 4309.1129\\
46 & 256.8836 & 249.7706 & 248.8128 & 249.0104 & 249.2507 & 287.2192 & 292.4874 & 295.1427 & 296.3753 & 133.9464 & 147.7017\\
47 & 257.5749 & 249.8362 & 248.8453 & 249.0188 & 249.2525 & 287.2407 & 292.4949 & 295.1448 & 296.3758 & 133.9608 & 147.7279\\
48 & 264.1816 & 253.8831 & 252.3535 & 252.4533 & 252.6822 & 298.5472 & 303.2841 & 305.8112 & 307.0152 & 139.4436 & 155.9706\\
49 & 264.8909 & 254.0443 & 252.3887 & 252.4612 & 252.6840 & 299.0875 & 303.8888 & 306.4345 & 307.6444 & 139.5079 & 156.0858\\
50 & 267.7611 & 254.6446 & 252.5217 & 252.4942 & 252.6926 & 299.1595 & 303.9086 & 306.4395 & 307.6455 & 139.5297 & 156.0957\\
51 & 272.7128 & 265.0564 & 264.0247 & 264.1694 & 264.3826 & 310.0401 & 314.6733 & 317.1648 & 318.3560 & 147.8138 & 165.6158\\
52 & 272.8490 & 265.1042 & 264.0510 & 264.1879 & 264.3926 & 310.0897 & 314.6885 & 317.1689 & 318.3570 & 147.8345 & 165.6338\\
53 & 277.7098 & 266.3433 & 264.3641 & 264.2651 & 264.4221 & 316.1729 & 320.6657 & 323.1237 & 324.3086 & 150.8042 & 169.9161\\
54 &          & 394.7752 & 390.1352 & 389.4083 & 389.5618 & 327.7361 & 331.8182 & 334.1926 & 335.3604 & 155.4753 & 177.0926\\
55 &          & 395.0934 & 390.1476 & 389.4215 & 389.5826 & 327.8617 & 331.8541 & 334.2015 & 335.3642 & 155.4788 & 177.1176\\
56 &          & 395.8473 & 390.4751 & 389.5049 & 389.6008 & 327.9972 & 331.8904 & 334.2120 & 335.3674 & 155.5087 & 177.1272\\
57 &          & 397.0418 & 393.5677 & 393.5439 & 393.9012 & 329.4418 & 333.9688 & 336.4513 & 337.6500 & 155.7132 & 177.4824\\
58 &          & 397.2188 & 393.6368 & 393.5617 & 393.9057 & 329.4648 & 333.9747 & 336.4525 & 337.6503 & 155.7173 & 177.5169\\
59 &          & 398.9017 & 395.0404 & 395.1069 & 395.4978 & 330.2152 & 334.7970 & 337.3001 & 338.5050 & 155.8221 & 177.6701\\
60 &          & 406.0758 & 402.1123 & 401.8646 & 402.1570 & 384.3876 & 387.3094 & 389.2708 & 390.2912 & 186.3913 & 217.3709\\
    \end{tabular}
    }
    \caption{SG Normalized Eigenvalues}
    \label{tablethree2}
\end{table}

We also look at the case $\epsilon=0$, where the $3^{m}$ triangles in $\Omega_{m}$ intersect at single points.  Thus the interior of $\Omega_{m}$ consists of $3^{m}$ disjoint triangles, and if we interpret the Neumann Laplacian on $\Omega_{m}$ in the usual way, the spectrum would just be $3^{m}$ copies of the spectrum of $\Omega$.  This is nothing like the spectrum of SG, and also it is not what we get when we use the FEM.  The reason is that the spline space chosen consists of continuous functions, and this effectively couples the disjoint triangles at their junction points.  Effectively this means that we are not looking at the entire Sobolev space $H^{1}(\Omega_{m})$, but only the subspace $H_{0}^{1}(\Omega_{m})$ defined to be the closure of continuous functions in $H^{1}(\Omega_{m})$ in the Sobolev norm.  In fact, functions in $H_{0}^{1}(\Omega_{m})$ do not have to be continuous (or even bounded), since $H^{1}$ does not embed in continuous functions on $\R^{2}$. They do have to satisfy some integral continuity condition (see \cite{strichartz67} for analogous results for $H^{1/2}$ on a half-line).  The conclusion is that the Neumann eigenvalues (and eigenfunctions) that the FEM approximates are the stationary values (and associated functions) for the Rayleigh quotient
\begin{equation}\label{three1}
R(u)=\frac{\int_{\Omega_{m}}\abs{\nabla u}^{2}dx}{\int_{\Omega_{m}}\abs{u}^{2}dx}
\end{equation}
for some $u\in H_{0}^{1}(\Omega_{m})$.  Of course some of these eigenfunctions restrict to Neumann eigenfunctions on each triangle in $\Omega_{m}$ and are continuous functions at the junction points, but it is easy to see that there are not enough of these (in fact the smallest such eigenvalue must be on the order of magnitude $4^{m}$).  We claim that all the other eigenfunctions have poles at some junction points.  Indeed, consider the restriction of an eigenfunction to a triangle. Because it is a Neumann eigenfunction, it must have vanishing normal derivatives along the side of the triangle.  Choose a vertex of the triangle and reflect the eigenfunction evenly six times around.  This yields an eigenfunction in a deleted neighborhood of the vertex.  The removable singularities theorem yields the following dichotomy: either the function is unbounded or it satisfies the eigenvalue equation at the vertex.  If it satisfies the eigenvalue equation at all three vertices of the triangle, then the restriction to the triangle is a Neumann eigenvalue, contrary to our assumption.  It is not difficult to see that the singularities must be logarithmic poles.

With this in mind, we look at the eigenvalue data in Tables \ref{tablethree3} and \ref{tablethree4}.  In contrast to our preceding computations, we do not see an apparent convergence of eigenvalues on a fixed $\Omega_{m}$ when we increase the refinement of the triangulation.  In particular the numerical values in Table \ref{tablethree4} are even better than the data in Table \ref{tablethree2}.  In other words, the poor approximations by the FEM to the actual eigenvalues on $\Omega_{m}$ yield very good approximations to the relative eigenvalues on SG.  We can even extract rather decent estimates for $r=1.25$ from the data in Table \ref{tablethree3} if we pair off corresponding refinements at different levels.  For example, if we compute $\lambda_{n}^{(3)}/\lambda_{n}^{(4)}$ using 3 refinements on level 3 and 4 refinements on level 4, the first six distinct eigenvalues yield ratios $1.246, 1.233, 1.223, 1.158, 1.128, 1.112$.

Of course the eigenfunctions on $\Omega_{m}$ cannot approximate the eigenfunctions on SG, since the latter are bounded.  Since we are already getting more information than we deserve, we might speculate that the eigenfunction approximation might be accurate in the complement of a small neighborhood of the junction points.

\begin{table}[htbp!]
  \resizebox{\tabfivewd}{!}{
  \begin{tabular}{r | r*{21}{r}r}
    \multicolumn{7}{p{9cm}}{Sierpinski Gasket, No Overlap, Unnormalized}\\
    \hline
    Level: & 1 & 1 & 1 & 2 & 2 & 2 & 3 & 3 & 4 & 4 \\
    Refinement: & 2 & 3 & 4 & 2 & 3 & 4 & 3 & 4 & 3 & 4\\
    \hline
    n\\	
1 &   3.650 &   3.113 &   2.713 &   3.689 &   3.103 &   2.677 &   2.773 &   2.364 &   2.687 &   2.224\\
2 &   3.721 &   3.164 &   2.752 &   3.689 &   3.103 &   2.677 &   2.773 &   2.364 &   2.687 &   2.224\\
3 &  70.334 &  70.221 &  70.193 &  17.179 &  14.285 &  12.220 &  13.684 &  11.644 &  13.409 &  11.091\\
4 &  70.356 &  70.227 &  70.195 &  17.179 &  14.285 &  12.220 &  13.684 &  11.644 &  13.409 &  11.091\\
5 &  70.362 &  70.228 &  70.195 &  17.179 &  14.285 &  12.220 &  13.684 &  11.644 &  13.409 &  11.091\\
6 &  82.579 &  80.476 &  79.025 &  25.961 &  21.402 &  18.196 &  21.944 &  18.643 &  21.699 &  17.941\\
7 &  82.850 &  80.663 &  79.163 &  25.961 &  21.402 &  18.196 &  21.944 &  18.643 &  21.699 &  17.941\\
8 &  96.289 &  91.743 &  88.598 &  31.117 &  25.509 &  21.604 &  27.689 &  23.498 &  27.564 &  22.783\\
9 & 212.039 & 210.923 & 210.645 & 282.881 & 281.270 & 280.869 &  63.390 &  53.398 &  66.354 &  54.735\\
10 & 235.397 & 230.816 & 227.972 & 282.881 & 281.270 & 280.869 &  63.390 &  53.398 &  66.354 &  54.735\\
11 & 236.040 & 231.194 & 228.236 & 282.881 & 281.270 & 280.869 &  63.390 &  53.398 &  66.354 &  54.735\\
12 & 283.145 & 281.338 & 280.886 & 282.881 & 281.270 & 280.869 &  63.390 &  53.398 &  66.354 &  54.735\\
13 & 283.383 & 281.397 & 280.901 & 282.881 & 281.270 & 280.869 &  63.390 &  53.398 &  66.354 &  54.735\\
14 & 283.654 & 281.464 & 280.918 & 282.881 & 281.270 & 280.869 &  63.390 &  53.398 &  66.354 &  54.735\\
15 & 303.428 & 296.864 & 293.464 & 302.450 & 297.037 & 294.131 &  78.165 &  65.625 &  84.053 &  69.271\\
16 & 304.318 & 297.305 & 293.740 & 316.088 & 308.066 & 303.425 &  85.035 &  71.278 &  92.754 &  76.407\\
17 & 310.855 & 303.820 & 299.961 & 316.088 & 308.066 & 303.425 &  85.035 &  71.278 &  92.754 &  76.407\\
18 & 497.970 & 492.957 & 491.705 & 342.012 & 329.003 & 321.049 &  95.415 &  79.778 & 106.636 &  87.776\\
19 & 499.255 & 493.281 & 491.786 & 342.012 & 329.003 & 321.049 &  95.415 &  79.778 & 106.636 &  87.776\\
20 & 500.036 & 493.472 & 491.833 & 342.012 & 329.003 & 321.049 &  95.415 &  79.778 & 106.636 &  87.776\\
21 & 526.723 & 515.905 & 511.065 & 391.489 & 368.308 & 353.783 & 110.130 &  91.738 & 128.376 & 105.547\\
22 & 528.063 & 516.482 & 511.415 & 391.489 & 368.308 & 353.783 & 110.130 &  91.738 & 128.376 & 105.547\\
23 & 568.769 & 547.592 & 536.643 & 408.152 & 381.181 & 364.323 & 114.069 &  94.922 & 134.759 & 110.756\\
24 & 642.681 & 634.411 & 632.344 & 408.152 & 381.181 & 364.323 & 114.069 &  94.922 & 134.759 & 110.756\\
25 & 644.193 & 634.780 & 632.435 & 408.152 & 381.181 & 364.323 & 114.069 &  94.922 & 134.759 & 110.756\\
26 & 644.399 & 634.848 & 632.454 & 408.152 & 381.181 & 364.323 & 114.069 &  94.922 & 134.759 & 110.756\\
27 & 646.983 & 635.458 & 632.604 & 861.618 & 847.027 & 843.410 & 1127.762 & 1124.145 & 311.734 & 253.559\\
28 & 677.294 & 661.371 & 654.619 & 887.555 & 867.807 & 860.994 & 1127.762 & 1124.145 & 311.734 & 253.559\\
29 & 681.117 & 662.685 & 655.217 & 887.555 & 867.807 & 860.994 & 1127.762 & 1124.145 & 311.734 & 253.559\\
30 & 864.923 & 847.924 & 843.639 & 973.982 & 938.678 & 921.579 & 1127.762 & 1124.145 & 311.734 & 253.559\\
31 & 888.687 & 867.510 & 860.707 & 973.982 & 938.678 & 921.579 & 1127.762 & 1124.145 & 311.734 & 253.559\\
32 & 890.096 & 868.076 & 861.020 & 973.982 & 938.678 & 921.579 & 1127.762 & 1124.145 & 311.734 & 253.559\\
33 & 938.355 & 918.889 & 914.016 & 1028.867 & 984.272 & 960.390 & 1127.762 & 1124.145 & 311.734 & 253.559\\
34 & 941.159 & 919.527 & 914.172 & 1028.867 & 984.272 & 960.390 & 1127.762 & 1124.145 & 311.734 & 253.559\\
35 & 944.371 & 920.288 & 914.359 & 1063.734 & 1013.093 & 984.352 & 1127.762 & 1124.145 & 311.734 & 253.559\\
36 & 990.412 & 956.321 & 943.405 & 1157.651 & 1131.525 & 1125.082 & 1127.762 & 1124.145 & 311.734 & 253.559\\
37 & 995.980 & 958.075 & 944.225 & 1157.651 & 1131.525 & 1125.082 & 1127.762 & 1124.145 & 311.734 & 253.559\\
38 & 1007.032 & 969.717 & 955.860 & 1157.651 & 1131.525 & 1125.082 & 1127.762 & 1124.145 & 311.734 & 253.559\\
39 & 1161.208 & 1132.507 & 1125.335 & 1157.651 & 1131.525 & 1125.082 & 1127.762 & 1124.145 & 311.734 & 253.559\\
40 & 1163.930 & 1133.179 & 1125.504 & 1157.651 & 1131.525 & 1125.082 & 1127.762 & 1124.145 & 311.734 & 253.559\\
41 & 1167.498 & 1134.086 & 1125.730 & 1157.651 & 1131.525 & 1125.082 & 1127.762 & 1124.145 & 311.734 & 253.559\\
42 & 1184.368 & 1147.961 & 1137.349 & 1230.800 & 1174.160 & 1154.340 & 1198.824 & 1182.609 & 386.255 & 312.662\\
43 & 1186.310 & 1148.403 & 1137.598 & 1247.495 & 1189.429 & 1167.359 & 1198.824 & 1182.609 & 386.255 & 312.662\\
44 & 1203.984 & 1162.140 & 1148.937 & 1247.495 & 1189.429 & 1167.359 & 1198.824 & 1182.609 & 386.255 & 312.662\\
45 & 1384.681 & 1346.242 & 1336.679 & 1267.855 & 1209.137 & 1185.248 & 1198.824 & 1182.609 & 386.255 & 312.662\\
46 & 1392.804 & 1348.470 & 1337.250 & 1267.855 & 1209.137 & 1185.248 & 1209.019 & 1191.017 & 394.037 & 318.795\\
47 & 1398.054 & 1349.707 & 1337.549 & 1267.855 & 1209.137 & 1185.248 & 1209.019 & 1191.017 & 394.037 & 318.795\\
48 & 1418.418 & 1369.638 & 1355.732 & 1289.621 & 1231.534 & 1206.796 & 1248.438 & 1223.543 & 421.196 & 340.141\\
49 & 1424.922 & 1371.372 & 1356.319 & 1289.621 & 1231.534 & 1206.796 & 1248.438 & 1223.543 & 421.196 & 340.141\\
50 & 1488.198 & 1415.530 & 1390.277 & 1294.363 & 1236.613 & 1211.852 & 1248.438 & 1223.543 & 421.196 & 340.141\\
51 & 1540.507 & 1490.281 & 1477.951 & 1294.363 & 1236.613 & 1211.852 & 1277.622 & 1247.622 & 439.088 & 354.151\\
52 & 1544.858 & 1491.604 & 1478.303 & 1294.363 & 1236.613 & 1211.852 & 1277.622 & 1247.622 & 439.088 & 354.151\\
53 & 1547.571 & 1492.089 & 1478.407 & 1294.363 & 1236.613 & 1211.852 & 1297.625 & 1264.115 & 450.532 & 363.091\\
54 & 1562.171 & 1495.591 & 1479.271 &  & 1991.484 & 1971.705 & 1342.696 & 1301.218 & 474.374 & 381.660\\
55 & 1628.290 & 1554.758 & 1529.753 &  & 1991.484 & 1971.705 & 1342.696 & 1301.218 & 474.374 & 381.660\\
56 & 1648.776 & 1560.742 & 1531.934 &  & 1991.484 & 1971.705 & 1342.696 & 1301.218 & 474.374 & 381.660\\
57 & 1847.262 & 1777.702 & 1760.374 &  & 1991.484 & 1971.705 & 1342.696 & 1301.218 & 474.374 & 381.660\\
58 & 1852.666 & 1779.553 & 1760.911 &  & 1991.484 & 1971.705 & 1342.696 & 1301.218 & 474.374 & 381.660\\
59 & 1862.961 & 1782.208 & 1761.547 &  & 1991.484 & 1971.705 & 1342.696 & 1301.218 & 474.374 & 381.660\\
60 & 1872.896 & 1793.161 & 1771.954 &  &  & 1999.469 & 1465.787 & 1401.719 & 529.475 & 424.272\\
    \end{tabular}
    }
    \caption{Sierpinski Gasket, No Overlap, Unnormalized}
    \label{tablethree3}
\end{table}

\begin{table}[htbp!]
  \resizebox{\tabfivewd}{!}{
  \begin{tabular}{r | r*{21}{r}r}
    \multicolumn{7}{p{9cm}}{Sierpinski Gasket, No Overlap, Normalized}\\
    \hline
    Level: & 1 & 1 & 1 & 2 & 2 & 2 & 3 & 3 & 4 & 4 \\
    Refinement: & 2 & 3 & 4 & 2 & 3 & 4 & 3 & 4 & 3 & 4\\
    \hline
    n\\	
  1 &   1.000 &   1.000 &   1.000 &   1.000 &   1.000 &   1.000 &   1.000 &   1.000 &   1.000 &   1.000\\
  2 &   1.019 &   1.016 &   1.014 &   1.000 &   1.000 &   1.000 &   1.000 &   1.000 &   1.000 &   1.000\\
  3 &  19.268 &  22.557 &  25.875 &   4.657 &   4.604 &   4.564 &   4.935 &   4.925 &   4.990 &   4.987\\
  4 &  19.274 &  22.559 &  25.876 &   4.657 &   4.604 &   4.564 &   4.935 &   4.925 &   4.990 &   4.987\\
  5 &  19.275 &  22.559 &  25.876 &   4.657 &   4.604 &   4.564 &   4.935 &   4.925 &   4.990 &   4.987\\
  6 &  22.622 &  25.851 &  29.131 &   7.037 &   6.897 &   6.796 &   7.914 &   7.886 &   8.075 &   8.067\\
  7 &  22.697 &  25.911 &  29.182 &   7.037 &   6.897 &   6.796 &   7.914 &   7.886 &   8.075 &   8.067\\
  8 &  26.378 &  29.470 &  32.660 &   8.435 &   8.221 &   8.069 &   9.986 &   9.939 &  10.257 &  10.244\\
  9 &  58.087 &  67.754 &  77.650 &  76.680 &  90.646 & 104.909 &  22.861 &  22.586 &  24.692 &  24.612\\
 10 &  64.486 &  74.144 &  84.038 &  76.680 &  90.646 & 104.909 &  22.861 &  22.586 &  24.692 &  24.612\\
 11 &  64.662 &  74.266 &  84.135 &  76.680 &  90.646 & 104.909 &  22.861 &  22.586 &  24.692 &  24.612\\
 12 &  77.566 &  90.373 & 103.543 &  76.680 &  90.646 & 104.909 &  22.861 &  22.586 &  24.692 &  24.612\\
 13 &  77.632 &  90.392 & 103.549 &  76.680 &  90.646 & 104.909 &  22.861 &  22.586 &  24.692 &  24.612\\
 14 &  77.706 &  90.414 & 103.555 &  76.680 &  90.646 & 104.909 &  22.861 &  22.586 &  24.692 &  24.612\\
 15 &  83.123 &  95.361 & 108.180 &  81.985 &  95.727 & 109.862 &  28.189 &  27.758 &  31.278 &  31.148\\
 16 &  83.367 &  95.502 & 108.282 &  85.682 &  99.281 & 113.334 &  30.667 &  30.149 &  34.517 &  34.357\\
 17 &  85.158 &  97.595 & 110.575 &  85.682 &  99.281 & 113.334 &  30.667 &  30.149 &  34.517 &  34.357\\
 18 & 136.417 & 158.351 & 181.258 &  92.709 & 106.029 & 119.916 &  34.410 &  33.744 &  39.682 &  39.469\\
 19 & 136.769 & 158.455 & 181.287 &  92.709 & 106.029 & 119.916 &  34.410 &  33.744 &  39.682 &  39.469\\
 20 & 136.983 & 158.517 & 181.305 &  92.709 & 106.029 & 119.916 &  34.410 &  33.744 &  39.682 &  39.469\\
 21 & 144.294 & 165.723 & 188.394 & 106.121 & 118.696 & 132.143 &  39.717 &  38.804 &  47.772 &  47.459\\
 22 & 144.661 & 165.908 & 188.523 & 106.121 & 118.696 & 132.143 &  39.717 &  38.804 &  47.772 &  47.459\\
 23 & 155.812 & 175.901 & 197.823 & 110.638 & 122.844 & 136.080 &  41.138 &  40.150 &  50.148 &  49.802\\
 24 & 176.060 & 203.790 & 233.102 & 110.638 & 122.844 & 136.080 &  41.138 &  40.150 &  50.148 &  49.802\\
 25 & 176.474 & 203.908 & 233.135 & 110.638 & 122.844 & 136.080 &  41.138 &  40.150 &  50.148 &  49.802\\
 26 & 176.530 & 203.930 & 233.142 & 110.638 & 122.844 & 136.080 &  41.138 &  40.150 &  50.148 &  49.802\\
 27 & 177.238 & 204.126 & 233.197 & 233.558 & 272.973 & 315.026 & 406.713 & 475.493 & 116.005 & 114.013\\
 28 & 185.542 & 212.450 & 241.313 & 240.589 & 279.670 & 321.593 & 406.713 & 475.493 & 116.005 & 114.013\\
 29 & 186.589 & 212.872 & 241.533 & 240.589 & 279.670 & 321.593 & 406.713 & 475.493 & 116.005 & 114.013\\
 30 & 236.942 & 272.376 & 310.991 & 264.016 & 302.510 & 344.223 & 406.713 & 475.493 & 116.005 & 114.013\\
 31 & 243.452 & 278.668 & 317.283 & 264.016 & 302.510 & 344.223 & 406.713 & 475.493 & 116.005 & 114.013\\
 32 & 243.838 & 278.850 & 317.399 & 264.016 & 302.510 & 344.223 & 406.713 & 475.493 & 116.005 & 114.013\\
 33 & 257.059 & 295.172 & 336.935 & 278.894 & 317.203 & 358.719 & 406.713 & 475.493 & 116.005 & 114.013\\
 34 & 257.827 & 295.377 & 336.992 & 278.894 & 317.203 & 358.719 & 406.713 & 475.493 & 116.005 & 114.013\\
 35 & 258.707 & 295.621 & 337.061 & 288.345 & 326.492 & 367.669 & 406.713 & 475.493 & 116.005 & 114.013\\
 36 & 271.319 & 307.196 & 347.768 & 313.804 & 364.659 & 420.234 & 406.713 & 475.493 & 116.005 & 114.013\\
 37 & 272.845 & 307.760 & 348.071 & 313.804 & 364.659 & 420.234 & 406.713 & 475.493 & 116.005 & 114.013\\
 38 & 275.872 & 311.499 & 352.360 & 313.804 & 364.659 & 420.234 & 406.713 & 475.493 & 116.005 & 114.013\\
 39 & 318.108 & 363.792 & 414.834 & 313.804 & 364.659 & 420.234 & 406.713 & 475.493 & 116.005 & 114.013\\
 40 & 318.854 & 364.008 & 414.896 & 313.804 & 364.659 & 420.234 & 406.713 & 475.493 & 116.005 & 114.013\\
 41 & 319.831 & 364.299 & 414.979 & 313.804 & 364.659 & 420.234 & 406.713 & 475.493 & 116.005 & 114.013\\
 42 & 324.453 & 368.756 & 419.262 & 333.632 & 378.399 & 431.163 & 432.341 & 500.222 & 143.736 & 140.589\\
 43 & 324.985 & 368.898 & 419.354 & 338.157 & 383.320 & 436.025 & 432.341 & 500.222 & 143.736 & 140.589\\
 44 & 329.826 & 373.311 & 423.534 & 338.157 & 383.320 & 436.025 & 432.341 & 500.222 & 143.736 & 140.589\\
 45 & 379.328 & 432.449 & 492.741 & 343.676 & 389.671 & 442.707 & 432.341 & 500.222 & 143.736 & 140.589\\
 46 & 381.553 & 433.165 & 492.952 & 343.676 & 389.671 & 442.707 & 436.018 & 503.779 & 146.633 & 143.347\\
 47 & 382.991 & 433.562 & 493.062 & 343.676 & 389.671 & 442.707 & 436.018 & 503.779 & 146.633 & 143.347\\
 48 & 388.570 & 439.965 & 499.765 & 349.576 & 396.889 & 450.756 & 450.234 & 517.537 & 156.739 & 152.945\\
 49 & 390.352 & 440.522 & 499.981 & 349.576 & 396.889 & 450.756 & 450.234 & 517.537 & 156.739 & 152.945\\
 50 & 407.686 & 454.706 & 512.499 & 350.862 & 398.526 & 452.644 & 450.234 & 517.537 & 156.739 & 152.945\\
 51 & 422.016 & 478.719 & 544.819 & 350.862 & 398.526 & 452.644 & 460.759 & 527.721 & 163.397 & 159.245\\
 52 & 423.208 & 479.144 & 544.948 & 350.862 & 398.526 & 452.644 & 460.759 & 527.721 & 163.397 & 159.245\\
 53 & 423.951 & 479.299 & 544.986 & 350.862 & 398.526 & 452.644 & 467.972 & 534.698 & 167.656 & 163.265\\
 54 & 427.951 & 480.424 & 545.305 &  & 641.800 & 736.460 & 484.227 & 550.392 & 176.528 & 171.614\\
 55 & 446.064 & 499.430 & 563.915 &  & 641.800 & 736.460 & 484.227 & 550.392 & 176.528 & 171.614\\
 56 & 451.676 & 501.353 & 564.718 &  & 641.800 & 736.460 & 484.227 & 550.392 & 176.528 & 171.614\\
 57 & 506.050 & 571.046 & 648.928 &  & 641.800 & 736.460 & 484.227 & 550.392 & 176.528 & 171.614\\
 58 & 507.530 & 571.641 & 649.126 &  & 641.800 & 736.460 & 484.227 & 550.392 & 176.528 & 171.614\\
 59 & 510.351 & 572.493 & 649.361 &  & 641.800 & 736.460 & 484.227 & 550.392 & 176.528 & 171.614\\
 60 & 513.073 & 576.012 & 653.197 &  &  & 746.830 & 528.618 & 592.902 & 197.033 & 190.775\\
    \end{tabular}
    }
    \caption{Sierpinski Gasket, No Overlap, Normalized}
    \label{tablethree4}
\end{table}

\clearpage
\section{Non-PCF Fractals}
\label{secnonpcf}

Our first example is the octagasket, generated by eight contractive homotheties with contraction ratio $1-\sqrt{2}/2$ and fixed points $\{q_{i}\}$ the vertices of a regular octagon.  Then the consecutive images $F_{i}K$ and $F_{i+1}K$ intersect along a Cantor set.  As yet, there has been no construction of a self-similar Laplacian on this fractal, although it is reasonable to expect the probabilistic methods in \cite{barlow95} will work, given the high symmetry in this example.  It is natural to approximate from without by taking $\Omega$ to be the interior of the octagon with vertices $\{q_{i}\}$.  Then $\Omega_{m}$ consists of the interior of the union of $8^{m}$ octagons that meet along edges.

\begin{figure}
\includegraphics[width=\threewidth\textwidth]{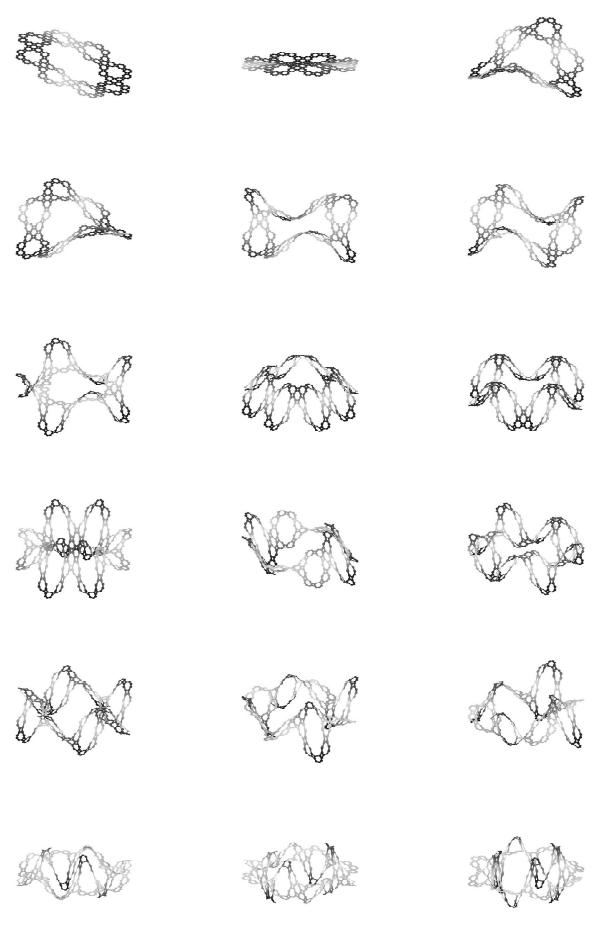}
\caption{Octagasket Eigenfunctions, Level 3}
\label{figfour1}
\end{figure}

\begin{figure}[htbp!]
\begin{center}
\includegraphics[width=\onewidth\textwidth]{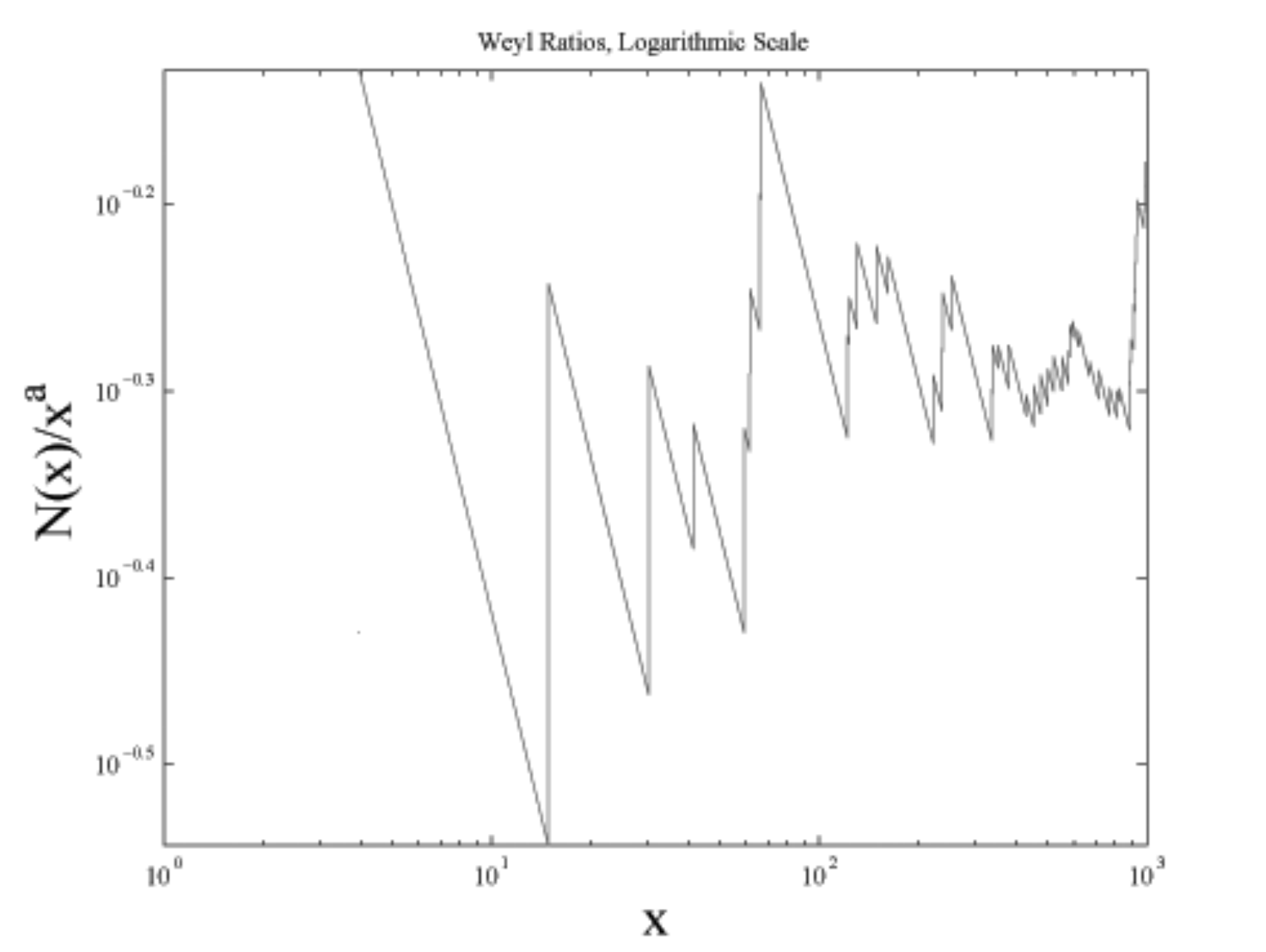}
\end{center}
\caption{Octagasket Weyl Ratios, Level 4, 1 Refinement, $\alpha=.71938$}
\label{figfour2}
\end{figure}

In table \ref{tablefour1} we give the eigenvalues on $\Omega_{m}$ for $m=0,1,2,3$ along with level-to-level ratios, suggesting a renormalization factor of about $r=1.2$.  In Table \ref{tablefour2} we normalize the eigenvalues by dividing by $\lambda_{1}^{(m)}$.  This suggests an eigenvalue renormalization factor of about $R=14.9476$ (the table indicates when a new eigenvalue appears that is approximately $R\lambda_{n}$ for an earlier value of $n$).  In the next section we will explain why this happens.  The tables show eigenvalues of multiplicities 1 and 2, but no higher multiplicities.  The $D_{8}$ symmetry forces multiplicity 2, since there are three irreducible representations of dimension 2.  There are a number of close coincidences (for example 910.5058 and 910.8645, each with multiplicity 2), but not close enough to be regarded as the same, in our judgement.  There is some evidence of large gaps in the spectrum, for example (66.45202,122.0411), (162.1709,223.2267) and (253.6123,336.1848).  However, there is not enough data to guess whether or not there are infinitely many gaps ($\lambda_{j+1}/\lambda_{j}\geq1+\epsilon$ for fixed $\epsilon$).  In Figure \ref{figfour1} we display the graphs of some eigenfunctions, and in Figure \ref{figfour2} we show the Weyl ratios.

The Weyl ratio is defined to be $W(x) = N(x)/x^{\alpha}$, where $N(x)=\#\{\lambda_{j} \leq x\}$ is the eigenvalue counting function, and $x^{\alpha}$ is its approximate growth rate. We determine $\alpha$ experimentally as the slope of the line of best fit to a log-log plot of $N(x)$. The Weyl ratio gives a nice ``snapshot'' of the spectrum. A question of interest is whether it tends to a limit, or shows periodic behavior for large $x$. Our experimental data does not give an indication of what answer to expect.

\begin{table}[htbp!]
  \resizebox{\tabsixwd}{!}{
  \begin{tabular}{r | r*{9}{r}r | r r r}
    \multicolumn{12}{p{8cm}}{Octagasket Unnormalized Eigenvalues}
    & \multicolumn{3}{p{2.3cm}}{Ratios $\lambda_{n}^{(j)}/\lambda_{n}^{(j+1)}$, highest refinements used}\\
    \hline
    Level: & 1 & 1 & 2 & 2 & 2 & 2 & 3 & 3 & 3 & 3 & 4  \\
    Refinement: & 0 & 1 & 1 & 2 & 3 & 4 & 1 & 2 & 3 & 4 & 1 \\
    \hline		
    n\\										
1 & 12.87 & 12.81 & 6.28 & 6.14 & 6.08 & 6.06 & 5.07 & 4.86 & 4.78 & 4.74 & 3.95 & 2.11 & 1.28 & 1.20\\
2 & 12.87 & 12.81 & 6.30 & 6.15 & 6.09 & 6.06 & 5.07 & 4.86 & 4.78 & 4.74 & 3.95 & 2.11 & 1.28 & 1.20\\
3 & 35.55 & 35.14 & 23.98 & 23.37 & 23.14 & 23.04 & 19.15 & 18.37 & 18.04 & 17.90 & 14.93 & 1.53 & 1.29 & 1.20\\
4 & 35.56 & 35.15 & 24.00 & 23.38 & 23.14 & 23.04 & 19.15 & 18.37 & 18.04 & 17.90 & 14.93 & 1.53 & 1.29 & 1.20\\
5 & 57.06 & 55.93 & 47.97 & 46.60 & 46.08 & 45.88 & 38.75 & 37.15 & 36.47 & 36.20 & 30.19 & 1.22 & 1.27 & 1.20\\
6 & 67.47 & 66.00 & 48.24 & 46.70 & 46.12 & 45.90 & 38.75 & 37.15 & 36.47 & 36.20 & 30.19 & 1.44 & 1.27 & 1.20\\
7 & 67.49 & 66.00 & 62.46 & 60.28 & 59.47 & 59.17 & 53.10 & 50.89 & 49.96 & 49.58 & 41.32 & 1.12 & 1.19 & 1.20\\
8 & 99.03 & 95.78 & 151.54 & 149.80 & 149.30 & 149.17 & 75.96 & 72.75 & 71.40 & 70.86 & 59.09 & 0.64 & 2.11 & 1.20\\
9 & 112.63 & 108.51 & 151.61 & 149.81 & 149.31 & 149.17 & 75.96 & 72.75 & 71.40 & 70.86 & 59.09 & 0.73 & 2.11 &	 1.20\\
10 & 112.78 & 108.55 & 154.94 & 152.84 & 152.22 & 152.04 & 78.99 & 75.65 & 74.25 & 73.68 & 61.42 & 0.71 & 2.06 & 1.20\\
11 & 124.95 & 120.29 & 155.27 & 152.92 & 152.24 & 152.04 & 78.99 & 75.65 & 74.25 & 73.68 & 61.42 & 0.79 & 2.06 & 1.20\\
12 & 166.13 & 157.78 & 161.09 & 158.27 & 157.40 & 157.14 & 84.90 & 81.30 & 79.79 & 79.18 & 65.96 & 1.00 & 1.98 & 1.20\\
13 & 166.21 & 157.80 & 161.44 & 158.37 & 157.43 & 157.14 & 84.90 & 81.30 & 79.79 & 79.18 & 65.96 & 1.00 & 1.98 & 1.20\\
14 & 179.85 & 169.80 & 162.62 & 159.42 & 158.46 & 158.17 & 85.54 & 81.91 & 80.39 & 79.78 & 66.45 & 1.07 & 1.98 & 1.20\\
15 & 180.29 & 169.87 & 163.28 & 159.63 & 158.52 & 158.19 & 85.54 & 81.91 & 80.39 & 79.78 & 66.45 & 1.07 & 1.98 & 1.20\\
16 & 201.34 & 188.62 & 223.56 & 213.46 & 210.19 & 209.10 & 157.34 & 150.33 & 147.46 & 146.31 & 122.04 & 0.90 &  1.43 & 1.20\\
17 & 237.11 & 220.91 & 225.86 & 217.24 & 214.55 & 213.67 & 157.34 & 150.33 & 147.46 & 146.31 & 122.04 & 1.03 &  1.46 & 1.20\\
18 & 237.65 & 221.00 & 227.27 & 217.60 & 214.66 & 213.72 & 159.54 & 152.43 & 149.51 & 148.34 & 123.73 & 1.03 &  1.44 & 1.20\\
19 & 274.23 & 251.90 & 251.17 & 242.41 & 239.73 & 238.90 & 168.29 & 160.76 & 157.67 & 156.43 & 130.46 & 1.05 &  1.53 & 1.20\\
20 & 274.43 & 251.95 & 251.31 & 242.44 & 239.75 & 238.90 & 168.29 & 160.76 & 157.67 & 156.43 & 130.46 & 1.05 &  1.53 & 1.20\\
21 & 289.57 & 265.40 & 290.09 & 278.88 & 275.54 & 274.51 & 193.38 & 184.60 & 181.03 & 179.60 & 149.71 & 0.97 &  1.53 & 1.20\\
22 & 290.74 & 265.65 & 292.73 & 279.78 & 275.84 & 274.61 & 193.38 & 184.60 & 181.03 & 179.60 & 149.71 & 0.97 &  1.53 & 1.20\\
23 & 335.78 & 304.52 & 314.75 & 299.70 & 295.20 & 293.80 & 209.59 & 199.98 & 196.08 & 194.53 & 162.17 & 1.04 &  1.51 & 1.20\\
24 & 336.40 & 304.62 & 428.94 & 413.83 & 409.52 & 408.37 & 290.88 & 277.07 & 271.53 & 269.34 & 223.23 & 0.75 &  1.52 & 1.21\\
25 & 360.47 & 322.62 & 431.35 & 414.35 & 409.64 & 408.40 & 290.88 & 277.07 & 271.53 & 269.34 & 223.23 & 0.79 &  1.52 & 1.21\\
26 & 379.58 & 339.03 & 448.72 & 432.14 & 427.30 & 425.95 & 309.23 & 294.38 & 288.46 & 286.11 & 237.13 & 0.80 &  1.49 & 1.21\\
27 & 400.07 & 356.68 & 452.24 & 432.96 & 427.49 & 425.99 & 309.23 & 294.38 & 288.46 & 286.11 & 237.13 & 0.84 &  1.49 & 1.21\\
28 & 428.17 & 381.56 & 453.90 & 433.38 & 427.60 & 426.02 & 310.63 & 295.70 & 289.75 & 287.39 & 238.19 & 0.90 &  1.48 & 1.21\\
29 & 441.65 & 390.22 & 459.24 & 434.84 & 427.99 & 426.12 & 310.63 & 295.70 & 289.75 & 287.39 & 238.19 & 0.92 &  1.48 & 1.21\\
30 & 445.04 & 390.90 & 479.16 & 456.52 & 450.02 & 448.16 & 331.23 & 315.14 & 308.75 & 306.23 & 253.61 & 0.87 &  1.46 & 1.21\\
31 & 486.74 & 425.77 & 488.40 & 458.88 & 450.63 & 448.32 & 331.23 & 315.14 & 308.75 & 306.23 & 253.61 & 0.95 &  1.46 & 1.21\\
32 & 500.94 & 434.79 & 562.07 & 535.94 & 528.39 & 526.24 & 437.21 & 414.29 & 405.37 & 401.88 & 336.18 & 0.83 &  1.31 & 1.20\\
33 & 503.98 & 435.36 & 565.73 & 537.94 & 529.17 & 526.46 & 437.21 & 414.29 & 405.37 & 401.88 & 336.18 & 0.83 &  1.31 & 1.20\\
34 & 517.82 & 450.69 & 571.43 & 539.01 & 530.14 & 527.96 & 441.53 & 418.44 & 409.46 & 405.95 & 338.61 & 0.85 &  1.30 & 1.20\\
35 & 519.43 & 451.06 & 573.83 & 539.58 & 530.53 & 528.06 & 441.53 & 418.44 & 409.46 & 405.95 & 338.61 & 0.85 &  1.30 & 1.20\\
36 & 607.09 & 518.07 & 584.89 & 549.09 & 539.29 & 536.57 & 459.56 & 435.33 & 425.94 & 422.28 & 351.80 & 0.97 &  1.27 & 1.20\\
37 & 611.05 & 518.98 & 643.32 & 605.90 & 595.21 & 592.22 & 493.31 & 466.56 & 456.26 & 452.25 & 379.14 & 0.88 &  1.31 & 1.19\\
38 & 638.35 & 538.29 & 651.28 & 607.91 & 595.75 & 592.37 & 493.31 & 466.56 & 456.26 & 452.25 & 379.14 & 0.91 &  1.31 & 1.19\\
39 & 644.93 & 539.72 & 700.98 & 661.74 & 651.18 & 648.45 & 557.35 & 525.93 & 513.96 & 509.33 & 428.58 & 0.83 &  1.27 & 1.19\\
40 & 659.58 & 557.54 & 749.55 & 703.59 & 690.24 & 686.48 & 588.09 & 554.32 & 541.50 & 536.55 & 451.74 & 0.81 &  1.28 & 1.19\\
41 & 665.39 & 558.99 & 759.91 & 706.48 & 690.87 & 686.62 & 588.09 & 554.32 & 541.50 & 536.55 & 451.74 & 0.81 &  1.28 & 1.19\\
42 & 669.03 & 569.71 & 764.98 & 710.86 & 695.42 & 689.38 & 616.83 & 580.78 & 567.16 & 561.91 & 473.92 & 0.83 &  1.23 & 1.19\\
43 & 671.06 & 569.92 & 768.33 & 711.94 & 696.68 & 692.85 & 616.83 & 580.78 & 567.16 & 561.91 & 473.92 & 0.82 &  1.23 & 1.19\\
44 & 752.63 & 614.76 & 780.01 & 716.52 & 696.98 & 692.92 & 649.02 & 610.27 & 595.70 & 590.10 & 499.69 & 0.89 &  1.17 & 1.18\\
45 & 799.56 & 664.07 & 790.04 & 724.17 & 709.42 & 705.55 & 649.02 & 610.27 & 595.70 & 590.10 & 499.69 & 0.94 &  1.20 & 1.18\\
46 & 802.16 & 664.43 & 792.39 & 727.47 & 710.39 & 705.82 & 673.48 & 632.54 & 617.22 & 611.34 & 519.17 & 0.94 &  1.15 & 1.18\\
47 & 815.03 & 679.96 & 806.00 & 740.88 & 723.51 & 718.84 & 673.48 & 632.54 & 617.22 & 611.34 & 519.17 & 0.95 &  1.18 & 1.18\\
48 & 832.53 & 685.05 & 847.00 & 785.98 & 769.26 & 764.87 & 712.51 & 668.11 & 651.58 & 645.25 & 551.07 & 0.90 &  1.19 & 1.17\\
49 & 840.23 & 686.13 & 849.40 & 786.33 & 769.32 & 764.88 & 712.51 & 668.11 & 651.58 & 645.25 & 551.07 & 0.90 &  1.19 & 1.17\\
50 & 870.15 & 712.75 & 889.09 & 821.77 & 803.21 & 798.22 & 737.80 & 690.37 & 672.83 & 666.14 & 576.94 & 0.89 &  1.20 & 1.15\\
51 & 895.83 & 718.89 & 899.51 & 824.78 & 803.97 & 798.40 & 737.80 & 690.37 & 672.83 & 666.14 & 576.94 & 0.90 &  1.20 & 1.15\\
52 & 904.07 & 722.76 & 951.70 & 867.27 & 844.89 & 838.85 & 749.87 & 701.93 & 684.17 & 677.40 & 582.26 & 0.86 &  1.24 & 1.16\\
53 & 951.97 & 779.21 & 956.27 & 869.23 & 845.43 & 838.99 & 749.87 & 701.93 & 684.17 & 677.40 & 582.26 & 0.93 &  1.24 & 1.16\\
54 & 977.47 & 790.31 & 984.78 & 894.52 & 869.56 & 862.77 & 763.18 & 713.92 & 695.71 & 688.77 & 593.76 & 0.92 &  1.25 & 1.16\\
55 & 983.56 & 791.59 & 987.21 & 895.16 & 869.79 & 862.84 & 769.50 & 718.12 & 699.26 & 692.09 & 619.00 & 0.92 &  1.25 & 1.12\\
56 & 1008.08 & 813.30 & 1099.26 & 990.90 & 961.30 & 953.18 & 968.08 & 909.14 & 887.65 & 879.56 & 668.25 & 0.85 & 1.08 & 1.32\\
58 & 1035.98 & 828.30 & 1101.07 & 996.03 & 962.98 & 953.67 & 1005.76 & 941.74 & 918.58 & 909.87 & 712.82 & 0.87	& 1.05 & 1.28\\
59 & 1059.87 & 841.04 & 1103.77 & 999.31 & 971.51 & 964.09 & 1005.76 & 941.74 & 918.58 & 909.87 & 712.82 & 0.87	& 1.06 & 1.28\\
60 & 1131.15 & 901.70 & 1117.52 & 1000.43 & 971.72 & 964.14 & 1078.51 & 1005.57 & 979.42 & 969.63 & 769.87 & 0.94 & 0.99 & 1.26\\
    \end{tabular}
    }
    \caption{Octagasket Unnormalized Eigenvalues and Ratios}
    \label{tablefour1}
\end{table}

\begin{table}[htbp!]
  \resizebox{\tabthreewda}{!}{
  \begin{tabular}{r | r*{13}{r}}
    \multicolumn{8}{p{8cm}}{Octagasket Normalized Eigenvalues}\\
    \hline
    Level: & 1 & 1 & 2 & 2 & 2 & 2 & 3 & 3 & 3 & 3 & 4  \\
    Refinement: & 0 & 1 & 1 & 2 & 3 & 4 & 1 & 2 & 3 & 4 & 1 \\
    \hline
    n\\
 1 &    1.00 &    1.00 &    1.00 &    1.00 &    1.00 &    1.00 &    1.00 &    1.00 &    1.00 &    \textbf{1.00} &    1.00\\
 2 &    1.00 &    1.00 &    1.00 &    1.00 &    1.00 &    1.00 &    1.00 &    1.00 &    1.00 &    \textbf{1.00} &    1.00\\
 3 &    2.76 &    2.74 &    3.82 &    3.81 &    3.80 &    3.80 &    3.78 &    3.78 &    3.78 &    \textbf{3.78} &    3.78\\
 4 &    2.76 &    2.74 &    3.82 &    3.81 &    3.80 &    3.80 &    3.78 &    3.78 &    3.78 &    \textbf{3.78} &    3.78\\
 5 &    4.43 &    4.37 &    7.64 &    7.59 &    7.57 &    7.57 &    7.64 &    7.64 &    7.64 &    \textbf{7.64} &    7.64\\
 6 &    5.24 &    5.15 &    7.68 &    7.60 &    7.58 &    7.57 &    7.64 &    7.64 &    7.64 &    \textbf{7.64} &    7.64\\
 7 &    5.25 &    5.15 &    9.95 &    9.82 &    9.77 &    9.76 &   10.48 &   10.46 &   10.46 &   10.46 &   10.45\\
 8 &    7.70 &    7.48 &   24.13 &   24.39 &   24.54 &   24.61 &   14.99 &   14.96 &   14.95 &   14.95 &   \textbf{14.95}\\
 9 &    8.75 &    8.47 &   24.14 &   24.39 &   24.54 &   24.61 &   14.99 &   14.96 &   14.95 &   14.95 &   \textbf{14.95}\\
10 &    8.77 &    8.47 &   24.67 &   24.89 &   25.02 &   25.08 &   15.58 &   15.55 &   15.54 &   15.54 &   15.54\\
11 &    9.71 &    9.39 &   24.72 &   24.90 &   25.02 &   25.08 &   15.58 &   15.55 &   15.54 &   15.54 &   15.54\\
12 &   12.91 &   12.31 &   25.65 &   25.77 &   25.87 &   25.92 &   16.75 &   16.71 &   16.70 &   16.70 &   16.69\\
13 &   12.92 &   12.32 &   25.70 &   25.79 &   25.87 &   25.92 &   16.75 &   16.71 &   16.70 &   16.70 &   16.69\\
14 &   13.98 &   13.25 &   25.89 &   25.96 &   26.04 &   26.09 &   16.88 &   16.84 &   16.83 &   16.83 &   16.81\\
15 &   14.01 &   13.26 &   26.00 &   25.99 &   26.05 &   26.09 &   16.88 &   16.84 &   16.83 &   16.83 &   16.81\\
16 &   15.65 &   14.72 &   35.59 &   34.76 &   34.54 &   34.49 &   31.04 &   30.91 &   30.87 &   30.86 &   30.87\\
17 &   18.43 &   17.24 &   35.96 &   35.37 &   35.26 &   35.25 &   31.04 &   30.91 &   30.87 &   30.86 &   30.87\\
18 &   18.47 &   17.25 &   36.18 &   35.43 &   35.28 &   35.25 &   31.47 &   31.34 &   31.30 &   31.29 &   31.30\\
19 &   21.31 &   19.66 &   39.99 &   39.47 &   39.40 &   39.41 &   33.20 &   33.05 &   33.01 &   32.99 &   33.00\\
20 &   21.33 &   19.66 &   40.01 &   39.48 &   39.40 &   39.41 &   33.20 &   33.05 &   33.01 &   32.99 &   33.00\\
21 &   22.51 &   20.71 &   46.19 &   45.41 &   45.28 &   45.28 &   38.15 &   37.95 &   37.90 &   37.88 &   37.87\\
22 &   22.60 &   20.73 &   46.61 &   45.56 &   45.33 &   45.30 &   38.15 &   37.95 &   37.90 &   37.88 &   37.87\\
23 &   26.10 &   23.77 &   50.11 &   48.80 &   48.51 &   48.46 &   41.35 &   41.11 &   41.05 &   41.03 &   41.02\\
24 &   26.15 &   23.78 &   68.29 &   67.39 &   67.30 &   67.36 &   57.38 &   56.96 &   56.84 &   56.81 &   \textbf{56.47}\\
25 &   28.02 &   25.18 &   68.68 &   67.47 &   67.32 &   67.37 &   57.38 &   56.96 &   56.84 &   56.81 &   \textbf{56.47}\\
26 &   29.50 &   26.46 &   71.44 &   70.37 &   70.22 &   70.26 &   61.00 &   60.52 &   60.39 &   60.35 &   59.98\\
27 &   31.09 &   27.84 &   72.00 &   70.50 &   70.26 &   70.27 &   61.00 &   60.52 &   60.39 &   60.35 &   59.98\\
28 &   33.28 &   29.78 &   72.27 &   70.57 &   70.27 &   70.27 &   61.28 &   60.79 &   60.66 &   60.62 &   60.25\\
29 &   34.33 &   30.46 &   73.12 &   70.81 &   70.34 &   70.29 &   61.28 &   60.79 &   60.66 &   60.62 &   60.25\\
30 &   34.59 &   30.51 &   76.29 &   74.34 &   73.96 &   73.92 &   65.34 &   64.79 &   64.64 &   64.59 &   64.15\\
31 &   37.83 &   33.23 &   77.76 &   74.72 &   74.06 &   73.95 &   65.34 &   64.79 &   64.64 &   64.59 &   64.15\\
32 &   38.93 &   33.93 &   89.49 &   87.27 &   86.84 &   86.80 &   86.25 &   85.17 &   84.86 &   84.76 &   85.04\\
33 &   39.17 &   33.98 &   90.07 &   87.59 &   86.96 &   86.84 &   86.25 &   85.17 &   84.86 &   84.76 &   85.04\\
34 &   40.25 &   35.18 &   90.98 &   87.77 &   87.12 &   87.09 &   87.10 &   86.02 &   85.72 &   85.62 &   85.65\\
35 &   40.37 &   35.20 &   91.36 &   87.86 &   87.19 &   87.10 &   87.10 &   86.02 &   85.72 &   85.62 &   85.65\\
36 &   47.18 &   40.43 &   93.12 &   89.41 &   88.63 &   88.51 &   90.66 &   89.50 &   89.17 &   89.07 &   88.99\\
37 &   47.49 &   40.51 &  102.42 &   98.66 &   97.82 &   97.69 &   97.32 &   95.92 &   95.52 &   95.39 &   95.91\\
38 &   49.61 &   42.01 &  103.69 &   98.99 &   97.91 &   97.71 &   97.32 &   95.92 &   95.52 &   95.39 &   95.91\\
39 &   50.13 &   42.12 &  111.60 &  107.75 &  107.02 &  106.96 &  109.95 &  108.12 &  107.60 &  107.43 &  108.41\\
40 &   51.26 &   43.52 &  119.34 &  114.57 &  113.44 &  113.23 &  116.01 &  113.96 &  113.36 &  113.17 &  \textbf{114.27}\\
41 &   51.72 &   43.63 &  120.99 &  115.04 &  113.54 &  113.26 &  116.01 &  113.96 &  113.36 &  113.17 &  \textbf{114.27}\\
42 &   52.00 &   44.47 &  121.79 &  115.75 &  114.29 &  113.71 &  121.68 &  119.40 &  118.73 &  118.52 &  119.88\\
43 &   52.16 &   44.48 &  122.33 &  115.93 &  114.49 &  114.29 &  121.68 &  119.40 &  118.73 &  118.52 &  119.88\\
44 &   58.50 &   47.98 &  124.19 &  116.67 &  114.54 &  114.30 &  128.03 &  125.46 &  124.71 &  124.46 &  126.40\\
45 &   62.14 &   51.83 &  125.78 &  117.92 &  116.59 &  116.38 &  128.03 &  125.46 &  124.71 &  124.46 &  126.40\\
46 &   62.35 &   51.86 &  126.16 &  118.46 &  116.75 &  116.42 &  132.86 &  130.04 &  129.21 &  128.94 &  131.33\\
47 &   63.34 &   53.07 &  128.32 &  120.64 &  118.90 &  118.57 &  132.86 &  130.04 &  129.21 &  128.94 &  131.33\\
48 &   64.70 &   53.47 &  134.85 &  127.98 &  126.42 &  126.16 &  140.56 &  137.35 &  136.40 &  136.09 &  139.40\\
49 &   65.30 &   53.55 &  135.23 &  128.04 &  126.43 &  126.17 &  140.56 &  137.35 &  136.40 &  136.09 &  139.40\\
50 &   67.63 &   55.63 &  141.55 &  133.81 &  132.00 &  131.67 &  145.55 &  141.93 &  140.85 &  140.50 &  145.94\\
51 &   69.62 &   56.11 &  143.21 &  134.30 &  132.13 &  131.70 &  145.55 &  141.93 &  140.85 &  140.50 &  145.94\\
52 &   70.27 &   56.41 &  151.52 &  141.22 &  138.85 &  138.37 &  147.93 &  144.30 &  143.23 &  142.87 &  147.29\\
53 &   73.99 &   60.82 &  152.25 &  141.54 &  138.94 &  138.39 &  147.93 &  144.30 &  143.23 &  142.87 &  147.29\\
54 &   75.97 &   61.68 &  156.79 &  145.66 &  142.90 &  142.31 &  150.55 &  146.77 &  145.64 &  145.27 &  150.20\\
55 &   76.44 &   61.78 &  157.17 &  145.76 &  142.94 &  142.32 &  151.80 &  147.63 &  146.39 &  145.97 &  156.58\\
56 &   78.35 &   63.48 &  175.01 &  161.35 &  157.98 &  157.23 &  190.97 &  186.90 &  185.83 &  185.51 &  169.04\\
57 &   80.52 &   64.65 &  175.30 &  162.19 &  158.26 &  157.31 &  198.41 &  193.61 &  192.30 &  191.91 &  180.32\\
58 &   82.37 &   65.64 &  175.73 &  162.72 &  159.66 &  159.03 &  198.41 &  193.61 &  192.30 &  191.91 &  180.32\\
59 &   87.91 &   70.38 &  177.92 &  162.90 &  159.69 &  159.03 &  212.76 &  206.73 &  205.04 &  204.51 &  194.75\\
60 &   88.87 &   70.49 &  195.56 &  174.63 &  169.25 &  167.91 &  212.76 &  206.73 &  205.04 &  204.51 &  194.75\\
    \end{tabular}
    }
    \caption{Octagasket Normalized Eigenvalues.  Eigenvalues in boldface on level $4$ are
approximately $R$ $(14.95)$ times the eigenvalues in boldface on level $3$}
    \label{tablefour2}
\end{table}

\begin{figure}
$$
\begin{array}{ccc}
\includegraphics[width=\threewidth\textwidth]{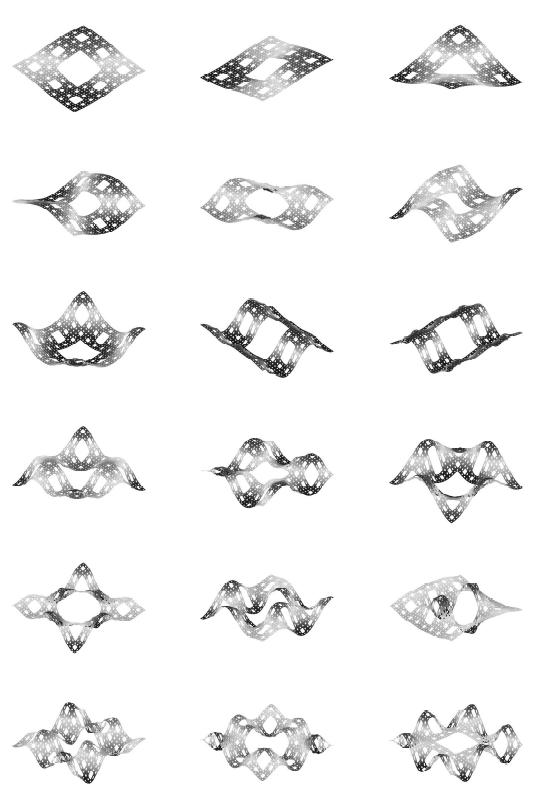}
\end{array}
$$
\caption{Sierpinski Carpet (SC) Eigenfunctions, Level 4}
\label{figfour3}
\end{figure}

\begin{figure}[htbp!]
\begin{center}
\includegraphics[width=\onewidth\textwidth]{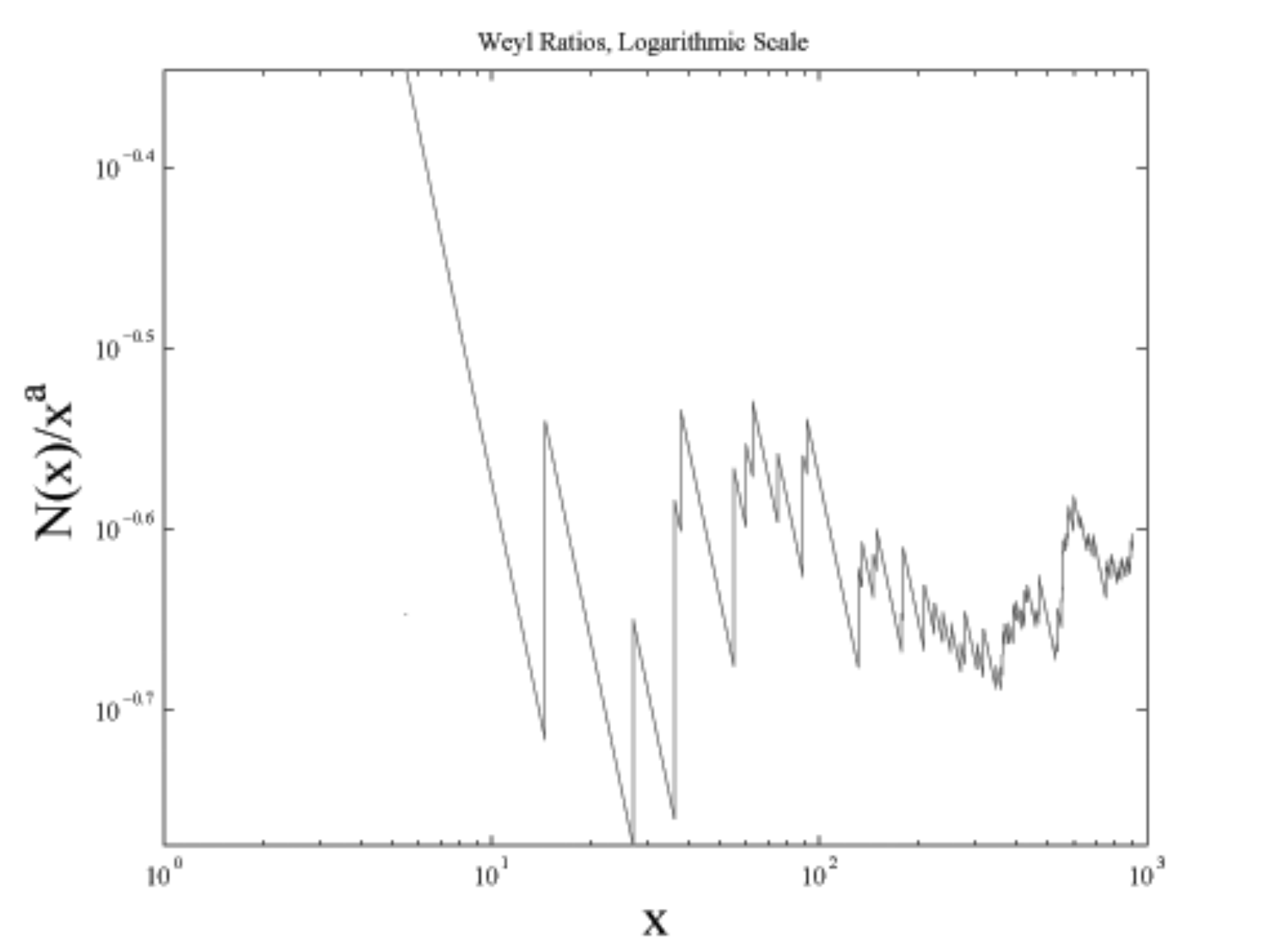}
\end{center}
\caption{SC Weyl Ratios, Level 3, 3 Refinements, $\alpha=.87392$}
\label{figfour4}
\end{figure}

The next example we consider is the standard SC generated by eight contractions of ratio $\frac{1}{3}$ (omitting the middle tic-tac-toe square).  Here the existence of a self-similar Laplacian is known, and as stated above, uniqueness is established in \cite{barlow08}.  Here it is natural to choose $\Omega$ to be the interior of the square that just contains SC, so $\Omega_{m}$ contains $8^{m}$ squares of side length $3^{-m}$ intersecting along edges.  In Tables \ref{tablefour3} and \ref{tablefour4} we report unnormalized and normalized eigenvalue data, as before.  In Table \ref{tablefour5} we describe the $D_{4}$ representation type associated to the eigenspace.  There is one 2-dimensional representation (denoted 2) and four 1-dimensional representations ($1++,1+-,1-+$, and $1--$) described in more detail in the next section.  Again we only see eigenvalue multiplicities of 1 or 2.  There is an apparent eigenvalue renormalization factor of about $R=10.0081$, which is consistent with computations in \cite{barlow90}.  In the next section we will give an explanation of this behavior.  Spectral gaps are consistent with the data.  Figure \ref{figfour3} shows some eigenfunctions and Figure \ref{figfour4} shows the Weyl ratios.

\begin{table}[htbp!]
  \resizebox{\tabfourwd}{!}{
  \begin{tabular}{r | r*{21}{r}r}
    \multicolumn{7}{p{9cm}}{SC Unnormalized Eigenvalue Data}\\
    \hline
    Level: & 1 & 1 & 1 & 1 & 1 & 2 & 2 & 2 & 2 & 3 & 3 & 3 & 3\\
    Refinement: & 0 & 1 & 2 & 3 & 4 & 0 & 1 & 2 & 3 & 0 & 1 & 2 & 3\\
    \hline
    n\\	
1 & 6.9095 & 6.8043 & 6.7653 & 6.7505 & 6.7449 & 6.4313 & 6.2251 & 6.1354 & 6.10 & 5.87 & 5.6310 & 5.5325 & 5.4936\\
2 & 6.9151 & 6.8070 & 6.7664 & 6.7510 & 6.7451 & 6.4323 & 6.2254 & 6.1355 & 6.10 & 5.87 & 5.6310 & 5.5325 & 5.4936\\
3 & 18.9925 & 18.8600 & 18.8243 & 18.8150 & 18.8127 & 17.0488 & 16.5574 & 16.3441 & 16.26 & 15.57 & 14.9425 & 14.6841 & 14.5823\\
4 & 34.3954 & 33.3153 & 32.9506 & 32.8218 & 32.7746 & 32.1097 & 30.7611 & 30.2096 & 30.00 & 28.99 & 27.7633 & 27.2621 & 27.0651\\
5 & 47.0332 & 45.9328 & 45.6120 & 45.5142 & 45.4829 & 42.8446 & 41.1639 & 40.4941 & 40.24 & 38.71 & 37.0938 & 36.4341 & 36.1754\\
6 & 47.1243 & 45.9546 & 45.6175 & 45.5157 & 45.4834 & 42.8729 & 41.1729 & 40.4965 & 40.24 & 38.71 & 37.0938 & 36.4341 & 36.1754\\
7 & 52.1193 & 51.1387 & 50.8862 & 50.8221 & 50.8059 & 44.8902 & 43.1147 & 42.3923 & 42.12 & 40.48 & 38.7875 & 38.0972 & 37.8262\\
8 & 92.8584 & 89.8444 & 89.0822 & 88.8905 & 88.8425 & 66.1417 & 62.7984 & 61.4619 & 60.97 & 58.99 & 56.4159 & 55.3699 & 54.9598\\
9 & 93.0125 & 89.8865 & 89.0933 & 88.8933 & 88.8432 & 66.2982 & 62.8473 & 61.4753 & 60.97 & 58.99 & 56.4159 & 55.3699 & 54.9598\\
10 & 99.0749 & 95.5391 & 94.6410 & 94.4145 & 94.3576 & 71.9960 & 68.3327 & 66.9090 & 66.39 & 64.16 & 61.3551 & 60.2173 & 59.7724\\
11 & 103.3472 & 99.2335 & 98.1700 & 97.8905 & 97.8154 & 75.6967 & 71.6324 & 70.0686 & 69.50 & 67.31 & 64.3456 & 63.1414 & 62.6712\\
12 & 116.5071 & 111.4375 & 110.1583 & 109.8362 & 109.7554 & 90.1313 & 85.5167 & 83.8352 & 83.24 & 80.14 & 76.6848 & 75.2891 & 74.7457\\
13 & 138.4478 & 129.5749 & 127.0433 & 126.2799 & 126.0367 & 108.9329 & 102.2554 & 99.8578 & 99.01 & 96.14 & 91.8300 & 90.0930 & 89.4182\\
14 & 139.6143 & 132.0182 & 130.0016 & 129.4432 & 129.2825 & 109.1286 & 102.3168 & 99.8748 & 99.01 & 96.14 & 91.8300 & 90.0930 & 89.4182\\
15 & 140.8159 & 132.2970 & 130.0662 & 129.4575 & 129.2854 & 113.3883 & 104.7730 & 101.6793 & 100.57 & 98.96 & 94.3752 & 92.5219 & 91.7984\\
16 & 176.7080 & 165.3371 & 162.4409 & 161.6916 & 161.4947 & 161.7668 & 152.1424 & 148.7128 & 147.54 & 141.59 & 135.1528 & 132.6025 & 131.6229\\
17 & 177.1511 & 165.4399 & 162.4692 & 161.6997 & 161.4971 & 162.4681 & 152.3584 & 148.7706 & 147.56 & 141.59 & 135.1528 & 132.6025 & 131.6229\\
18 & 186.3468 & 174.2256 & 171.1654 & 170.3923 & 170.1977 & 168.6568 & 155.5422 & 151.3127 & 149.93 & 145.60 & 138.8378 & 136.1567 & 135.1289\\
19 & 196.9752 & 182.4076 & 178.8404 & 177.9499 & 177.7272 & 178.1670 & 171.8544 & 169.9944 & 169.49 & 156.72 & 150.0903 & 147.4883 & 146.4968\\
20 & 207.4874 & 192.1334 & 188.2268 & 187.2187 & 186.9544 & 185.7692 & 174.5895 & 170.8667 & 169.66 & 161.04 & 153.7589 & 150.8882 & 149.7908\\
21 & 255.1933 & 232.4075 & 226.4589 & 224.8390 & 224.3803 & 224.8288 & 207.6838 & 200.9126 & 198.63 & 193.66 & 184.0351 & 180.2588 & 178.8197\\
22 & 258.1854 & 233.1381 & 226.6316 & 224.8783 & 224.3886 & 225.7046 & 211.4197 & 207.1909 & 205.92 & 194.18 & 185.2580 & 181.7748 & 180.4537\\
23 & 281.0056 & 253.6079 & 246.3701 & 244.3644 & 243.7824 & 228.0738 & 211.4581 & 207.2012 & 205.93 & 194.18 & 185.2580 & 181.7748 & 180.4537\\
24 & 294.4730 & 267.9779 & 261.2521 & 259.5384 & 259.1060 & 267.3355 & 245.9265 & 238.5531 & 236.14 & 226.58 & 215.4968 & 211.1887 & 209.5562\\
25 & 296.0104 & 268.7574 & 261.9548 & 260.2517 & 259.8246 & 268.9097 & 246.0942 & 238.5881 & 236.15 & 226.58 & 215.4968 & 211.1887 & 209.5562\\
26 & 299.1848 & 269.4117 & 262.0976 & 260.2854 & 259.8329 & 287.3782 & 265.0656 & 258.6422 & 256.75 & 242.91 & 231.2213 & 226.7026 & 225.0039\\
27 & 331.3359 & 297.2510 & 288.6697 & 286.5081 & 285.9648 & 308.4704 & 283.9907 & 277.1404 & 275.23 & 256.13 & 244.2915 & 239.7481 & 238.0328\\
28 & 375.1390 & 328.3141 & 316.2114 & 312.9576 & 312.0478 & 327.3493 & 299.1712 & 290.6243 & 288.04 & 273.39 & 259.6175 & 254.3202 & 252.3353\\
29 & 411.4908 & 363.1051 & 350.9457 & 347.8623 & 347.0751 & 358.6256 & 317.3075 & 302.3844 & 297.35 & 296.49 & 279.3908 & 272.8143 & 270.3526\\
30 & 415.5502 & 364.2623 & 351.2403 & 347.9359 & 347.0933 & 366.8938 & 332.0219 & 321.4484 & 318.28 & 302.48 & 286.9005 & 280.9431 & 278.7202\\
31 & 422.8709 & 372.1289 & 359.5109 & 356.3576 & 355.5688 & 370.3870 & 332.2614 & 321.5045 & 318.29 & 302.48 & 286.9005 & 280.9431 & 278.7202\\
32 & 428.1018 & 373.1788 & 359.7477 & 356.4153 & 355.5831 & 402.0958 & 361.9981 & 350.5678 & 347.27 & 328.13 & 311.5196 & 305.2009 & 302.8469\\
33 & 436.6960 & 380.8205 & 367.3222 & 363.9749 & 363.1377 & 426.2507 & 379.0875 & 363.7003 & 358.88 & 346.65 & 327.5608 & 320.3207 & 317.6237\\
34 & 439.9719 & 381.6613 & 367.5250 & 364.0248 & 363.1501 & 431.0532 & 380.0237 & 363.9265 & 358.94 & 346.65 & 327.5608 & 320.3207 & 317.6237\\
35 & 451.1346 & 393.6453 & 379.2330 & 375.6239 & 374.7192 & 469.0618 & 415.1885 & 398.1174 & 392.86 & 378.00 & 357.2698 & 349.4391 & 346.5403\\
36 & 483.1916 & 416.8176 & 400.4065 & 396.3120 & 395.2878 & 481.9094 & 429.2890 & 414.4408 & 409.86 & 390.66 & 369.4654 & 361.4892 & 358.5420\\
37 & 493.4493 & 422.8265 & 405.5962 & 401.2680 & 400.1626 & 496.4015 & 433.4022 & 416.1715 & 411.27 & 394.69 & 373.5203 & 365.5479 & 362.6043\\
38 & 517.9593 & 439.4004 & 419.9247 & 414.9014 & 413.5579 & 504.2032 & 436.6132 & 416.8702 & 411.39 & 394.69 & 373.5203 & 365.5479 & 362.6043\\
39 & 521.0245 & 440.3721 & 420.1850 & 414.9666 & 413.5741 & 515.6273 & 445.0090 & 422.9119 & 416.05 & 404.24 & 380.8049 & 372.0055 & 368.7444\\
40 & 554.2167 & 471.2427 & 450.8932 & 445.8230 & 444.5550 & 536.0207 & 478.3656 & 461.6796 & 455.99 & 412.24 & 390.5909 & 382.4977 & 379.4817\\
41 & 561.4218 & 472.5244 & 451.1818 & 445.8931 & 444.5724 & 552.9508 & 482.3672 & 461.9683 & 456.00 & 425.16 & 402.1807 & 393.5960 & 390.4050\\
42 & 605.5574 & 504.9343 & 479.8486 & 473.4384 & 471.7517 & 558.1539 & 482.9689 & 463.2601 & 458.71 & 425.16 & 402.1807 & 393.5960 & 390.4050\\
43 & 616.5028 & 518.5646 & 493.7396 & 487.4812 & 485.9106 & 567.3412 & 488.3777 & 465.6960 & 459.52 & 434.77 & 411.2604 & 402.4804 & 399.2239\\
44 & 628.0957 & 524.1471 & 498.5358 & 492.2513 & 490.6859 & 574.2506 & 492.2321 & 468.0827 & 460.80 & 454.40 & 428.3600 & 418.6438 & 415.0663\\
45 & 634.2437 & 531.7247 & 505.9813 & 499.5675 & 497.9617 & 598.2522 & 519.5769 & 494.7583 & 487.00 & 466.46 & 438.7119 & 428.4159 & 424.6321\\
46 & 641.5337 & 532.7246 & 506.1940 & 499.6197 & 497.9748 & 608.7209 & 521.6171 & 495.0149 & 487.05 & 466.46 & 438.7119 & 428.4159 & 424.6321\\
47 & 667.4687 & 547.9734 & 517.8253 & 510.1141 & 508.0784 & 630.2738 & 531.2947 & 504.3023 & 496.86 & 473.55 & 445.8369 & 435.5522 & 431.7830\\
48 & 698.1628 & 573.9857 & 543.1162 & 535.4598 & 533.5479 & 664.4216 & 562.5890 & 529.0856 & 518.77 & 503.18 & 473.2766 & 462.2661 & 458.2523\\
49 & 700.7259 & 574.8348 & 543.7180 & 535.8214 & 533.7649 & 670.4290 & 576.8794 & 551.2608 & 543.26 & 514.10 & 483.1208 & 471.5134 & 467.2072\\
50 & 752.4778 & 618.6779 & 584.2646 & 575.3579 & 572.9848 & 683.4854 & 579.9359 & 551.3653 & 543.30 & 514.10 & 483.1208 & 471.6932 & 467.5056\\
51 & 772.0233 & 624.3244 & 585.6444 & 575.6898 & 573.0624 & 701.0202 & 589.9260 & 560.3374 & 553.09 & 514.95 & 483.2307 & 471.6932 & 467.5056\\
52 & 818.5900 & 672.0718 & 631.6665 & 621.5085 & 618.9558 & 803.3161 & 678.8308 & 641.9561 & 631.51 & 585.09 & 547.0977 & 533.1663 & 528.0563\\
53 & 833.2754 & 676.1826 & 632.6191 & 621.7422 & 619.0143 & 805.8113 & 679.9199 & 642.3952 & 631.64 & 590.98 & 552.3832 & 538.3562 & 533.2348\\
54 & 846.5291 & 681.5713 & 641.2196 & 631.1787 & 628.6572 & 832.7409 & 709.9328 & 668.0265 & 655.94 & 592.82 & 554.2646 & 540.1749 & 535.0612\\
55 & 861.5834 & 694.2754 & 650.1452 & 639.0805 & 636.2458 & 848.3769 & 711.2532 & 670.9435 & 659.33 & 592.82 & 554.2646 & 540.1749 & 535.0612\\
56 & 868.5235 & 702.4581 & 659.4855 & 648.8736 & 646.2178 & 866.9036 & 733.0618 & 700.4748 & 690.67 & 610.41 & 569.6585 & 554.8341 & 549.3760\\
57 & 869.8747 & 703.6618 & 659.8472 & 649.0127 & 646.3092 & 881.1987 & 754.4762 & 717.1605 & 706.72 & 610.41 & 569.6585 & 554.8341 & 549.3760\\
58 & 874.0871 & 704.4428 & 659.9817 & 649.0347 & 646.3133 & 885.7117 & 755.2010 & 718.9525 & 708.32 & 611.55 & 570.8047 & 555.9232 & 550.4958\\
59 & 969.8144 & 782.3749 & 728.5923 & 715.1103 & 711.7366 & 901.5093 & 767.8076 & 720.6769 & 708.45 & 611.66 & 572.3659 & 558.0957 & 552.8390\\
60 & 991.5708 & 787.4866 & 733.0238 & 718.9878 & 715.2396 & 916.7462 & 781.9384 & 731.5970 & 718.66 & 614.27 & 573.8453 & 559.0770 & 553.6998\\
    \end{tabular}
    }
    \caption{SC Unnormalized Eigenvalues}
    \label{tablefour3}
\end{table}

\begin{table}[htbp!]
  \resizebox{\tabfourwd}{!}{
  \begin{tabular}{r | r*{21}{r}r}
    \multicolumn{13}{p{21cm}}{SC Normalized Eigenvalue Data}\\
    \hline
    Level: & 1 & 1 & 1 & 1 & 1 & 2 & 2 & 2 & 2 & 3 & 3 & 3 & 3\\
    Refinement: & 0 & 1 & 2 & 3 & 4 & 0 & 1 & 2 & 3 & 0 & 1 & 2 & 3\\
    \hline
    n\\	
1 & 1.0000 & 1.0000 & 1.0000 & 1.0000 & 1.0000 & 1.0000 & 1.0000 & 1.0000 & 1.0000 & 1.0000 & 1.0000 & 1.0000 & 1.0000\\
2 & 1.0008 & 1.0004 & 1.0002 & 1.0001 & 1.0000 & 1.0002 & 1.0001 & 1.0000 & 1.0000 & 1.0000 & 1.0000 & 1.0000 & 1.0000\\
3 & 2.7488 & 2.7718 & 2.7825 & 2.7872 & 2.7892 & 2.6509 & 2.6598 & 2.6639 & 2.6657 & 2.6530 & 2.6536 & 2.6541 & 2.6544\\
4 & 4.9780 & 4.8962 & 4.8706 & 4.8621 & 4.8592 & 4.9927 & 4.9414 & 4.9238 & 4.9181 & 4.9383 & 4.9305 & 4.9276 & 4.9267\\
5 & 6.8070 & 6.7505 & 6.7421 & 6.7423 & 6.7433 & 6.6619 & 6.6126 & 6.6001 & 6.5973 & 6.5950 & 6.5875 & 6.5855 & 6.5850\\
6 & 6.8202 & 6.7537 & 6.7429 & 6.7425 & 6.7434 & 6.6663 & 6.6140 & 6.6005 & 6.5974 & 6.5950 & 6.5875 & 6.5855 & 6.5850\\
7 & 7.5431 & 7.5156 & 7.5217 & 7.5286 & 7.5325 & 6.9800 & 6.9259 & 6.9095 & 6.9052 & 6.8963 & 6.8882 & 6.8861 & 6.8855\\
8 & 13.4392 & 13.2040 & 13.1676 & 13.1679 & 13.1719 & 10.2844 & 10.0879 & 10.0176 & 9.9946 & 10.0497 & 10.0189 & 10.0081 & 10.0043\\
9 & 13.4615 & 13.2102 & 13.1692 & 13.1684 & 13.1720 & 10.3087 & 10.0958 & 10.0198 & 9.9951 & 10.0497 & 10.0189 & 10.0081 & 10.0043\\
10 & 14.3389 & 14.0409 & 13.9892 & 13.9862 & 13.9895 & 11.1947 & 10.9769 & 10.9055 & 10.8833 & 10.9307 & 10.8960 & 10.8843 & 10.8804\\
11 & 14.9572 & 14.5839 & 14.5109 & 14.5012 & 14.5022 & 11.7701 & 11.5070 & 11.4205 & 11.3931 & 11.4684 & 11.4271 & 11.4128 & 11.4080\\
12 & 16.8618 & 16.3774 & 16.2829 & 16.2708 & 16.2724 & 14.0146 & 13.7374 & 13.6643 & 13.6459 & 13.6534 & 13.6184 & 13.6085 & 13.6059\\
13 & 20.0373 & 19.0430 & 18.7787 & 18.7067 & 18.6863 & 16.9380 & 16.4263 & 16.2758 & 16.2307 & 16.3794 & 16.3081 & 16.2843 & 16.2768\\
14 & 20.2061 & 19.4020 & 19.2160 & 19.1753 & 19.1675 & 16.9684 & 16.4361 & 16.2786 & 16.2314 & 16.3794 & 16.3081 & 16.2843 & 16.2768\\
15 & 20.3800 & 19.4430 & 19.2256 & 19.1774 & 19.1680 & 17.6308 & 16.8307 & 16.5727 & 16.4879 & 16.8602 & 16.7600 & 16.7233 & 16.7100\\
16 & 25.5746 & 24.2988 & 24.0110 & 23.9524 & 23.9433 & 25.1532 & 24.4401 & 24.2386 & 24.1877 & 24.1225 & 24.0017 & 23.9679 & 23.9593\\
17 & 25.6387 & 24.3139 & 24.0152 & 23.9536 & 23.9437 & 25.2622 & 24.4748 & 24.2481 & 24.1901 & 24.1225 & 24.0017 & 23.9679 & 23.9593\\
18 & 26.9696 & 25.6051 & 25.3006 & 25.2413 & 25.2337 & 26.2245 & 24.9863 & 24.6624 & 24.5786 & 24.8067 & 24.6561 & 24.6103 & 24.5975\\
19 & 28.5078 & 26.8075 & 26.4351 & 26.3609 & 26.3500 & 27.7032 & 27.6067 & 27.7073 & 27.7851 & 26.7006 & 26.6545 & 26.6585 & 26.6668\\
20 & 30.0292 & 28.2369 & 27.8225 & 27.7340 & 27.7180 & 28.8853 & 28.0460 & 27.8495 & 27.8134 & 27.4371 & 27.3060 & 27.2730 & 27.2664\\
21 & 36.9336 & 34.1558 & 33.4737 & 33.3069 & 33.2668 & 34.9587 & 33.3623 & 32.7467 & 32.5627 & 32.9950 & 32.6827 & 32.5818 & 32.5505\\
22 & 37.3666 & 34.2631 & 33.4993 & 33.3127 & 33.2680 & 35.0949 & 33.9624 & 33.7700 & 33.7589 & 33.0824 & 32.8999 & 32.8558 & 32.8479\\
23 & 40.6694 & 37.2715 & 36.4169 & 36.1993 & 36.1434 & 35.4633 & 33.9686 & 33.7716 & 33.7593 & 33.0824 & 32.8999 & 32.8558 & 32.8479\\
24 & 42.6185 & 39.3834 & 38.6166 & 38.4471 & 38.4153 & 41.5681 & 39.5056 & 38.8817 & 38.7131 & 38.6022 & 38.2700 & 38.1724 & 38.1455\\
25 & 42.8410 & 39.4979 & 38.7205 & 38.5528 & 38.5218 & 41.8128 & 39.5325 & 38.8874 & 38.7146 & 38.6022 & 38.2700 & 38.1724 & 38.1455\\
26 & 43.3004 & 39.5941 & 38.7416 & 38.5578 & 38.5230 & 44.6845 & 42.5801 & 42.1560 & 42.0914 & 41.3849 & 41.0625 & 40.9765 & 40.9574\\
27 & 47.9536 & 43.6855 & 42.6693 & 42.4424 & 42.3974 & 47.9642 & 45.6202 & 45.1710 & 45.1209 & 43.6369 & 43.3836 & 43.3345 & 43.3291\\
28 & 54.2931 & 48.2507 & 46.7404 & 46.3605 & 46.2645 & 50.8997 & 48.0588 & 47.3687 & 47.2201 & 46.5777 & 46.1053 & 45.9684 & 45.9325\\
29 & 59.5542 & 53.3637 & 51.8746 & 51.5312 & 51.4576 & 55.7628 & 50.9722 & 49.2855 & 48.7473 & 50.5131 & 49.6169 & 49.3112 & 49.2122\\
30 & 60.1417 & 53.5338 & 51.9181 & 51.5421 & 51.4603 & 57.0484 & 53.3359 & 52.3928 & 52.1778 & 51.5339 & 50.9505 & 50.7805 & 50.7354\\
31 & 61.2012 & 54.6899 & 53.1406 & 52.7896 & 52.7169 & 57.5916 & 53.3744 & 52.4019 & 52.1804 & 51.5339 & 50.9505 & 50.7805 & 50.7354\\
32 & 61.9583 & 54.8442 & 53.1756 & 52.7982 & 52.7190 & 62.5220 & 58.1513 & 57.1389 & 56.9307 & 55.9042 & 55.3226 & 55.1651 & 55.1271\\
33 & 63.2021 & 55.9673 & 54.2952 & 53.9180 & 53.8391 & 66.2778 & 60.8965 & 59.2794 & 58.8346 & 59.0590 & 58.1714 & 57.8980 & 57.8170\\
34 & 63.6762 & 56.0909 & 54.3252 & 53.9254 & 53.8409 & 67.0246 & 61.0469 & 59.3162 & 58.8439 & 59.0590 & 58.1714 & 57.8980 & 57.8170\\
35 & 65.2918 & 57.8521 & 56.0558 & 55.6437 & 55.5562 & 72.9346 & 66.6958 & 64.8890 & 64.4045 & 64.4013 & 63.4474 & 63.1611 & 63.0807\\
36 & 69.9313 & 61.2576 & 59.1856 & 58.7083 & 58.6057 & 74.9322 & 68.9609 & 67.5495 & 67.1924 & 66.5581 & 65.6132 & 65.3392 & 65.2653\\
37 & 71.4159 & 62.1407 & 59.9527 & 59.4425 & 59.3284 & 77.1856 & 69.6217 & 67.8316 & 67.4236 & 67.2439 & 66.3333 & 66.0728 & 66.0048\\
38 & 74.9632 & 64.5765 & 62.0706 & 61.4621 & 61.3144 & 78.3987 & 70.1375 & 67.9455 & 67.4430 & 67.2439 & 66.3333 & 66.0728 & 66.0048\\
39 & 75.4068 & 64.7193 & 62.1091 & 61.4718 & 61.3168 & 80.1751 & 71.4862 & 68.9302 & 68.2071 & 68.8714 & 67.6269 & 67.2400 & 67.1225\\
40 & 80.2106 & 69.2562 & 66.6482 & 66.0427 & 65.9101 & 83.3460 & 76.8446 & 75.2490 & 74.7537 & 70.2347 & 69.3648 & 69.1365 & 69.0770\\
41 & 81.2534 & 69.4446 & 66.6908 & 66.0531 & 65.9127 & 85.9785 & 77.4874 & 75.2960 & 74.7554 & 72.4361 & 71.4231 & 71.1425 & 71.0653\\
42 & 87.6411 & 74.2077 & 70.9282 & 70.1336 & 69.9423 & 86.7875 & 77.5840 & 75.5066 & 75.2006 & 72.4361 & 71.4231 & 71.1425 & 71.0653\\
43 & 89.2252 & 76.2108 & 72.9815 & 72.2138 & 72.0415 & 88.2161 & 78.4529 & 75.9036 & 75.3333 & 74.0737 & 73.0355 & 72.7483 & 72.6706\\
44 & 90.9030 & 77.0313 & 73.6904 & 72.9205 & 72.7495 & 89.2904 & 79.0721 & 76.2926 & 75.5423 & 77.4183 & 76.0722 & 75.6699 & 75.5544\\
45 & 91.7928 & 78.1449 & 74.7910 & 74.0043 & 73.8282 & 93.0224 & 83.4647 & 80.6405 & 79.8375 & 79.4725 & 77.9106 & 77.4362 & 77.2957\\
46 & 92.8479 & 78.2919 & 74.8224 & 74.0120 & 73.8301 & 94.6502 & 83.7925 & 80.6823 & 79.8453 & 79.4725 & 77.9106 & 77.4362 & 77.2957\\
47 & 96.6014 & 80.5329 & 76.5417 & 75.5666 & 75.3281 & 98.0015 & 85.3471 & 82.1960 & 81.4535 & 80.6806 & 79.1759 & 78.7261 & 78.5974\\
48 & 101.0437 & 84.3558 & 80.2800 & 79.3212 & 79.1042 & 103.3111 & 90.3742 & 86.2354 & 85.0454 & 85.7282 & 84.0489 & 83.5546 & 83.4156\\
49 & 101.4146 & 84.4806 & 80.3690 & 79.3748 & 79.1364 & 104.2452 & 92.6698 & 89.8498 & 89.0616 & 87.5885 & 85.7972 & 85.2261 & 85.0456\\
50 & 108.9046 & 90.9240 & 86.3623 & 85.2316 & 84.9512 & 106.2754 & 93.1608 & 89.8668 & 89.0675 & 87.5885 & 85.7972 & 85.2586 & 85.0999\\
51 & 111.7333 & 91.7538 & 86.5663 & 85.2808 & 84.9627 & 109.0018 & 94.7656 & 91.3292 & 90.6730 & 87.7330 & 85.8167 & 85.2586 & 85.0999\\
52 & 118.4728 & 98.7710 & 93.3690 & 92.0682 & 91.7669 & 124.9079 & 109.0473 & 104.6322 & 103.5291 & 99.6834 & 97.1588 & 96.3698 & 96.1220\\
53 & 120.5982 & 99.3752 & 93.5097 & 92.1028 & 91.7756 & 125.2959 & 109.2222 & 104.7037 & 103.5495 & 100.6863 & 98.0975 & 97.3079 & 97.0646\\
54 & 122.5164 & 100.1671 & 94.7810 & 93.5007 & 93.2052 & 129.4831 & 114.0435 & 108.8814 & 107.5333 & 101.0009 & 98.4316 & 97.6366 & 97.3971\\
55 & 124.6952 & 102.0342 & 96.1004 & 94.6712 & 94.3303 & 131.9144 & 114.2556 & 109.3568 & 108.0895 & 101.0009 & 98.4316 & 97.6366 & 97.3971\\
56 & 125.6996 & 103.2367 & 97.4810 & 96.1220 & 95.8088 & 134.7951 & 117.7589 & 114.1701 & 113.2266 & 103.9967 & 101.1654 & 100.2863 & 100.0028\\
57 & 125.8952 & 103.4137 & 97.5344 & 96.1426 & 95.8223 & 137.0179 & 121.1989 & 116.8897 & 115.8587 & 103.9967 & 101.1654 & 100.2863 & 100.0028\\
58 & 126.5048 & 103.5284 & 97.5543 & 96.1458 & 95.8229 & 137.7196 & 121.3154 & 117.1818 & 116.1207 & 104.1921 & 101.3689 & 100.4831 & 100.2066\\
59 & 140.3592 & 114.9817 & 107.6959 & 105.9340 & 105.5226 & 140.1760 & 123.3405 & 117.4628 & 116.1421 & 104.2096 & 101.6462 & 100.8758 & 100.6332\\
60 & 143.5080 & 115.7330 & 108.3509 & 106.5085 & 106.0420 & 142.5452 & 125.6105 & 119.2427 & 117.8163 & 104.6552 & 101.9089 & 101.0532 & 100.7898\\
    \end{tabular}
    }
    \caption{SC Normalized Eigenvalues}
    \label{tablefour4}
\end{table}

The last two examples we consider are alternate carpets.  We subdivide the unit square into 16 subsquares of side length $\frac{1}{4}$, and retain all but the 4 inner squares ($\frac{12}{16}$ carpet) or all but 3 of the inner squares ($\frac{13}{16}$ carpet).  The $\frac{12}{16}$ carpet has $D_{4}$ symmetry and it is known that a self-similar Laplacian exists.  The $\frac{13}{16}$ carpet has no symmetry, and the methods used to construct a Laplacian on SC do not work on this example.  So the situation is even worse than for the octagasket.

\begin{figure}[htbp!]
\begin{center}
\includegraphics[width=\onewidth\textwidth]{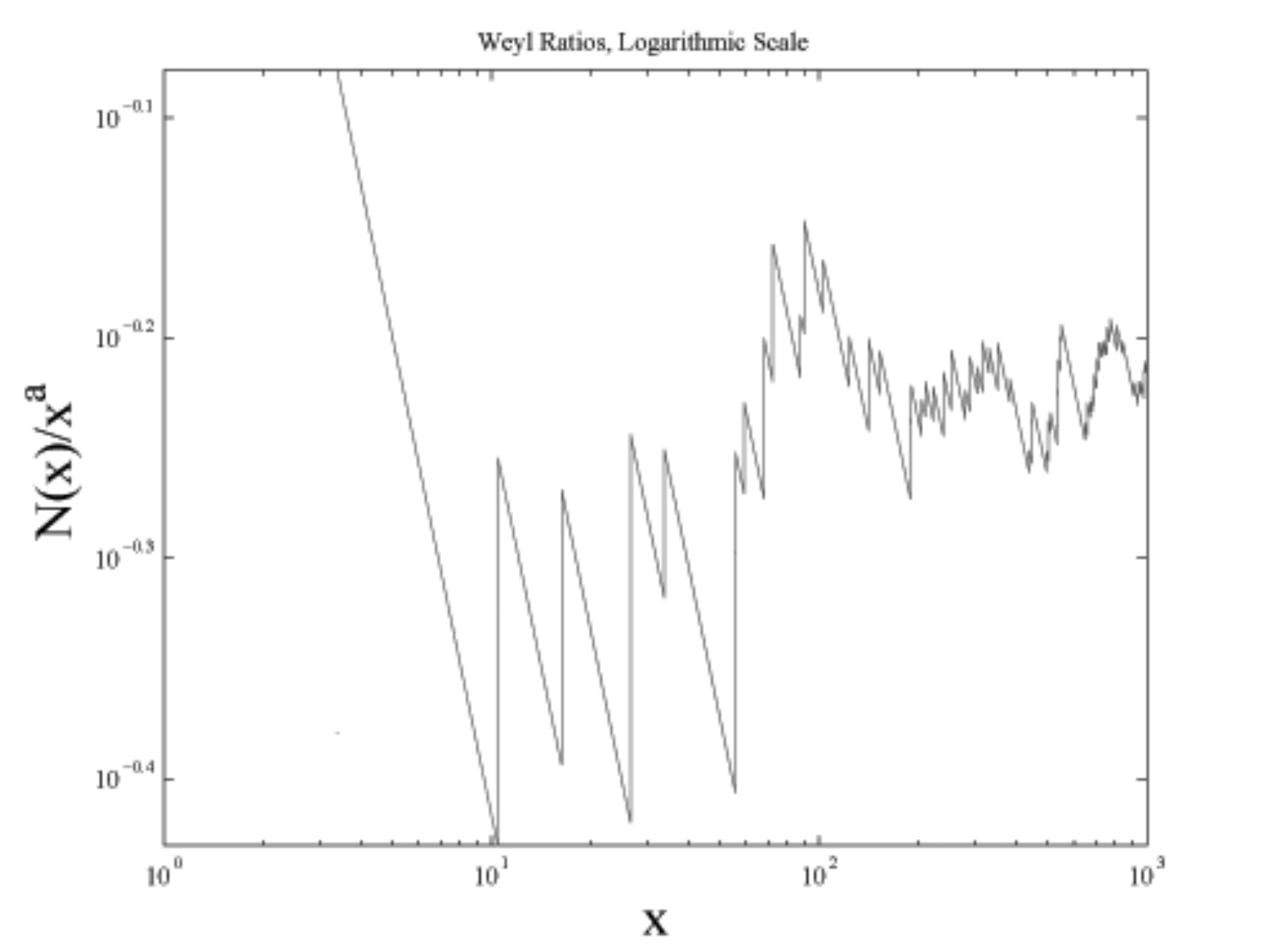}
\end{center}
\caption{$\frac{12}{16}$ Carpet Weyl Ratios, Level 3, 0 Refinements, $\alpha=.71738$}
\label{figfour5}
\end{figure}

\begin{figure}[htbp!]
\begin{center}
\includegraphics[width=\onewidth\textwidth]{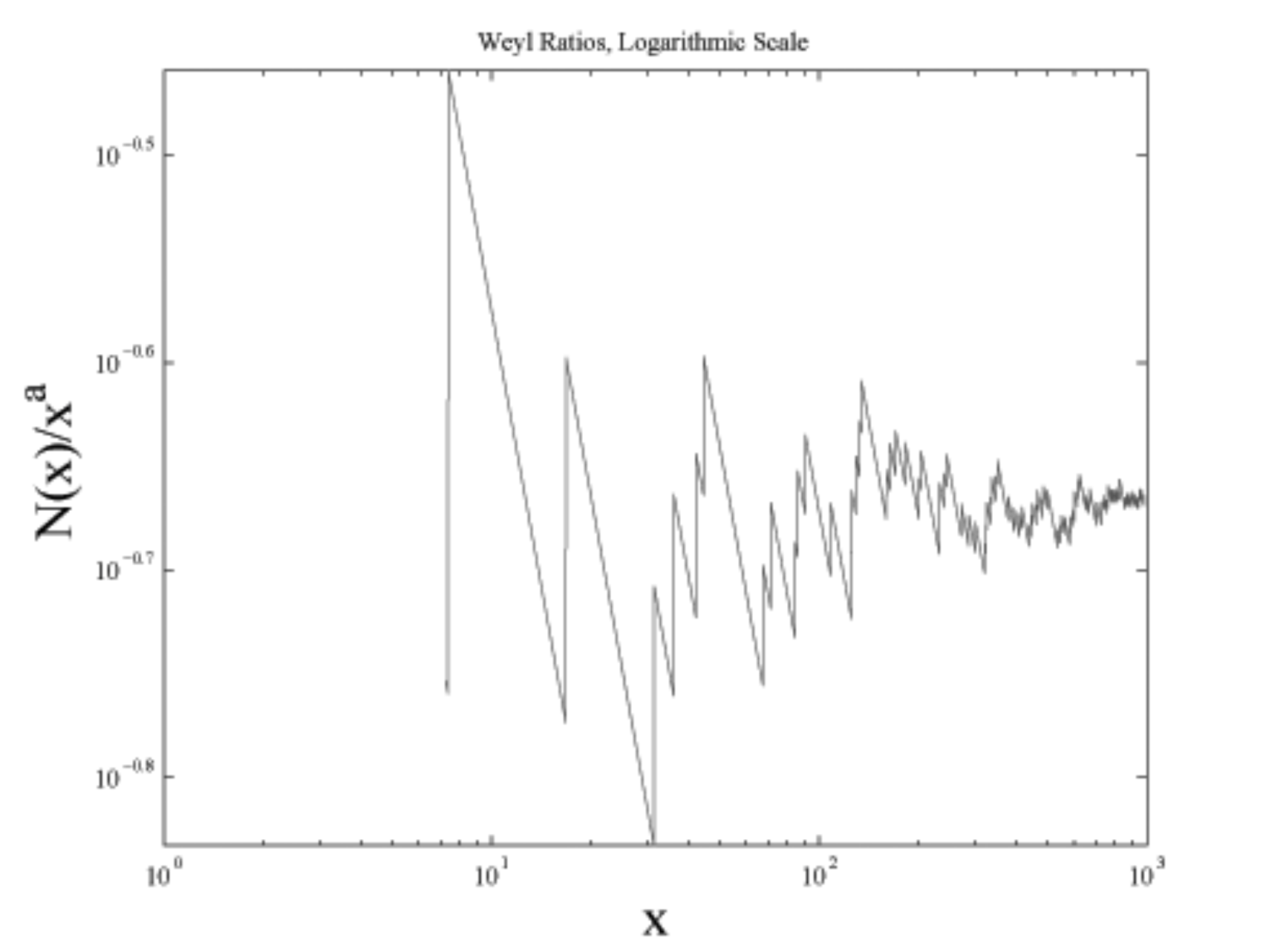}
\end{center}
\caption{$\frac{13}{16}$ Carpet Weyl Ratios, Level 3, 0 Refinements, $\alpha=.87537$}
\label{figfour6}
\end{figure}

In tables \ref{tablefour6} and \ref{tablefour7} we present unnormalized and normalized eigenvalues for the $\frac{12}{16}$ carpet, and in Tables \ref{tablefour8} and \ref{tablefour9} the same data for the $\frac{13}{16}$ carpet.  (Again we use the interior of the square for $\Omega$).  In Figure \ref{figfour5} we show the Weyl ratios for the $\frac{12}{16}$ carpet, and in Figure \ref{figfour6} we show those for the $\frac{13}{16}$ carpet.  The evidence for convergence is strong in both cases.  But the nature of the spectrum is quite different.  In the symmetric $\frac{12}{16}$ carpet, we see multiplicities of 1 or 2, and an eigenvalue renormalization factor of about $R=20.123$.  For the $\frac{13}{16}$ carpet we do not see any multiplicities above 1, and there is no apparent eigenvalue renormalization factor.  The evidence for spectral gaps is also weaker for the $\frac{13}{16}$ carpet, but this is not conclusive.

\begin{table}[htbp!]
  \resizebox{\tabsevenwd}{!}{
  \begin{tabular}{p{1.3cm}p{1.5cm}p{2.1cm}||p{1.3cm}p{1.5cm}p{2.1cm}}
    \multicolumn{6}{p{11cm}}{Sierpinski Carpet, Level 3, 3 Refinements, D4 Representation Type}\\
    \hline
    \mbox{Number} & \mbox{Eigenvalue} & $\substack{\mbox{Eigenfunction} \\ \mbox{Type}}$
    & \mbox{Number} & \mbox{Eigenvalue} & $\substack{\mbox{Eigenfunction} \\ \mbox{Type}}$\\
    \hline
1	& 5.4936	& 2	& 48	& 458.2523	& 1+ -\\
2	& 5.4936	& 2	& 49	& 467.2072	& 1- +\\
3	& 14.5823	& 1+ -	& 50	& 467.5056	& 2\\
4	& 27.0651	& 1- +	& 51	& 467.5056	& 2\\
5	& 36.1754	& 2	& 52	& 528.0563	& 1- +\\
6	& 36.1754	& 2	& 53	& 533.2348	& 1+ +\\
7	& 37.8262	& 1+ +	& 54	& 535.0612	& 2\\
8	& 54.9598	& 2	& 55	& 535.0612	& 2\\
9	& 54.9598	& 2	& 56	& 549.3760	& 2\\
10	& 59.7724	& 1+ -	& 57	& 549.3760	& 2\\
11	& 62.6712	& 1- -	& 58	& 550.4958	& 1+ -\\
12	& 74.7457	& 1+ +	& 59	& 552.8390	& 1+ +\\
13	& 89.4182	& 2	& 60	& 553.6998	& 2\\
14	& 89.4182	& 2	& & & \\
15	& 91.7984	& 1- +	& & & \\
16	& 131.6229	& 2	& & & \\
17	& 131.6229	& 2	& & & \\
18	& 135.1289	& 1- +	& & & \\
19	& 146.4968	& 1+ -	& & & \\
20	& 149.7908	& 1+ +	& & & \\
21	& 178.8197	& 1- -	& & & \\
22	& 180.4537	& 2	& & & \\
23	& 180.4537	& 2	& & & \\
24	& 209.5562	& 2	& & & \\
25	& 209.5562	& 2	& & & \\
26	& 225.0039	& 1+ +	& & & \\
27	& 238.0328	& 1+ -	& & & \\
28	& 252.3353	& 1- -	& & & \\
29	& 270.3526	& 1- +	& & & \\
30	& 278.7202	& 2	& & & \\
31	& 278.7202	& 2	& & & \\
32	& 302.8469	& 1+ -	& & & \\
33	& 317.6237	& 2	& & & \\
34	& 317.6237	& 2	& & & \\
35	& 346.5403	& 1+ +	& & & \\
36	& 358.5420	& 1- +	& & & \\
37	& 362.6043	& 2	& & & \\
38	& 362.6043	& 2	& & & \\
39	& 368.7444	& 1- -	& & & \\
40	& 379.4817	& 1+ +	& & & \\
41	& 390.4050	& 2	& & & \\
42	& 390.4050	& 2	& & & \\
43	& 399.2239	& 1+ -	& & & \\
44	& 415.0663	& 1+ +	& & & \\
45	& 424.6321	& 2	& & & \\
46	& 424.6321	& 2	& & & \\
47	& 431.7830	& 1- -	& & & \\
    \end{tabular}
    }
    \caption{$D_{4}$ Representation Type}
    \label{tablefour5}
\end{table}

\begin{table}[htbp!]
  \resizebox{\tabeightwd}{!}{
  \begin{tabular}{r | r*{6}{r}r}
    \multicolumn{6}{p{9cm}}{12/16 Symmetric Carpet Unnormalized Eigenvalues}\\
    \hline
    Level: & 1 & 1 & 1 & 2 & 2 & 2 & 3 \\
    Refinement: & 0 & 1 & 2 & 0 & 1 & 2 & 0 \\
    \hline
    n\\	
1 & 5.184 & 5.108 & 5.079 & 4.246 & 4.119 & 4.065 & 3.38 \\
2 & 5.194 & 5.113 & 5.081 & 4.246 & 4.119 & 4.065 & 3.38 \\
3 & 16.788 & 16.699 & 16.673 & 13.176 & 12.805 & 12.649 & 10.47\\
4 & 25.029 & 24.128 & 23.805 & 20.638 & 19.903 & 19.593 & 16.37 \\
5 & 43.231 & 42.181 & 41.846 & 33.584 & 32.459 & 31.992 & 26.58 \\
6 & 43.301 & 42.210 & 41.858 & 33.584 & 32.459 & 31.992 & 26.58 \\
7 & 58.008 & 56.942 & 56.648 & 42.636 & 41.255 & 40.689 & 33.70 \\
8 & 93.773 & 89.292 & 87.948 & 70.623 & 68.005 & 66.881 & 55.56 \\
9 & 101.704 & 98.004 & 96.965 & 70.623 & 68.005 & 66.961 & 55.56 \\
10 & 102.074 & 98.095 & 96.985 & 70.890 & 68.023 & 66.961 & 55.62 \\
11 & 107.886 & 104.162 & 103.184 & 75.068 & 72.356 & 71.284 & 59.05 \\
12 & 168.673 & 160.716 & 158.630 & 86.447 & 83.122 & 81.802 & 67.90 \\
13 & 168.747 & 160.725 & 158.632 & 86.472 & 83.124 & 81.803 & 67.90 \\
14 & 175.179 & 166.385 & 164.087 & 92.260 & 88.575 & 87.127 & 72.27 \\
15 & 175.369 & 166.449 & 164.106 & 92.260 & 88.575 & 87.127 & 72.27 \\
16 & 195.384 & 183.514 & 180.284 & 112.216 & 107.801 & 106.110 & 87.76 \\
17 & 195.913 & 183.653 & 180.315 & 115.148 & 110.400 & 108.576 & 90.08 \\
18 & 198.538 & 186.390 & 183.291 & 115.148 & 110.400 & 108.576 & 90.08 \\
19 & 226.387 & 208.148 & 203.133 & 132.107 & 125.817 & 123.405 & 102.89 \\
20 & 250.061 & 230.500 & 225.375 & 159.110 & 152.124 & 149.550 & 123.50 \\
21 & 258.904 & 237.074 & 231.163 & 183.241 & 173.226 & 169.527 & 141.49 \\
22 & 275.083 & 252.049 & 245.674 & 183.241 & 173.226 & 169.527 & 141.49 \\
23 & 276.210 & 252.346 & 245.734 & 199.505 & 186.143 & 181.305 & 153.30 \\
24 & 293.995 & 270.945 & 264.860 & 251.553 & 235.159 & 229.326 & 189.94 \\
25 & 301.675 & 274.485 & 267.363 & 251.553 & 235.159 & 229.326 & 189.94 \\
26 & 303.099 & 274.895 & 267.479 & 253.119 & 235.612 & 229.495 & 191.29 \\
27 & 357.036 & 327.150 & 318.681 & 274.092 & 255.344 & 248.670 & 205.68 \\
28 & 390.429 & 332.905 & 318.781 & 276.411 & 269.536 & 267.090 & 210.72 \\
29 & 411.828 & 361.426 & 347.970 & 299.289 & 275.422 & 267.421 & 223.12 \\
30 & 414.711 & 366.359 & 353.561 & 318.345 & 301.347 & 294.251 & 239.50 \\
31 & 416.004 & 366.831 & 353.679 & 318.345 & 301.347 & 294.251 & 239.50 \\
32 & 434.933 & 389.600 & 377.634 & 341.808 & 315.583 & 308.149 & 253.11 \\
33 & 468.215 & 408.316 & 392.905 & 341.808 & 315.583 & 308.149 & 253.11 \\
34 & 476.368 & 409.745 & 393.225 & 382.105 & 347.140 & 334.977 & 279.87 \\
35 & 540.090 & 462.398 & 442.531 & 387.080 & 356.422 & 342.840 & 288.43 \\
36 & 588.083 & 490.445 & 465.059 & 396.599 & 366.152 & 359.525 & 288.50 \\
37 & 598.813 & 516.054 & 493.601 & 410.635 & 385.344 & 373.006 & 304.18 \\
38 & 605.365 & 517.742 & 493.954 & 437.436 & 389.213 & 373.006 & 315.47 \\
39 & 619.523 & 521.657 & 496.133 & 437.436 & 389.213 & 377.237 & 315.47 \\
40 & 621.643 & 532.628 & 510.347 & 461.115 & 405.563 & 387.438 & 329.91 \\
41 & 677.697 & 568.112 & 538.900 & 479.305 & 446.788 & 436.606 & 350.60 \\
42 & 685.080 & 569.402 & 539.138 & 479.305 & 446.788 & 436.606 & 350.60 \\
43 & 744.709 & 619.345 & 585.659 & 525.984 & 488.838 & 477.254 & 381.79 \\
44 & 806.313 & 678.293 & 641.135 & 611.631 & 562.499 & 546.976 & 437.50 \\
45 & 817.160 & 681.832 & 643.705 & 622.996 & 575.895 & 561.102 & 445.75 \\
46 & 825.398 & 684.435 & 644.671 & 623.144 & 575.895 & 561.102 & 445.75 \\
47 & 830.539 & 688.865 & 650.019 & 623.144 & 580.403 & 566.833 & 447.11 \\
48 & 842.319 & 689.219 & 650.609 & 688.848 & 646.061 & 633.294 & 494.74 \\
49 & 856.140 & 694.608 & 653.557 & 699.940 & 653.837 & 639.759 & 501.09 \\
50 & 864.861 & 696.009 & 654.689 & 699.940 & 653.837 & 639.759 & 501.09 \\
51 & 882.370 & 724.078 & 681.934 & 716.995 & 663.525 & 646.620 & 509.19 \\
52 & 899.459 & 729.338 & 683.412 & 752.776 & 692.509 & 673.230 & 532.66 \\
53 & 913.961 & 730.414 & 685.271 & 759.692 & 694.576 & 673.230 & 532.94 \\
54 & 949.231 & 751.437 & 701.637 & 759.692 & 694.576 & 674.760 & 532.94 \\
55 & 957.747 & 754.911 & 702.461 & 766.267 & 699.657 & 677.267 & 536.83 \\
56 &          & 802.767 & 748.028 & 771.554 & 700.249 & 677.322 & 536.83 \\
57 &          & 804.804 & 748.705 & 782.463 & 709.971 & 687.690 & 544.75 \\
58 &          & 821.999 & 758.507 & 782.463 & 709.971 & 687.690 & 544.75 \\
59 &          & 839.997 & 780.140 & 789.713 & 715.217 & 692.308 & 548.14 \\
60 &          & 870.863 & 810.228 & 951.336 & 854.519 & 825.552 & 648.03 \\
    \end{tabular}
    }
    \caption{12/16 Carpet Unnormalized Eigenvalues}
    \label{tablefour6}
\end{table}
\begin{table}[htbp!]
  \resizebox{\tabeightwd}{!}{
  \begin{tabular}{r | r*{6}{r}r}
    \multicolumn{6}{p{9cm}}{12/16 Symmetric Carpet Normalized Eigenvalues}\\
    \hline
    Level: & 1 & 1 & 1 & 2 & 2 & 2 & 3 \\
    Refinement: & 0 & 1 & 2 & 0 & 1 & 2 & 0 \\
    \hline
    n\\	
1 & 1.000 & 1.000 & 1.000 & 1.000 & 1.000 & 1.000 & 1.000\\
2 & 1.002 & 1.001 & 1.000 & 1.000 & 1.000 & 1.000 & 1.000\\
3 & 3.239 & 3.269 & 3.283 & 3.103 & 3.109 & 3.112 & 3.097\\
4 & 4.829 & 4.724 & 4.687 & 4.861 & 4.831 & 4.820 & 4.841\\
5 & 8.340 & 8.258 & 8.239 & 7.910 & 7.879 & 7.870 & 7.859\\
6 & 8.354 & 8.263 & 8.242 & 7.910 & 7.879 & 7.870 & 7.859\\
7 & 11.191 & 11.147 & 11.154 & 10.042 & 10.015 & 10.009 & 9.967\\
8 & 18.091 & 17.481 & 17.317 & 16.633 & 16.508 & 16.453 & 16.431\\
9 & 19.621 & 19.186 & 19.092 & 16.633 & 16.508 & 16.472 & 16.431\\
10 & 19.692 & 19.204 & 19.096 & 16.696 & 16.513 & 16.472 & 16.449\\
11 & 20.813 & 20.392 & 20.316 & 17.680 & 17.565 & 17.536 & 17.463\\
12 & 32.540 & 31.463 & 31.233 & 20.360 & 20.178 & 20.123 & 20.081\\
13 & 32.555 & 31.465 & 31.234 & 20.366 & 20.178 & 20.124 & 20.081\\
14 & 33.795 & 32.573 & 32.308 & 21.729 & 21.502 & 21.433 & 21.373\\
15 & 33.832 & 32.586 & 32.312 & 21.729 & 21.502 & 21.433 & 21.373\\
16 & 37.693 & 35.927 & 35.497 & 26.429 & 26.169 & 26.103 & 25.953\\
17 & 37.795 & 35.954 & 35.503 & 27.120 & 26.800 & 26.710 & 26.640\\
18 & 38.302 & 36.490 & 36.089 & 27.120 & 26.800 & 26.710 & 26.640\\
19 & 43.674 & 40.749 & 39.996 & 31.114 & 30.542 & 30.358 & 30.428\\
20 & 48.242 & 45.125 & 44.375 & 37.474 & 36.928 & 36.789 & 36.523\\
21 & 49.948 & 46.412 & 45.515 & 43.157 & 42.051 & 41.704 & 41.843\\
22 & 53.069 & 49.344 & 48.372 & 43.157 & 42.051 & 41.704 & 41.843\\
23 & 53.286 & 49.402 & 48.384 & 46.988 & 45.186 & 44.601 & 45.338\\
24 & 56.717 & 53.043 & 52.150 & 59.246 & 57.085 & 56.414 & 56.173\\
25 & 58.199 & 53.736 & 52.642 & 59.246 & 57.085 & 56.414 & 56.173\\
26 & 58.474 & 53.816 & 52.665 & 59.615 & 57.195 & 56.456 & 56.572\\
27 & 68.879 & 64.046 & 62.747 & 64.554 & 61.985 & 61.173 & 60.827\\
28 & 75.321 & 65.173 & 62.766 & 65.101 & 65.430 & 65.704 & 62.317\\
29 & 79.450 & 70.756 & 68.514 & 70.489 & 66.859 & 65.786 & 65.987\\
30 & 80.006 & 71.722 & 69.614 & 74.977 & 73.152 & 72.386 & 70.830\\
31 & 80.255 & 71.815 & 69.638 & 74.977 & 73.152 & 72.386 & 70.830\\
32 & 83.907 & 76.272 & 74.354 & 80.503 & 76.608 & 75.805 & 74.855\\
33 & 90.328 & 79.936 & 77.361 & 80.503 & 76.608 & 75.805 & 74.855\\
34 & 91.901 & 80.216 & 77.424 & 89.993 & 84.268 & 82.405 & 82.767\\
35 & 104.194 & 90.524 & 87.132 & 91.165 & 86.521 & 84.339 & 85.300\\
36 & 113.453 & 96.015 & 91.568 & 93.407 & 88.884 & 88.443 & 85.321\\
37 & 115.523 & 101.028 & 97.188 & 96.713 & 93.542 & 91.760 & 89.958\\
38 & 116.787 & 101.359 & 97.257 & 103.025 & 94.482 & 91.760 & 93.297\\
39 & 119.518 & 102.125 & 97.686 & 103.025 & 94.482 & 92.801 & 93.297\\
40 & 119.927 & 104.273 & 100.485 & 108.602 & 98.451 & 95.310 & 97.568\\
41 & 130.741 & 111.219 & 106.107 & 112.886 & 108.458 & 107.405 & 103.688\\
42 & 132.165 & 111.472 & 106.154 & 112.886 & 108.458 & 107.405 & 103.688\\
43 & 143.669 & 121.249 & 115.313 & 123.880 & 118.666 & 117.405 & 112.911\\
44 & 155.554 & 132.790 & 126.236 & 144.052 & 136.547 & 134.556 & 129.385\\
45 & 157.646 & 133.482 & 126.742 & 146.728 & 139.799 & 138.032 & 131.825\\
46 & 159.236 & 133.992 & 126.933 & 146.763 & 139.799 & 138.032 & 131.825\\
47 & 160.227 & 134.859 & 127.986 & 146.763 & 140.893 & 139.441 & 132.228\\
48 & 162.500 & 134.929 & 128.102 & 162.238 & 156.832 & 155.791 & 146.313\\
49 & 165.166 & 135.984 & 128.682 & 164.850 & 158.719 & 157.381 & 148.192\\
50 & 166.849 & 136.258 & 128.905 & 164.850 & 158.719 & 157.381 & 148.192\\
51 & 170.226 & 141.753 & 134.269 & 168.867 & 161.071 & 159.069 & 150.587\\
52 & 173.523 & 142.783 & 134.560 & 177.294 & 168.107 & 165.615 & 157.527\\
53 & 176.321 & 142.993 & 134.927 & 178.923 & 168.609 & 165.615 & 157.610\\
54 & 183.125 & 147.109 & 138.149 & 178.923 & 168.609 & 165.991 & 157.610\\
55 & 184.768 & 147.789 & 138.311 & 180.472 & 169.842 & 166.608 & 158.761\\
56 &          & 157.158 & 147.283 & 181.717 & 169.986 & 166.622 & 158.761\\
57 &          & 157.557 & 147.416 & 184.286 & 172.346 & 169.172 & 161.104\\
58 &          & 160.923 & 149.346 & 184.286 & 172.346 & 169.172 & 161.104\\
59 &          & 164.446 & 153.606 & 185.994 & 173.619 & 170.308 & 162.107\\
60 &          & 170.489 & 159.530 & 224.059 & 207.435 & 203.086 & 191.647\\
    \end{tabular}
    }
    \caption{12/16 Carpet Normalized Eigenvalues}
    \label{tablefour7}
\end{table}

\begin{table}[htbp!]
  \resizebox{\tabeightwd}{!}{
  \begin{tabular}{r | r*{6}{r}r}
    \multicolumn{6}{p{9cm}}{13/16 Alternate Carpet Unnormalized Eigenvalues}\\
    \hline
    Level: & 1 & 1 & 1 & 2 & 2 & 2 & 3 \\
    Refinement: & 0 & 1 & 2 & 0 & 1 & 2 & 0 \\
    \hline
    n\\	
1 & 8.141 & 8.025 & 7.981 & 7.777 & 7.590 & 7.512 & 7.261\\
2 & 8.328 & 8.201 & 8.151 & 7.920 & 7.724 & 7.643 & 7.375\\
3 & 18.904 & 18.728 & 18.674 & 18.046 & 17.656 & 17.498 & 16.889\\
4 & 36.040 & 35.189 & 34.898 & 33.783 & 32.859 & 32.488 & 31.352\\
5 & 41.026 & 40.392 & 40.209 & 38.618 & 37.681 & 37.307 & 36.009\\
6 & 47.908 & 46.956 & 46.686 & 45.136 & 44.028 & 43.592 & 42.108\\
7 & 51.147 & 50.262 & 50.022 & 47.827 & 46.543 & 46.037 & 44.385\\
8 & 78.834 & 76.085 & 75.308 & 72.746 & 70.651 & 69.849 & 67.394\\
9 & 83.837 & 79.862 & 78.600 & 77.387 & 74.736 & 73.703 & 71.197\\
10 & 98.313	& 95.303 & 94.500 & 91.043 & 88.502 & 87.534 & 84.498\\
11 & 101.603 & 97.886 & 96.879 & 93.070 & 90.247 & 89.174 & 86.066\\
12 & 110.047 & 105.754 & 104.599 & 98.936 & 95.459 & 94.134 & 90.796\\
13 & 130.410 & 123.766 & 121.928 & 118.044 & 113.854 & 112.296 & 108.495\\
14 & 158.249 & 149.670 & 147.378 & 137.776 & 132.030 & 129.857 & 125.544\\
15 & 168.441 & 160.642 & 158.609 & 138.092 & 132.224 & 130.016 & 125.671\\
16 & 168.531 & 160.653 & 158.611 & 140.785 & 135.943 & 134.176 & 129.506\\
17 & 177.769 & 168.614 & 166.217 & 145.902 & 139.753 & 137.470 & 132.850\\
18 & 180.457 & 170.836 & 168.310 & 147.908 & 141.626 & 139.298 & 134.689\\
19 & 212.057 & 199.474 & 196.150 & 176.061 & 168.656 & 165.974 & 160.416\\
20 & 218.578 & 203.366 & 199.309 & 180.743 & 173.034 & 170.242 & 164.492\\
21 & 224.839 & 208.408 & 204.094 & 187.949 & 180.161 & 177.368 & 171.216\\
22 & 244.328 & 220.798 & 214.229 & 202.900 & 192.952 & 189.318 & 183.189\\
23 & 276.284 & 252.641 & 245.740 & 223.491 & 212.361 & 208.388 & 201.776\\
24 & 279.689 & 253.095 & 246.152 & 227.445 & 216.001 & 211.926 & 204.606\\
25 & 282.172 & 260.464 & 254.842 & 254.339 & 244.861 & 241.567 & 232.848\\
26 & 284.558 & 261.350 & 255.511 & 254.348 & 245.378 & 242.312 & 233.531\\
27 & 296.378 & 271.155 & 264.532 & 266.098 & 253.406 & 249.096 & 240.709\\
28 & 315.919 & 285.291 & 277.333 & 269.589 & 258.570 & 254.766 & 244.795\\
29 & 345.168 & 315.927 & 308.136 & 306.237 & 287.886 & 281.712 & 271.460\\
30 & 356.577 & 326.898 & 318.697 & 309.963 & 300.347 & 297.247 & 286.438\\
31 & 387.016 & 335.984 & 322.765 & 327.873 & 315.925 & 311.747 & 299.839\\
32 & 412.864 & 370.586 & 359.404 & 355.748 & 337.374 & 331.115 & 319.702\\
33 & 419.313 & 374.241 & 362.670 & 356.922 & 341.061 & 335.720 & 323.243\\
34 & 429.690 & 377.769 & 364.163 & 361.714 & 344.265 & 338.169 & 325.064\\
35 & 431.820 & 383.436 & 371.330 & 372.504 & 353.870 & 347.561 & 333.930\\
36 & 449.266 & 398.314 & 384.230 & 379.569 & 356.686 & 349.061 & 336.157\\
37 & 473.522 & 412.916 & 397.013 & 385.570 & 365.766 & 359.457 & 346.095\\
38 & 476.703 & 413.305 & 397.076 & 385.668 & 368.320 & 363.064 & 349.714\\
39 & 493.685 & 424.342 & 407.726 & 419.014 & 399.057 & 392.766 & 377.470\\
40 & 519.345 & 450.491 & 432.843 & 429.524 & 414.512 & 410.168 & 394.171\\
41 & 543.384 & 466.699 & 447.185 & 456.574 & 433.048 & 423.859 & 408.458\\
42 & 585.649 & 501.099 & 478.186 & 467.753 & 439.744 & 432.360 & 416.050\\
43 & 591.644 & 514.912 & 489.882 & 491.355 & 464.059 & 455.055 & 437.295\\
44 & 601.176 & 517.679 & 496.841 & 498.651 & 470.180 & 461.081 & 443.572\\
45 & 620.159 & 520.471 & 497.666 & 506.206 & 476.228 & 466.560 & 448.465\\
46 & 622.937 & 532.290 & 509.206 & 509.421 & 485.929 & 478.706 & 460.467\\
47 & 651.834 & 544.330 & 515.901 & 513.937 & 489.455 & 481.930 & 463.403\\
48 & 660.203 & 555.040 & 530.133 & 539.091 & 507.854 & 498.336 & 479.947\\
49 & 675.397 & 567.572 & 538.424 & 541.151 & 512.129 & 503.378 & 483.585\\
50 & 694.404 & 583.198 & 554.867 & 557.873 & 526.650 & 516.944 & 497.031\\
51 & 753.195 & 627.341 & 589.310 & 604.534 & 567.907 & 556.668 & 534.699\\
52 & 771.635 & 629.002 & 595.643 & 617.337 & 576.543 & 563.962 & 543.660\\
53 & 804.767 & 672.510 & 629.396 & 631.441 & 587.316 & 572.032 & 548.612\\
54 & 811.516 & 677.324 & 643.350 & 637.710 & 596.322 & 583.937 & 561.722\\
55 & 836.030 & 680.492 & 644.141 & 660.452 & 620.202 & 606.507 & 581.113\\
56 & 844.566 & 691.307 & 652.911 & 663.900 & 621.142 & 608.571 & 585.293\\
57 & 850.130 & 694.214 & 653.524 & 669.942 & 629.911 & 617.199 & 591.844\\
58 & 866.855 & 704.415 & 661.306 & 685.466 & 639.886 & 625.486 & 600.271\\
59 & 873.875 & 716.879 & 674.548 & 686.218 & 640.602 & 626.171 & 601.676\\
60 & 887.151 & 725.353 & 682.655 & 687.249 & 652.785 & 638.914 & 613.170\\
    \end{tabular}
    }
    \caption{13/16 Carpet Unnormalized Eigenvalues}
    \label{tablefour8}
\end{table}
\begin{table}[htbp!]
  \resizebox{\tabeightwd}{!}{
  \begin{tabular}{r | r*{6}{r}r}
    \multicolumn{6}{p{9cm}}{13/16 Alternate Carpet Normalized Eigenvalues}\\
    \hline
    Level: & 1 & 1 & 1 & 2 & 2 & 2 & 3 \\
    Refinement: & 0 & 1 & 2 & 0 & 1 & 2 & 0 \\
    \hline
    n\\
1 & 1.000 & 1.000 & 1.000 & 1.000 & 1.000 & 1.000 & 1.000\\
2 & 1.023 & 1.022 & 1.021 & 1.018 & 1.018 & 1.017 & 1.016\\
3 & 2.322 & 2.334 & 2.340 & 2.320 & 2.326 & 2.329 & 2.326\\
4 & 4.427 & 4.385 & 4.373 & 4.344 & 4.329 & 4.325 & 4.318\\
5 & 5.040 & 5.033 & 5.038 & 4.966 & 4.965 & 4.966 & 4.960\\
6 & 5.885 & 5.851 & 5.850 & 5.804 & 5.801 & 5.803 & 5.799\\
7 & 6.283 & 6.263 & 6.268 & 6.150 & 6.132 & 6.128 & 6.113\\
8 & 9.684 & 9.481 & 9.436 & 9.354 & 9.309 & 9.298 & 9.282\\
9 & 10.299 & 9.952 & 9.849 & 9.951 & 9.847 & 9.811 & 9.806\\
10 & 12.077 & 11.876 & 11.841 & 11.706 & 11.661 & 11.652 & 11.638\\
11 & 12.481 & 12.198 & 12.139 & 11.967 & 11.891 & 11.870 & 11.854\\
12 & 13.518 & 13.178 & 13.107 & 12.721 & 12.577 & 12.531 & 12.505\\
13 & 16.020 & 15.423 & 15.278 & 15.178 & 15.001 & 14.948 & 14.943\\
14 & 19.440 & 18.651 & 18.467 & 17.715 & 17.396 & 17.286 & 17.291\\
15 & 20.692 & 20.018 & 19.874 & 17.756 & 17.421 & 17.307 & 17.309\\
16 & 20.703 & 20.020 & 19.874 & 18.102 & 17.911 & 17.861 & 17.837\\
17 & 21.837 & 21.012 & 20.827 & 18.760 & 18.413 & 18.299 & 18.297\\
18 & 22.168 & 21.289 & 21.090 & 19.018 & 18.660 & 18.543 & 18.551\\
19 & 26.049 & 24.857 & 24.578 & 22.638 & 22.222 & 22.094 & 22.094\\
20 & 26.851 & 25.342 & 24.974 & 23.240 & 22.798 & 22.662 & 22.655\\
21 & 27.620 & 25.970 & 25.574 & 24.167 & 23.738 & 23.610 & 23.581\\
22 & 30.014 & 27.514 & 26.843 & 26.089 & 25.423 & 25.201 & 25.230\\
23 & 33.939 & 31.482 & 30.792 & 28.737 & 27.980 & 27.740 & 27.790\\
24 & 34.358 & 31.539 & 30.844 & 29.245 & 28.460 & 28.211 & 28.180\\
25 & 34.663 & 32.457 & 31.932 & 32.703 & 32.262 & 32.156 & 32.070\\
26 & 34.956 & 32.568 & 32.016 & 32.704 & 32.330 & 32.255 & 32.164\\
27 & 36.408 & 33.790 & 33.147 & 34.215 & 33.388 & 33.158 & 33.153\\
28 & 38.808 & 35.551 & 34.751 & 34.664 & 34.068 & 33.913 & 33.715\\
29 & 42.401 & 39.369 & 38.610 & 39.376 & 37.931 & 37.500 & 37.388\\
30 & 43.803 & 40.736 & 39.934 & 39.855 & 39.573 & 39.568 & 39.451\\
31 & 47.542 & 41.868 & 40.443 & 42.158 & 41.625 & 41.498 & 41.297\\
32 & 50.717 & 46.180 & 45.034 & 45.742 & 44.451 & 44.076 & 44.032\\
33 & 51.509 & 46.635 & 45.444 & 45.893 & 44.937 & 44.689 & 44.520\\
34 & 52.784 & 47.075 & 45.631 & 46.510 & 45.359 & 45.015 & 44.771\\
35 & 53.046 & 47.781 & 46.529 & 47.897 & 46.625 & 46.266 & 45.992\\
36 & 55.189 & 49.635 & 48.145 & 48.805 & 46.996 & 46.465 & 46.299\\
37 & 58.168 & 51.455 & 49.747 & 49.577 & 48.192 & 47.849 & 47.667\\
38 & 58.559 & 51.503 & 49.755 & 49.589 & 48.529 & 48.329 & 48.166\\
39 & 60.645 & 52.879 & 51.089 & 53.877 & 52.579 & 52.283 & 51.989\\
40 & 63.797 & 56.137 & 54.236 & 55.229 & 54.615 & 54.600 & 54.289\\
41 & 66.750 & 58.157 & 56.034 & 58.707 & 57.057 & 56.422 & 56.256\\
42 & 71.942 & 62.444 & 59.918 & 60.144 & 57.939 & 57.554 & 57.302\\
43 & 72.679 & 64.165 & 61.384 & 63.179 & 61.143 & 60.575 & 60.228\\
44 & 73.850 & 64.510 & 62.256 & 64.117 & 61.949 & 61.377 & 61.093\\
45 & 76.182 & 64.858 & 62.359 & 65.088 & 62.746 & 62.106 & 61.767\\
46 & 76.523 & 66.331 & 63.805 & 65.502 & 64.025 & 63.723 & 63.420\\
47 & 80.073 & 67.831 & 64.644 & 66.082 & 64.489 & 64.152 & 63.824\\
48 & 81.101 & 69.165 & 66.427 & 69.317 & 66.913 & 66.336 & 66.103\\
49 & 82.967 & 70.727 & 67.466 & 69.582 & 67.477 & 67.007 & 66.604\\
50 & 85.302 & 72.674 & 69.526 & 71.732 & 69.390 & 68.813 & 68.456\\
51 & 92.524 & 78.175 & 73.842 & 77.732 & 74.826 & 74.101 & 73.644\\
52 & 94.789 & 78.382 & 74.636 & 79.378 & 75.964 & 75.072 & 74.878\\
53 & 98.859 & 83.804 & 78.865 & 81.191 & 77.383 & 76.146 & 75.560\\
54 & 99.688 & 84.404 & 80.614 & 81.997 & 78.570 & 77.731 & 77.365\\
55 & 102.700 & 84.798 & 80.713 & 84.922 & 81.716 & 80.735 & 80.036\\
56 & 103.748 & 86.146 & 81.812 & 85.365 & 81.840 & 81.010 & 80.612\\
57 & 104.432 & 86.508 & 81.888 & 86.142 & 82.995 & 82.159 & 81.514\\
58 & 106.486 & 87.780 & 82.863 & 88.138 & 84.309 & 83.262 & 82.675\\
59 & 107.349 & 89.333 & 84.523 & 88.235 & 84.404 & 83.353 & 82.868\\
60 & 108.979 & 90.389 & 85.539 & 88.367 & 86.009 & 85.049 & 84.451\\
    \end{tabular}
    }
    \caption{13/16 Carpet Normalized Eigenvalues}
    \label{tablefour9}
\end{table}

\clearpage
\section{Miniaturization}
\label{secmini}

In order to make the ideas clear, we begin by explaining the method of miniaturization on the unit interval $I$. Here we have a two element group of symmetries consisting of the identity and the reflection $\rho(x)=1-x$ about the midpoint.  Every Neumann eigenfunction is of the form $\cos \pi kx$.  When $k$ is even, the function is even under $\rho$, namely $u\circ\rho=u$, while if $k$ is odd then the function is odd, namely $u\circ\rho=-u$.  In this way all eigenspaces are sorted corresponding to the two irreducible representations of the symmetry group.  For every even eigenfunction $u$ (except the constant), we can miniaturize it by defining $u_{+}$ to be
\begin{equation}\label{five1}
u_{+}(x) =
\left\{
\begin{aligned}
u\circ F_{0}^{-1}\mbox{ on } F_{0}I\\
u\circ F_{1}^{-1}\mbox{ on } F_{1}I\\
\end{aligned}
\right.
\end{equation}
Note that $u\circ F_{0}^{-1}(\frac{1}{2})=u\circ F_{1}^{-1}(\frac{1}{2})$ because $u$ is even, and the derivative vanishes at $\frac{1}{2}$ because $u$ is a Neumann eigenfunction.  This shows that $u_{+}$ is also a Neumann eigenfunction, and indeed $u_{+}(x)=\cos 2\pi kx$.  On the other hand, if $u$ is an odd eigenfunction, then define $u_{-}$ by
\begin{equation}\label{five2}
u_{+}(x) =
\left\{
\begin{aligned}
\quad u\circ F_{0}^{-1}  & \mbox{ on } F_{0}I\\
-u\circ F_{1}^{-1} & \mbox{ on } F_{1}I\\
\end{aligned}
\right.
\end{equation}
Again $u\circ F_{0}^{-1}(\frac{1}{2})=-u\circ F_{1}^{-1}(\frac{1}{2})$ because $u$ is odd, so $u_{-}$ is also a Neumann eigenfunction, and again $u_{-}(x)=\cos 2\pi kx$.  We call $u_{+}$ or $u_{-}$ the miniaturization of $u$. Note that the representation type of the miniaturization is always even.  The eigenvalue of $u_{+}$ or $u_{-}$ is always 4 times the eigenvalue of $u$.  Thus $R=4$ is an eigenvalue renormalization factor.  (Of course $I$ has other eigenvalue renormalization factors, namely any square integer, but such luxuries do not generalize to other fractals).

Now consider a self-similar fractal with a finite group of symmetries $G$, and suppose the Laplacian is $G$ invariant.  Then each eigenspace splits according to the irreducible representations of $G$.  We seek to find a set of recipes, analogous to (\ref{five1}) and (\ref{five2}), to miniaturize eigenfunctions according to the corresponding irreducible representations of $G$.  In fact, our goal is to obtain recipes that make sense on the fractal and also on the outer approximating domains.  In the latter case the miniaturization of an eigenfunction on $\Omega_{m}$ will be an eigenfunction on $\Omega_{m+1}$.

It is by no means clear that this goal is always attainable.  We will show explicitly that it is possible for SC, the $\frac{12}{16}$ carpet, and the octagasket.  In the first two examples the symmetry group is $D_{4}$ (the dihedral symmetry group of the square), and in the last example it is $D_{8}$.  In contrast to the interval, the representation type of the miniaturized eigenfunctions is the same as the original one.

The referee has pointed out that it is also possible to explain miniaturization on carpets using local reflection maps introduced in \cite{barlow89} and \cite{barlow99} (see also Definition $2.12$ in \cite{barlow08}).

We mention in passing that a version of miniaturization is valid for SG, but the recipes are more complicated.  In particular, the multiplicities increase.  This is part of the story of spectral decimation (see \cite{strichartz06} for a description).  On the other hand, it is not clear how to extend the recipes for the approximating domains $\Omega_{m}$ with a positive $\epsilon$ overlap, although they are presumably valid in the zero overlap case.

The symmetry group $D_{4}$ has five irreducible representations.  Let $\rho_{H}$ and $\rho_{V}$ denote the reflections about the horizontal and vertical axes in $D_{4}$, and $\rho'_{D}$ and $\rho''_{D}$ denote the two diagonal reflections.  The four one-dimensional representations $1++,1+-,1-+$, and $1--$ are characterized by parity with respect to these reflections.  (Strictly speaking, we describe functions that transform according to the representations, rather than the abstract representations, since we are interested in eigenfunctions that transform according to representations).  Functions transforming according to $1++$ are even with respect to all reflections, and those transforming according to $1--$ are odd with respect to all reflections.  The $1+-$ functions are odd with respect to $\rho_{H}$ and $\rho_{V}$ and even with respect to $\rho'_{D}$ and $\rho''_{D}$, while for $1-+$ the reverse holds.

Now suppose $u$ is a Neumann eigenfunction on $\Omega_{m}$ of $1++$ or $1-+$ type.  Define the miniaturization
\begin{equation}\label{five3}
u_{+}=\{u\circ F_{i}^{-1}\mbox{ on } F_{i}\Omega_{m}\}\mbox{ on }\Omega_{m+1}
\end{equation}
for either the SC or $\frac{12}{16}$ carpet.  On the other hand, for an eigenfunction of $1+-$ or $1--$ type define
\begin{equation}\label{five4}
u_{-}=\{\pm u\circ F_{i}^{-1}\mbox{ on } F_{i}\Omega_{m}\}\mbox{ on }\Omega_{m+1}
\end{equation}
where we alternate the choice of $\pm$ on neighboring cells (see Figure\ref{figfive1}).  Because of the even or odd parity of $u$ with respect to the reflections $\rho_{H}$ and $\rho_{V}$, the miniaturized functions are continuous along the boundaries of the cells of order one.  Since $u$ satisfies Neumann boundary conditions, it follows that $u_{+}$ or $u_{-}$ satisfy matching conditions along these boundaries, hence they are Neumann eigenfunctions on $u_{m+1}$, and the eigenvalue is multiplied by $\lambda^{-2}$ where $\lambda$ denotes the contraction ratio of the $F_{i}$ mappings (so $\lambda=\frac{1}{3}$ for SC and $\lambda=\frac{1}{4}$ for the $\frac{12}{16}$ carpet).  Note that on the $\frac{12}{16}$ carpet, the miniaturized eigenfunction has the same representation type as $u$, while on SC, $u_{+}$ preserves representation type while $u_{-}$ maps $1+-$ to $1++$ and $1--$ to $1-+$.

\begin{figure}
\begin{center}
\begin{tabular}{cc}
\resizebox{3.75cm}{!}{
\begin{tabular}{|@{}r@{}|@{}r@{}|@{}r@{}|}
\hline
\makebox[.48cm][r]{$u$\:} & \makebox[.48cm][r]{-$u$\:} &  \makebox[.48cm][r]{$u$\:}\\
\hline
\makebox[.48cm][r]{-$u$\:} &  & \makebox[.48cm][r]{-$u$\:}\\
\hline
\makebox[.48cm][r]{$u$\:} & \makebox[.48cm][r]{-$u$\:} & \makebox[.48cm][r]{$u$\:}\\
\hline
\end{tabular}
}
&
\resizebox{5cm}{!}{
\begin{tabular}{|@{}r@{}|@{}r@{}|@{}r@{}|@{}r@{}|}
\hline
\makebox[.48cm][r]{$u$\:} & \makebox[.48cm][r]{-$u$\:} & \makebox[.48cm][r]{$u$\:} & \makebox[.48cm][r]{-$u$\:}\\
\hline
\makebox[.48cm][r]{-$u$\:} &  &  & \makebox[.48cm][r]{$u$\:}\\
\hline
\makebox[.48cm][r]{$u$\:} & & & \makebox[.48cm][r]{-$u$\:}\\
\hline
\makebox[.48cm][r]{-$u$\:} & \makebox[.48cm][r]{$u$\:} & \makebox[.48cm][r]{-$u$\:} & \makebox[.48cm][r]{$u$\:}\\
\hline
\end{tabular}
}
\\
\\
\mbox{(a)} & \mbox{(b)}
\end{tabular}
\end{center}
\caption{1-D Miniaturized Carpet Eigenfunctions}
\label{figfive1}
\end{figure}

There is also a two-dimensional representation of $D_{4}$, that we denote by 2.  The representation space is spanned by functions $u$ and $v$ satisfying $v=\rho_{H}u=-\rho_{V}u$ and $\rho''_{D}u=-\rho'_{D}u=u$, $\rho'_{D}v=-\rho''_{D}v=v$.  The miniaturized functions $u_{2}$ and $v_{2}$ are shown in Figure \ref{figfive2}. Once again we see that $u_{2}$ and $v_{2}$ are Neumann eigenfunctions on $\Omega_{m+1}$ with eigenvalue multiplied by $\lambda^{-2}$, and the pair transform according to the representation 2.

What does this tell us about the Neumann spectrum on the corresponding fractal?  If we believe (\ref{one4}) then there will be an eigenvalue renormalization factor $R=r\lambda^{-2}$.  For every eigenvalue $\lambda_{n}$, there will be an eigenvalue equal to $R\lambda_{n}$ with equal multiplicity, and the corresponding eigenfunctions will be miniaturizations as illustrated.

\begin{figure}
\begin{center}
\begin{tabular}{cc}
\resizebox{3.75cm}{!}{
\begin{tabular}{|@{}r@{}|@{}r@{}|@{}r@{}|}
\hline
\makebox[.48cm][r]{$u$\:} & \makebox[.48cm][r]{-$v$\:} &  \makebox[.48cm][r]{$u$\:}\\
\hline
\makebox[.48cm][r]{$v$\:} &  & \makebox[.48cm][r]{$v$\:}\\
\hline
\makebox[.48cm][r]{$u$\:} & \makebox[.48cm][r]{-$v$\:} & \makebox[.48cm][r]{$u$\:}\\
\hline
\end{tabular}
}
&
\resizebox{5cm}{!}{
\begin{tabular}{|@{}r@{}|@{}r@{}|@{}r@{}|@{}r@{}|}
\hline
\makebox[.48cm][r]{$u$\:} & \makebox[.48cm][r]{-$v$\:} & \makebox[.48cm][r]{$u$\:} & \makebox[.48cm][r]{-$v$\:}\\
\hline
\makebox[.48cm][r]{$v$\:} &  &  & \makebox[.48cm][r]{-$u$\:}\\
\hline
\makebox[.48cm][r]{$u$\:} & & & \makebox[.48cm][r]{-$v$\:}\\
\hline
\makebox[.48cm][r]{$v$\:} & \makebox[.48cm][r]{-$u$\:} & \makebox[.48cm][r]{$v$\:} & \makebox[.48cm][r]{-$u$\:}\\
\hline
\end{tabular}
}
\\
\\
\resizebox{3.75cm}{!}{
\begin{tabular}{|@{}r@{}|@{}r@{}|@{}r@{}|}
\hline
\makebox[.48cm][r]{$v$\:} & \makebox[.48cm][r]{-$u$\:} &  \makebox[.48cm][r]{$v$\:}\\
\hline
\makebox[.48cm][r]{$u$\:} &  & \makebox[.48cm][r]{$u$\:}\\
\hline
\makebox[.48cm][r]{$v$\:} & \makebox[.48cm][r]{-$u$\:} & \makebox[.48cm][r]{$v$\:}\\
\hline
\end{tabular}
}
&
\resizebox{5cm}{!}{
\begin{tabular}{|@{}r@{}|@{}r@{}|@{}r@{}|@{}r@{}|}
\hline
\makebox[.48cm][r]{$v$\:} & \makebox[.48cm][r]{-$u$\:} & \makebox[.48cm][r]{-$v$\:} & \makebox[.48cm][r]{-$u$\:}\\
\hline
\makebox[.48cm][r]{$u$\:} &  &  & \makebox[.48cm][r]{-$v$\:}\\
\hline
\makebox[.48cm][r]{$v$\:} & & & \makebox[.48cm][r]{-$u$\:}\\
\hline
\makebox[.48cm][r]{$u$\:} & \makebox[.48cm][r]{-$v$\:} & \makebox[.48cm][r]{-$u$\:} & \makebox[.48cm][r]{-$v$\:}\\
\hline
\end{tabular}
}
\\
\\
\mbox{(a)} & \mbox{(b)}
\end{tabular}
\end{center}
\caption{1-D Miniaturized Carpet Eigenfunctions}
\label{figfive2}
\end{figure}

But in fact we can run the same miniaturization argument directly on the fractal.  Indeed, in both cases we know that there exists a Laplacian $\Delta$ on the fractal satisfying a self-similar identity
\begin{equation}\label{five5}
\Delta(u\circ F_{i})=R^{-1}(\Delta u)\circ F_{i}
\end{equation}
for a certain constant $R$.  Then the miniaturization recipes given above create eigenfunctions with eigenvalue multiplied by $R$.  This is true independent of the validity of the outer approximation method.  Incidentally, the miniaturization recipes given above extend easily to any $D_{4}$ symmetric carpet type fractal.

In our last example, the octagasket, the symmetry group is $D_{8}$.  Here we have four one-dimensional representations.  Since $D_{4}\subset D_{8}$ we may sort the reflections in $D_{8}$ into those that are in $D_{4}$ and those that are not.  The representation $1++$ is described by functions even with respect to all reflections, and $1--$ by all functions odd with respect to all reflections.  Similarly, $1+-$ functions are odd with respect to $D_{4}$ reflections and even with respect to all other reflections, while for $1-+$ functions the situation is reversed.  The miniaturizations $u_{+}$ (for $1++$ or $1+-$ eigenfunctions) and $u_{-}$ (for $1-+$ or $1--$ eigenfunctions) are again given by (\ref{five3}) and (\ref{five4}), where the $\pm$ signs alternate along the eight small octagons.  We note that the representation type is preserved under miniaturization.

\begin{figure}
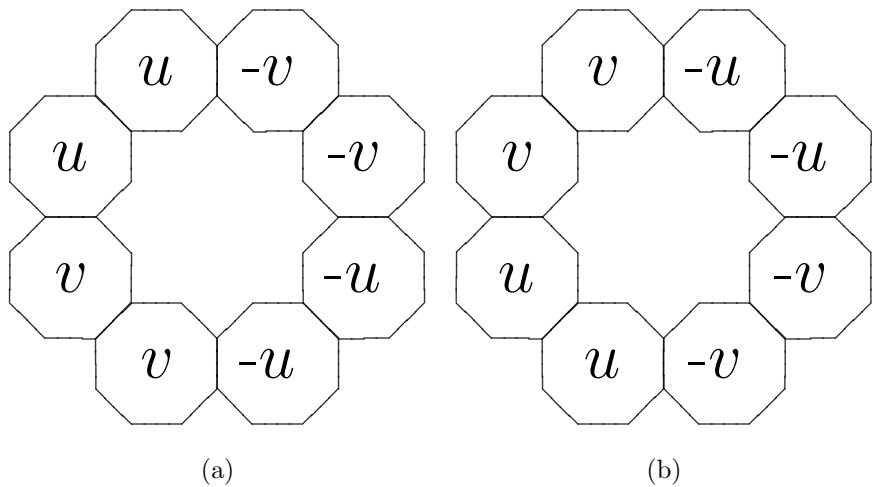

\begin{center}
\begin{tabular}{cc}
\xy{
\xypolygon8"A"{~:{/r5pc/:}~>{}
~<>{\xypolygon8{~:{/-2.071067pc/:}~={0}}}
~*{\alphanum}
}
}
\endxy
&
\xy
{
\xypolygon8"B"{~:{/r5pc/:}~>{}
~<>{\xypolygon8{~:{/-2.071067pc/:}~={0}}}
~*{\betanum}
}
}
\endxy
\\
\\
\mbox{(a)} & \mbox{(b)}
\end{tabular}
\end{center}
\caption{The miniaturizations (a) $u_{2}$ and (b) $v_{2}$ for a $2_{1}$ or $2_{3}$ eigenspace}
\label{figfive3}
\end{figure}

\begin{figure}
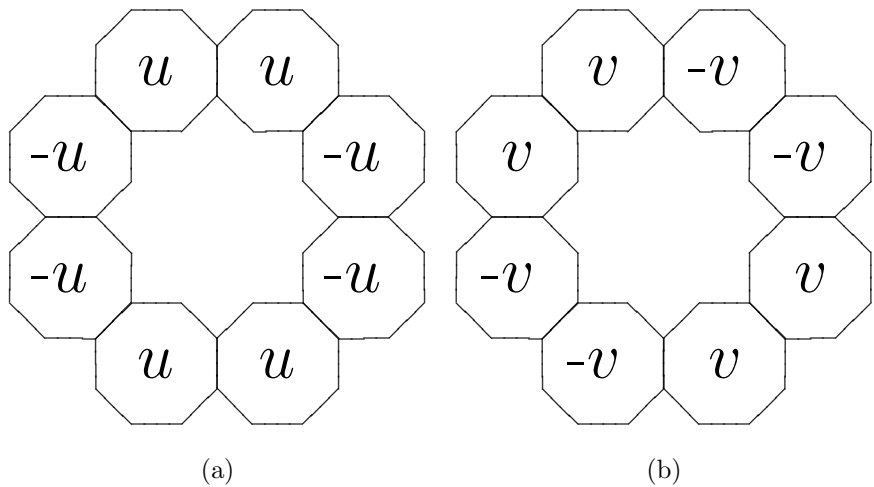

\begin{center}
\begin{tabular}{cc}
\xy{
\xypolygon8"A"{~:{/r5pc/:}~>{}
~<>{\xypolygon8{~:{/-2.071067pc/:}~={0}}}
~*{\gammanum}
}
}
\endxy
&
\xy
{
\xypolygon8"B"{~:{/r5pc/:}~>{}
~<>{\xypolygon8{~:{/-2.071067pc/:}~={0}}}
~*{\deltanum}
}
}
\endxy
\\
\\
\mbox{(a)} & \mbox{(b)}
\end{tabular}
\end{center}
\caption{The miniaturizations (a) $u_{2}'$ and (b) $v_{2}'$ for a $2_{2}$ eigenspace}
\label{figfive4}
\end{figure}

In this case there are three two-dimensional representations, denoted $2_{1},2_{2},2_{3}$.  In terms of complex valued functions on the circle, $2_{1}$ is spanned by $e^{\pm 2\pi i\theta/8}$, $2_{2}$ is spanned by $e^{\pm 2\pi i2\theta/8}$, and $2_{3}$ is spanned by $e^{\pm 2\pi 3i\theta/8}$.  If $x,y,z$ denote any consecutive points on an eight element orbit of $D_{8}$, then $2_{1}$ functions satisfy
\begin{equation}\label{five6}
u(y)=\frac{\sqrt{2}}{2}(u(x)+u(z)),
\end{equation}
$2_{2}$ functions satisfy
\begin{equation}\label{five7}
u(x)+u(z)=0,
\end{equation}
and $2_{3}$ functions satisfy
\begin{equation}\label{five8}
u(y)=-\frac{\sqrt{2}}{2}(u(x)+u(z)).
\end{equation}
The $2_{1}$ and $2_{3}$ representations have the property that restricted to $D_{4}$ they become the 2 representation.  So if $u,v$ are the basis described above, the miniaturization $u_{2},v_{2}$ are given in Figure \ref{figfive3}.  On the other hand, the restriction of $2_{2}$ to $D_{4}$ splits into a direct sum of a $1+-$ and a $1-+$ representation.  So we can choose a basis $u,v$ such that $\rho_{H}u=\rho_{V}u=u=-\rho'_{D}u=-\rho''_{D}u$ and $-\rho_{H}v=-\rho_{V}v=v=\rho'_{D}v=\rho''_{D}v$, and the miniaturization $u'_{2},v'_{2}$ is given in Figure \ref{figfive4}.  Again the representation type is preserved under miniaturization.

Some types of miniaturization on the pentagasket are described in \cite{adams01}.

\clearpage
\section{Random Carpets}
\label{secranc}

For $j\in\Z$, $j>1$, we partition the unit square into a grid of $j$ by $j$ smaller, equally sized squares of width $1/j$.  We then randomly remove $k$ of these smaller squares, where $k$ is a small positive integer, and the result is our level 1 domain $\Omega_{1}$.  To produce $\Omega_{2}$, we partition each square of width $1/j$ into a grid of $j$ by $j$ equally sized squares of width $1/j^{2}$, and we then randomly remove $m$ squares of width $1/j^{2}$ from each square of width $1/j$.  Iterating this process yields a sequence of nested compact domains $\{\Omega_{m}\}_{m=1}^{\infty}$ where $\Omega_{m}$ is a union of squares of side length $j^{-m}$.  Matlab's \verb!rand('state')! function, a modified version of Marsaglia's Subtract-with-Borrow algorithm, makes our random choices.  The number generator's state is set according to the exact date and time of the computation, so that the generator's own state is essentially randomly determined.  Also, to shorten  FEM computation time we triangulate $\Omega_{m}$ with the four sides and two diagonals of each square of side length $j^{-m}$.

\begin{figure}[htbp!]
\begin{center}
\includegraphics[width=\onewidth\textwidth]{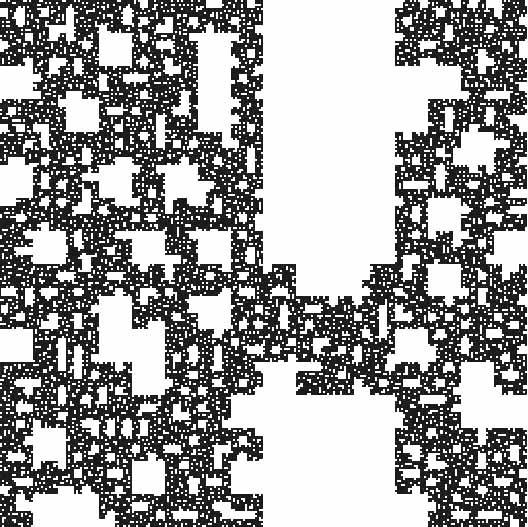}
\end{center}
\caption{Level 4 Domain $\Omega_{4}$ for $j=4$, $k=3$}
\label{figsix1}
\end{figure}

The problem we find with our FEM eigenvalue problem on these domains is connectivity.  How can we guarantee that each $\Omega_{m}$ has only one path component?  Also, if two squares are disjoint except at a common vertex, with no other squares in a neighborhood of that vertex, how can we avoid the problem we saw in Section \ref{secgasket}? Recall that in this case, the spline space of our finite element solver couples these squares at the common vertex.  For simplicity we resolve both questions by choosing small $k$ and altering the above algorithm so that this coupling problem is avoided, as follows.  When we pass from $\Omega_{m}$ to $\Omega_{m+1}$ we partition a square of side length $j^{-m}$ into squares of side length $j^{-m-1}$ and delete $k$ of the smaller squares randomly.  We then check if this deletion process has produced the above coupling problem.  If it has, then we go back and try again; otherwise, we move on to the next $m^{th}$ level square, and so on.  For $k$ small enough, the algorithm terminates.  Figure \ref{figsix1} shows a typical result of the above algorithm.  Notice that we have only one path component.

Now, we study our spectral information with the eigenvalue counting function $N:[0,\infty)\to\Z$, where $N(x)$ is the number of nonnegative eigenvalues less than or equal to $x$.  Then, we examine the Weyl ratio
\begin{equation}\label{six1}
W(x)=\frac{N(x)}{x^{\alpha}}
\end{equation}
where $x^{\alpha}$ is an approximate asymptotic bound for $N(x)$, i.e. we choose $\alpha\in\R$ so that $N(x)\sim x^{\alpha}$ in accordance with the experimental data.  So, finding $\alpha$ corresponds to finding the slope of a linear approximation of $N(x)$ on a log-log plot.  In fact, since we are dealing with domains in the plane, the Weyl asymptotic law implies that $\alpha=1$ is the correct value as $x\to\infty$.  The point is that we truncate our computations well before we reach the region where this asymptotic behavior is approximated, so we observe values of $\alpha$ considerably smaller than $1$.

\begin{figure}[htbp!]
$$
\begin{array}{cc}
\includegraphics[width=\foursqwidth\textwidth]{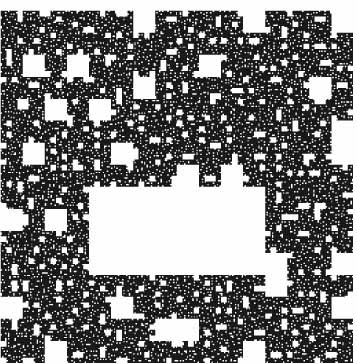} &
\includegraphics[width=\foursqwidth\textwidth]{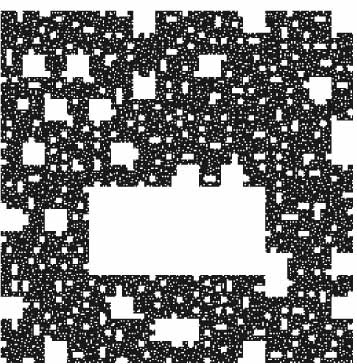}\\
\mbox{Original Carpet} & \mbox{Bifurcation at Level $4$}\\
\includegraphics[width=\foursqwidth\textwidth]{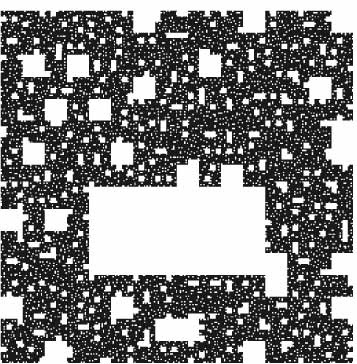} &
\includegraphics[width=\foursqwidth\textwidth]{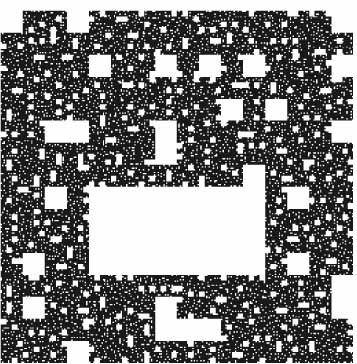}\\
\mbox{Bifurcation at Level $3$} & \mbox{Bifurcation at Level $2$}
\end{array}
$$
\caption{Carpet Bifurcations $\Omega_{4}$ for j=4, k=2}
\label{figsix2}
\end{figure}

\begin{figure}
$$
\begin{array}{cc}
\includegraphics[width=\foursqwidth\textwidth]{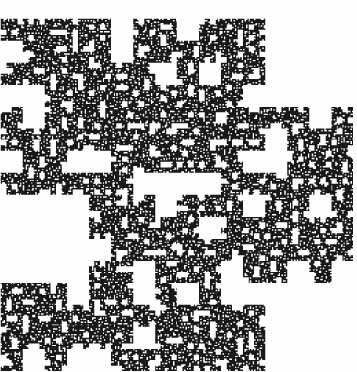} &
\includegraphics[width=\foursqwidth\textwidth]{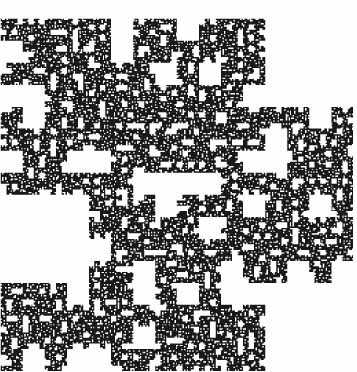}\\
\mbox{Original Carpet} & \mbox{Bifurcation at Level $4$}\\
\includegraphics[width=\foursqwidth\textwidth]{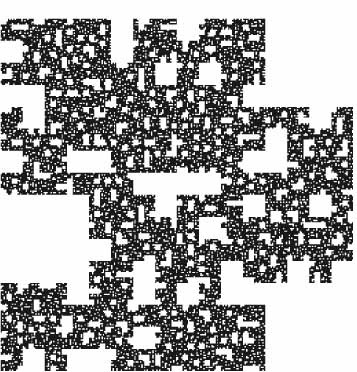} &
\includegraphics[width=\foursqwidth\textwidth]{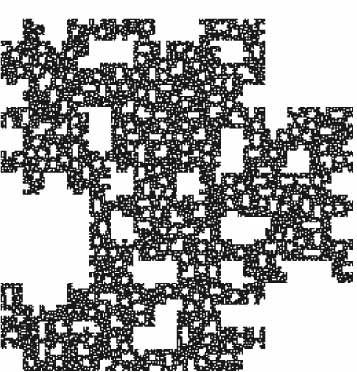}\\
\mbox{Bifurcation at Level $3$} & \mbox{Bifurcation at Level $2$}
\end{array}
$$
\caption{Carpet Bifurcations $\Omega_{4}$ for j=4, k=3}
\label{figsix4}
\end{figure}

In our first example, we let $j=4$ and $k=2$ and run our algorithm up to level $4$ to get $\{\Omega_{i}\}_{i=1}^{4}$, where $\Omega_{4}$ is the upper left carpet in Figure \ref{figsix2}.  From this initial carpet, we can restart our algorithm three separate times, beginning at $\Omega_{i}$ once for each $i=1,2,3$.  We then end the algorithm again at level $4$ and we call the resulting (level $4$) carpet which was started at $\Omega_{i}$ the bifurcation of $\Omega_{4}$ at level $i+1$.  The carpets are shown in Figure \ref{figsix2} and the eigenvalue data in Tables \ref{tablesix1} and \ref{tablesix2}.  Next, we let $j=4$ and $k=3$ and do the same bifurcation study.  The carpets are shown in Figure \ref{figsix4} and the eigenvalue data in Tables \ref{tablesix3} and \ref{tablesix4}.

\begin{figure}[htbp!]
\begin{center}
\includegraphics[width=\onewidth\textwidth]{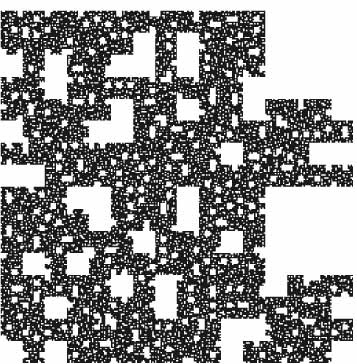}
\end{center}
\caption{Level 4 Domain $\Omega_{4}$ for $j=4$, D:2,3,2,3,2}
\label{figsix6}
\end{figure}

Finally, we fix $j=4$ and vary $k$ on different levels so that at level $1$ we set $k=2$, at level $2$ we set $k=3$, etc.  A similar procedure for gaskets rather than carpets is discussed in \cite{drenning08}.  Our sequence of $k$ values for the carpet in Figure \ref{figsix6} is $k=\{2,3,2,3,2\}$.  The eigenvalue data appears in Tables \ref{tablesix5} and \ref{tablesix6}.  The level-to-level eigenvalue ratios in Table \ref{tablesix5} appear to roughly alternate between the same ratios in Tables \ref{tablesix3} and \ref{tablesix1}.  This is the strongest evidence that the geometry of the domain at different scales is reflected in the spectrum of the Laplacian.  Such a correlation is more striking in \cite{drenning08}, but the fractals there have a more coherent structure.

\begin{figure}
$$
\begin{array}{cc}
\includegraphics[width=\foursqwidtha\textwidth]{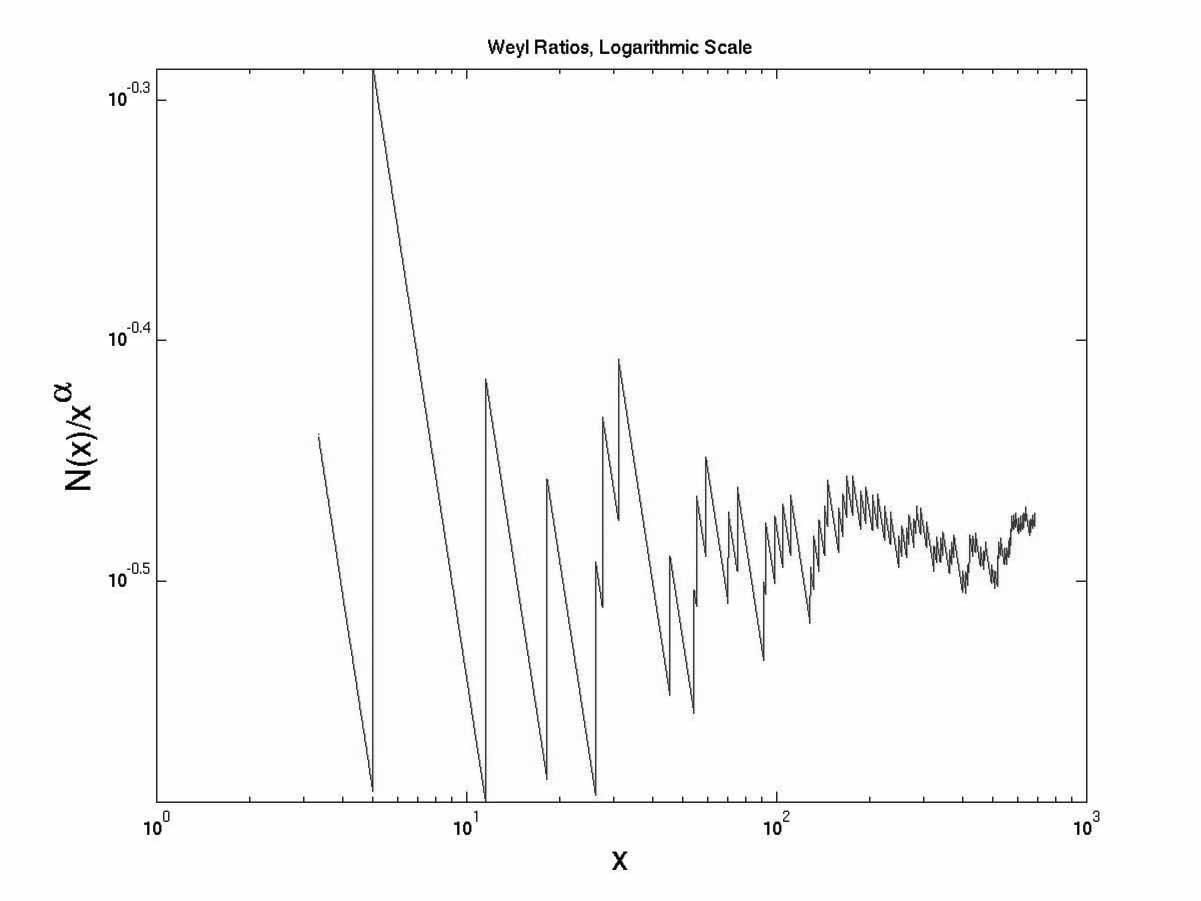} &
\includegraphics[width=\foursqwidtha\textwidth]{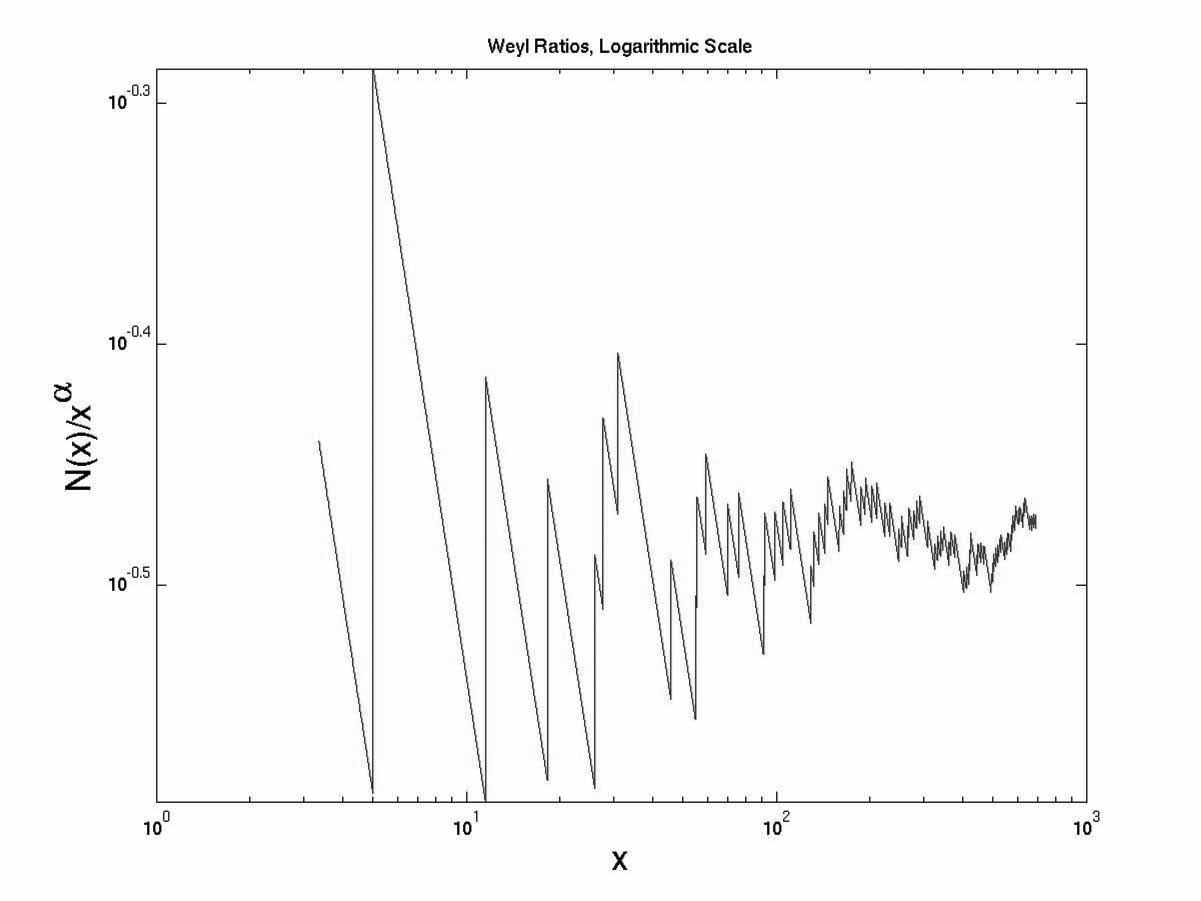}\\
\mbox{Original Carpet, $\alpha=.84032$} & \mbox{Bifurcation at Level $4$, $\alpha=.83853$}\\
\includegraphics[width=\foursqwidtha\textwidth]{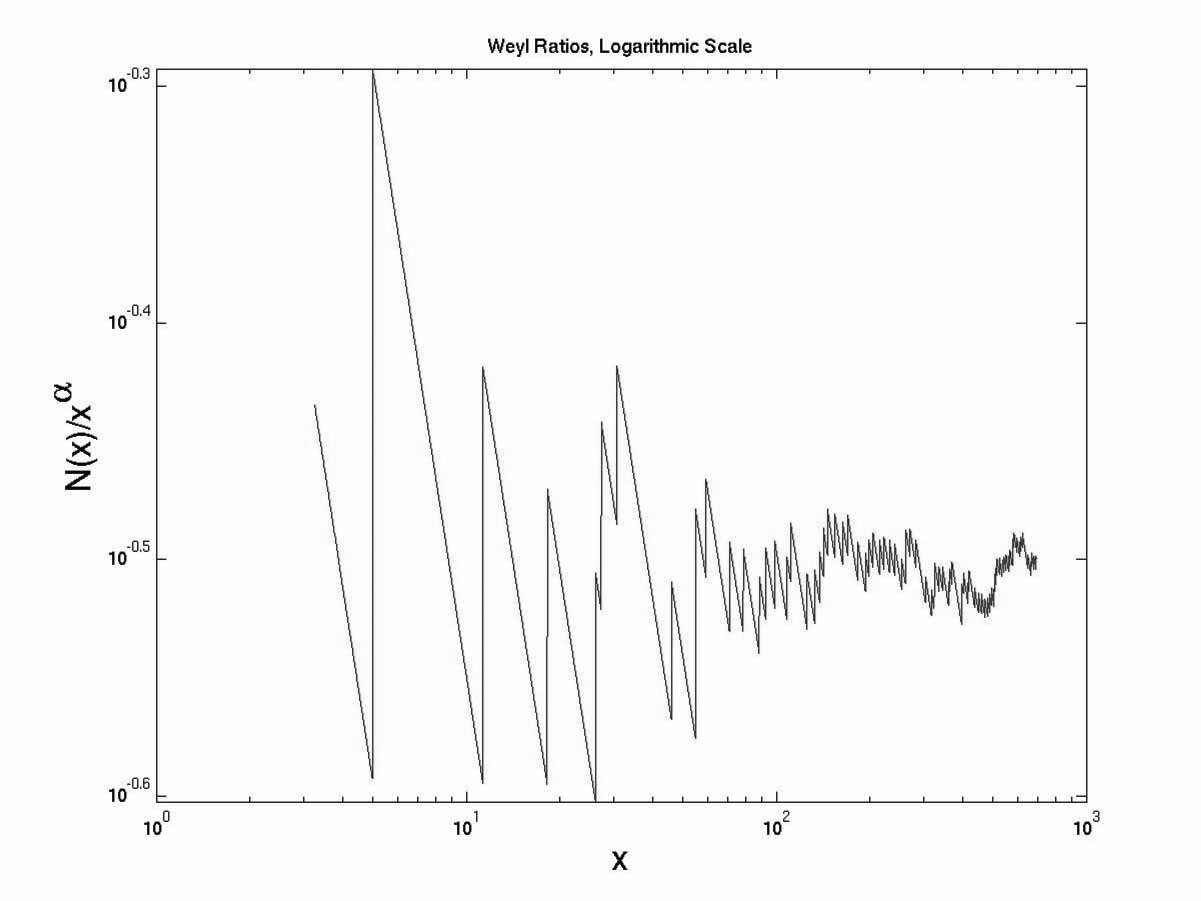} &
\includegraphics[width=\foursqwidtha\textwidth]{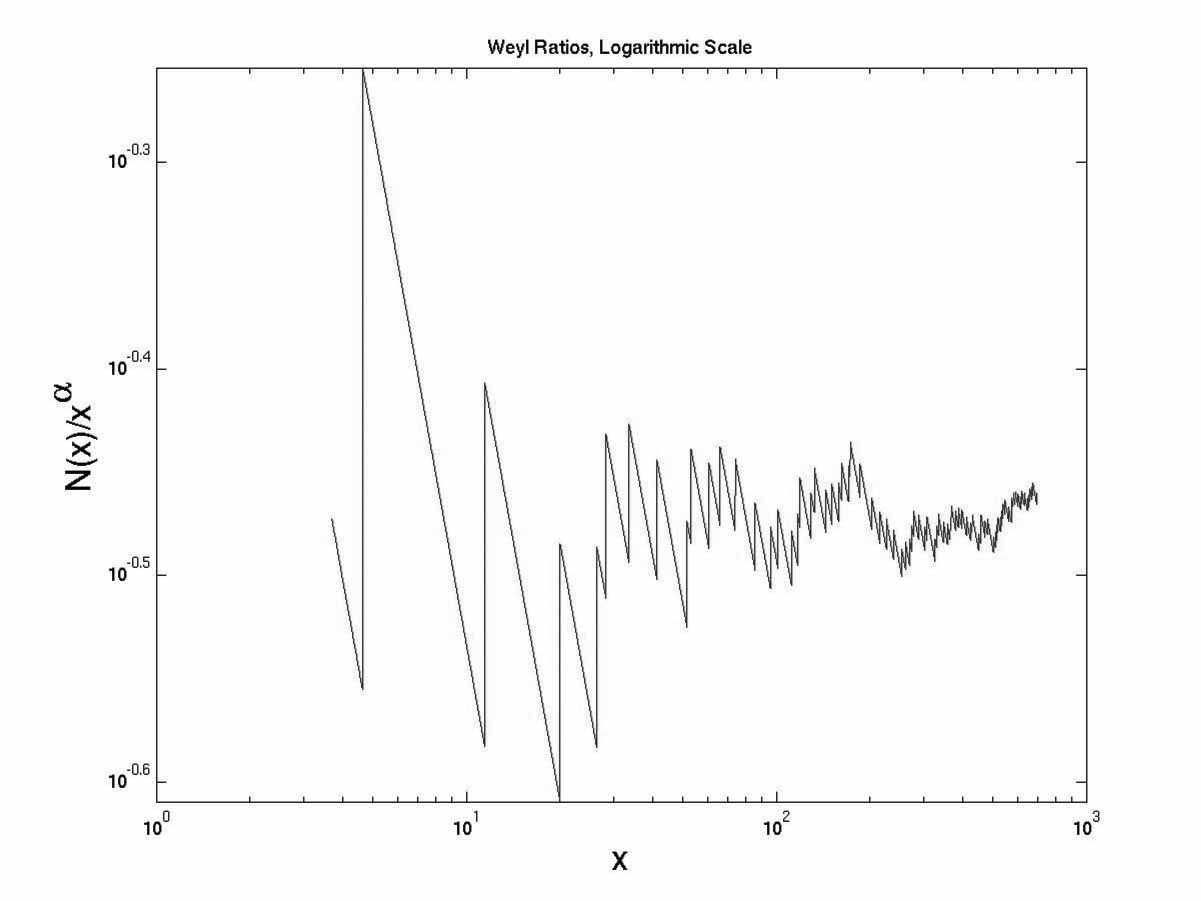}\\
\mbox{Bifurcation at Level $3$, $\alpha=.85007$} & \mbox{Bifurcation at Level $2$, $\alpha=.83383$}
\end{array}
$$
\caption{Weyl Ratios for $j=4$, $k=2$}
\label{figsix3}
\end{figure}

\begin{figure}
$$
\begin{array}{cc}
\includegraphics[width=\foursqwidtha\textwidth]{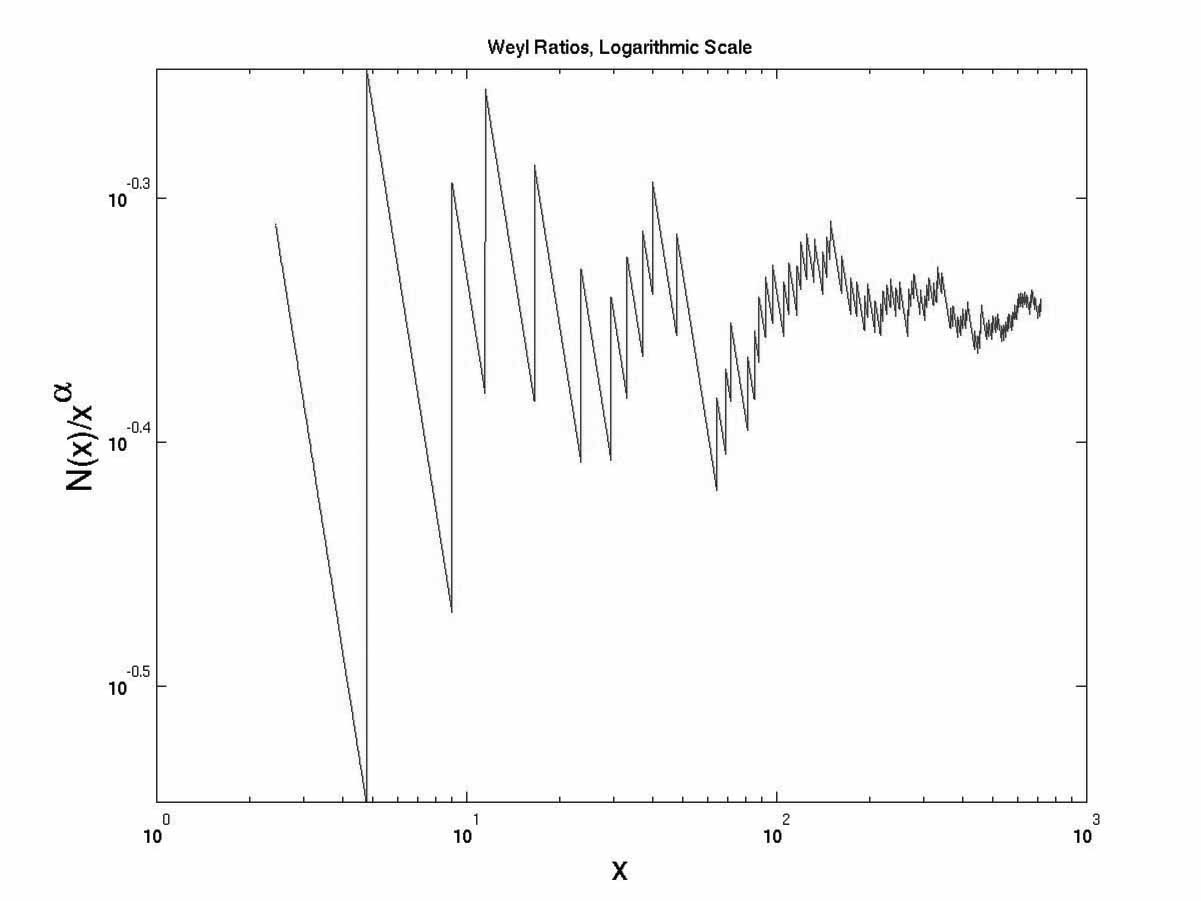} &
\includegraphics[width=\foursqwidtha\textwidth]{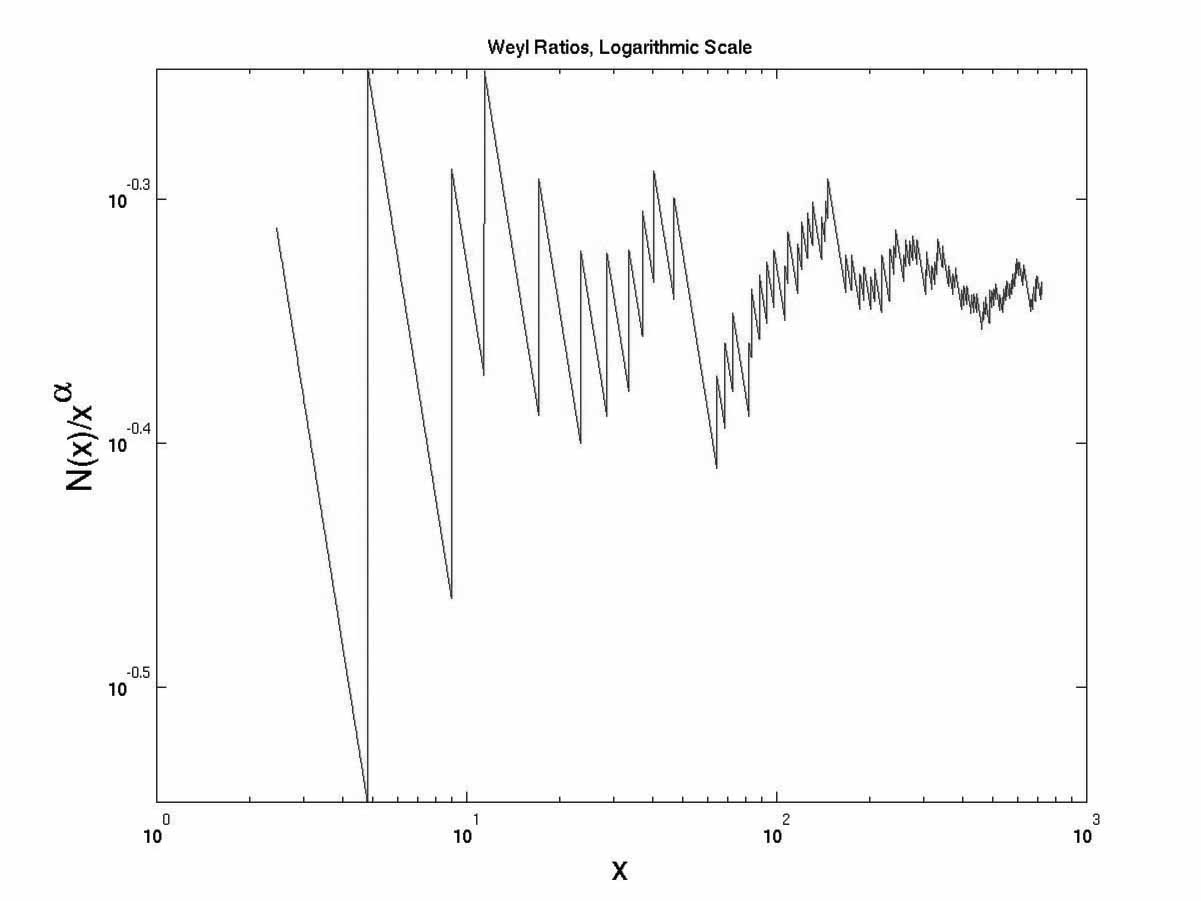}\\
\mbox{Original Carpet, $\alpha=.80788$} & \mbox{Bifurcation at Level $4$, $\alpha=.80253$}\\
\includegraphics[width=\foursqwidtha\textwidth]{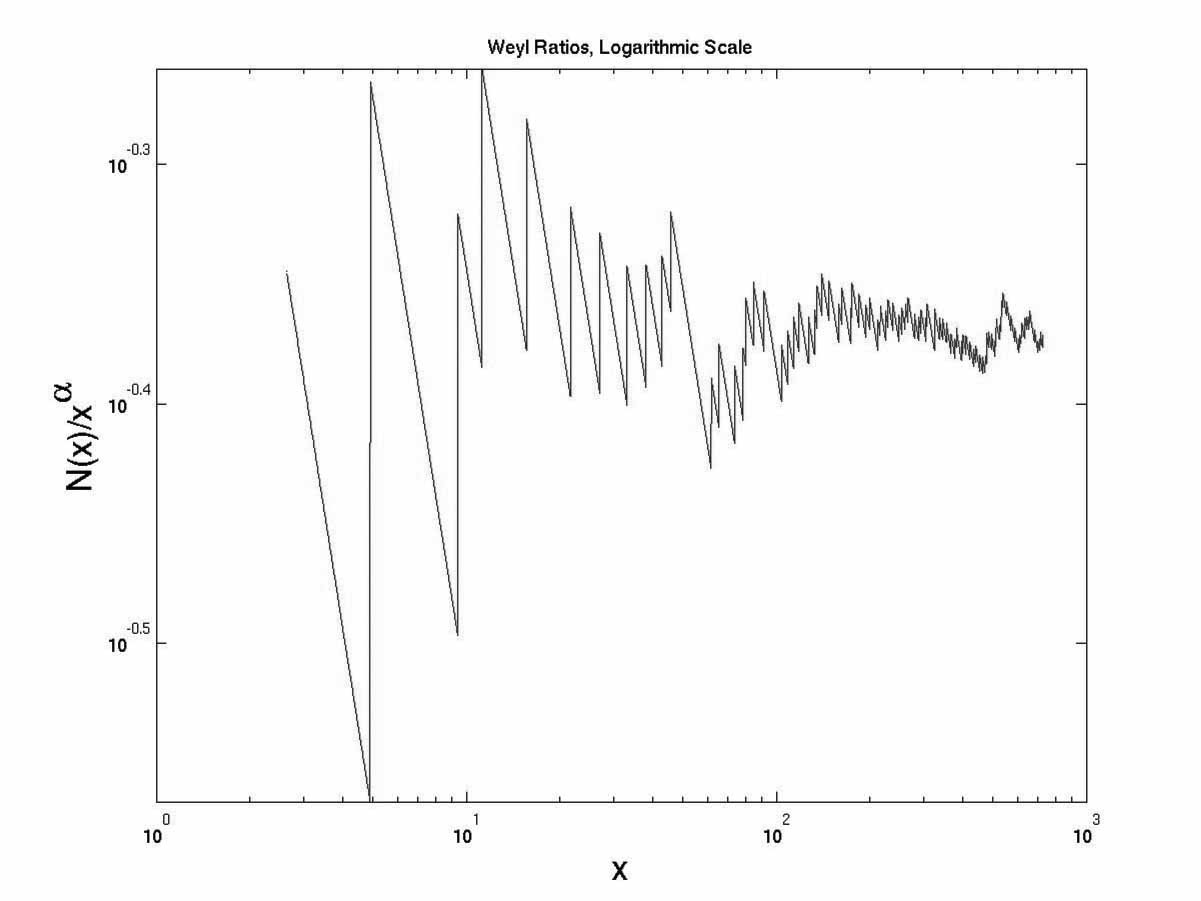} &
\includegraphics[width=\foursqwidtha\textwidth]{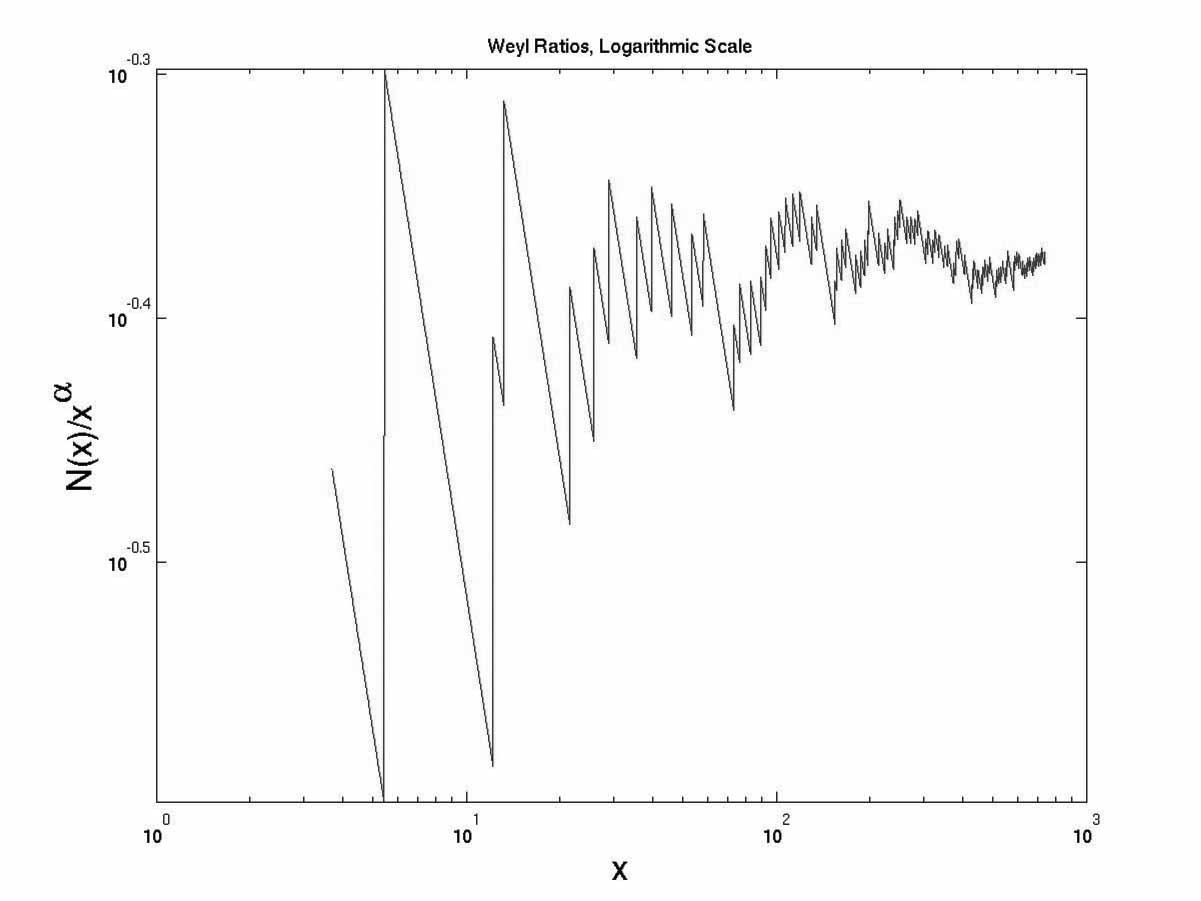}\\
\mbox{Bifurcation at Level $3$, $\alpha=.82004$} & \mbox{Bifurcation at Level $2$, $\alpha=.81408$}
\end{array}
$$
\caption{Weyl Ratios for $j=4$, $k=3$}
\label{figsix5}
\end{figure}

\begin{figure}
\includegraphics[width=\onewidth\textwidth]{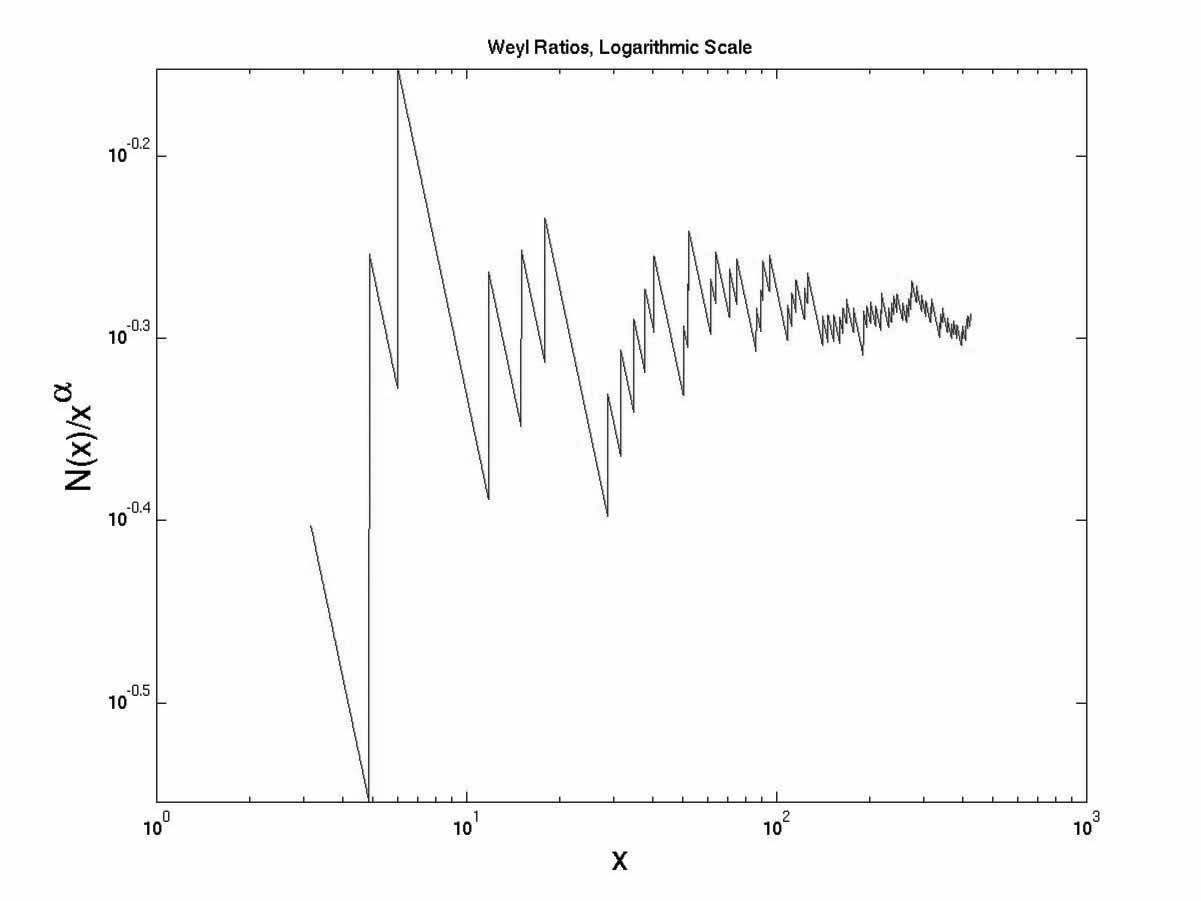}
\caption{Weyl Ratios for $j=4$, $k=\{2,3,2,3,2\}$, Level $5$ Carpet, $\alpha=.8071$}
\label{figsix7}
\end{figure}

The Weyl ratios of our first example (where $j=4$ and $k=2$) appear in Figure \ref{figsix3}.  We now look closely at the agreement of the graph of the original carpet to each individual bifurcation.  We see that the original agrees with the bifurcation at Level $4$ up to about $x=300$, the original agrees with that at Level $3$ up to around $x=65$, and it agrees with the Level $2$ bifurcation up to about $x=25$.  In our second example (where $j=4$ and $k=4$) we find the Weyl ratios in Figure \ref{figsix5}.  We do the same comparison.  The original agrees with the the Level $4$ bifurcation to around $x=150$, it agrees with the Level $3$ one up to approximately $x=30$, and it agrees with the Level $2$ bifurcation to approximately $x=10$.  In other words, the added detail at finer resolutions has only a minimal effect on some initial segment of the spectrum. This is consistent with
results in \cite{drenning08}.  Our final example's Weyl ratios (where $j=4$ and $k=\{2,3,2,3,2\}$) are found in Figure \ref{figsix7}.

\begin{figure}
$$
\begin{array}{cc}
\includegraphics[width=\foursqwidth\textwidth]{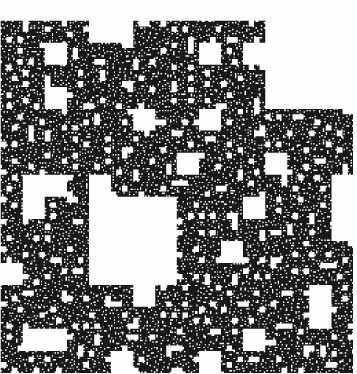} &
\includegraphics[width=\foursqwidth\textwidth]{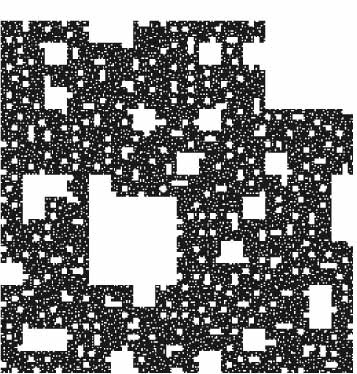}\\
\mbox{Original Carpet} & \mbox{Bifurcation at Level $4$}\\
\includegraphics[width=\foursqwidth\textwidth]{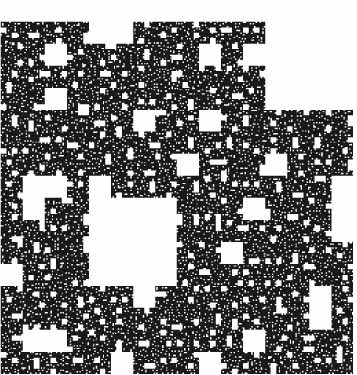} &
\includegraphics[width=\foursqwidth\textwidth]{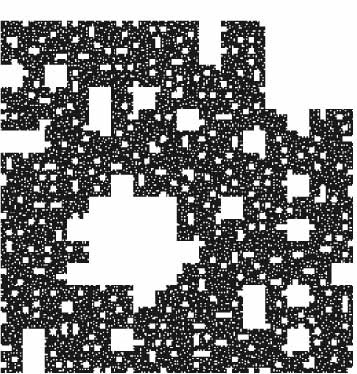}\\
\mbox{Bifurcation at Level $3$} & \mbox{Bifurcation at Level $2$}
\end{array}
$$
\caption{Carpet Bifurcations $\Omega_{4}$ for j=4, k=2}
\label{figsix8}
\end{figure}

\begin{figure}
$$
\begin{array}{cc}
\includegraphics[width=\foursqwidtha\textwidth]{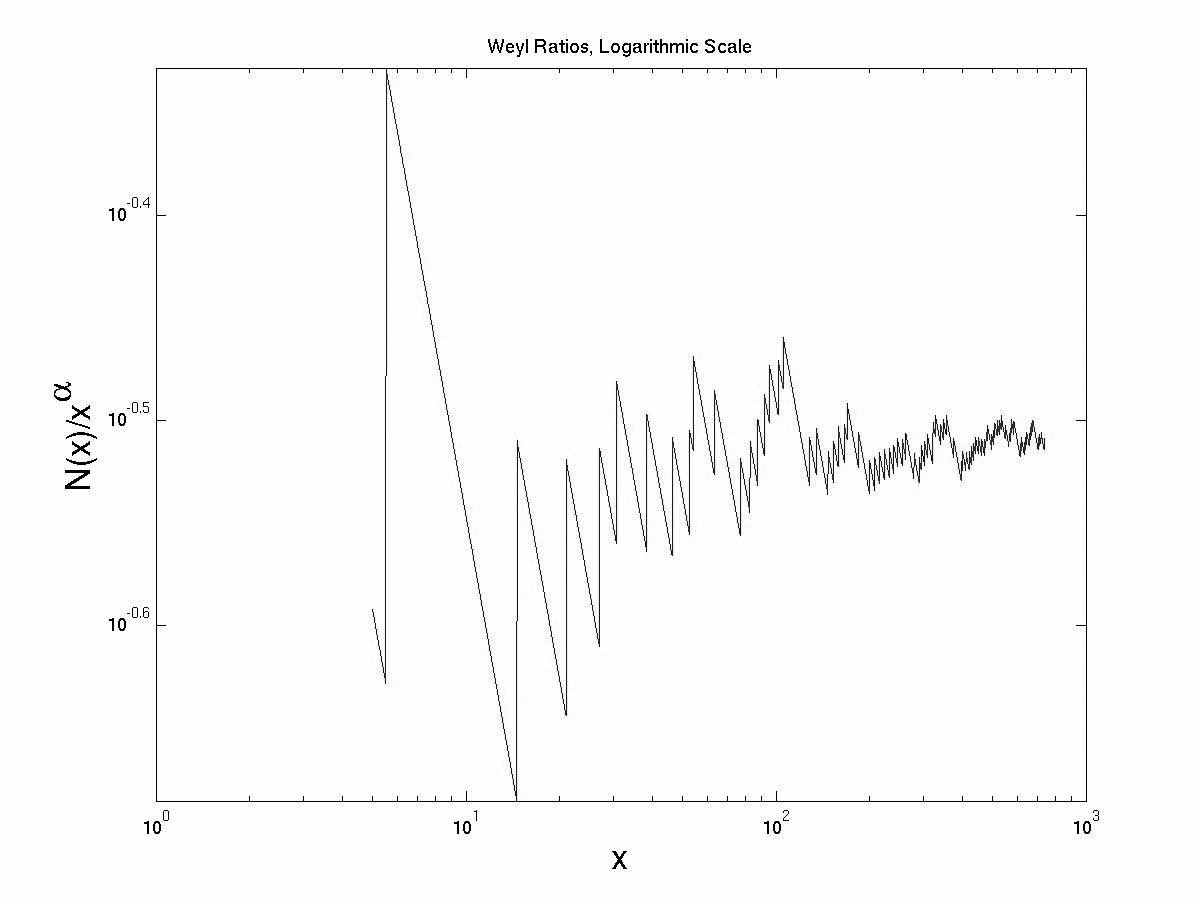} &
\includegraphics[width=\foursqwidtha\textwidth]{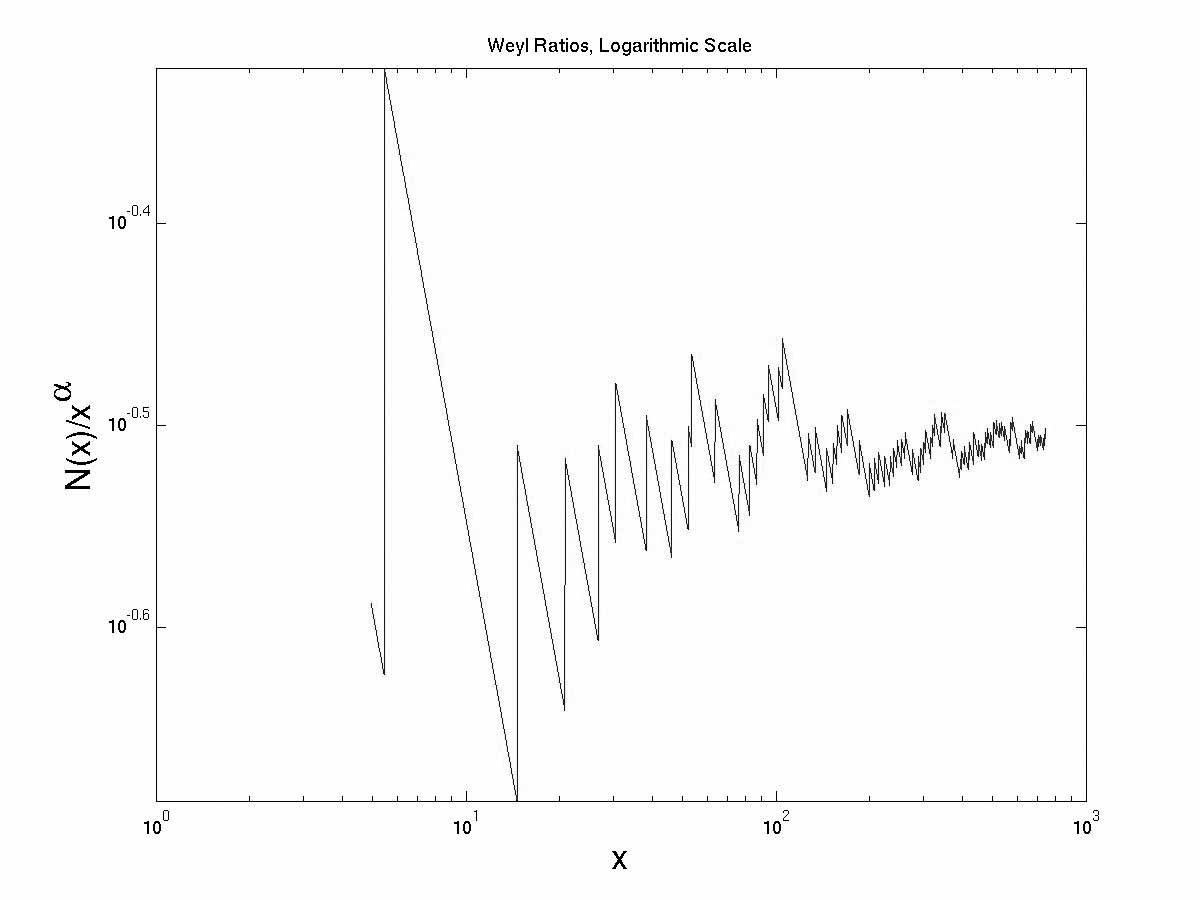}\\
\mbox{Original Carpet, $\alpha=.84747$} & \mbox{Bifurcation at Level $4$, $\alpha=.84753$}\\
\includegraphics[width=\foursqwidtha\textwidth]{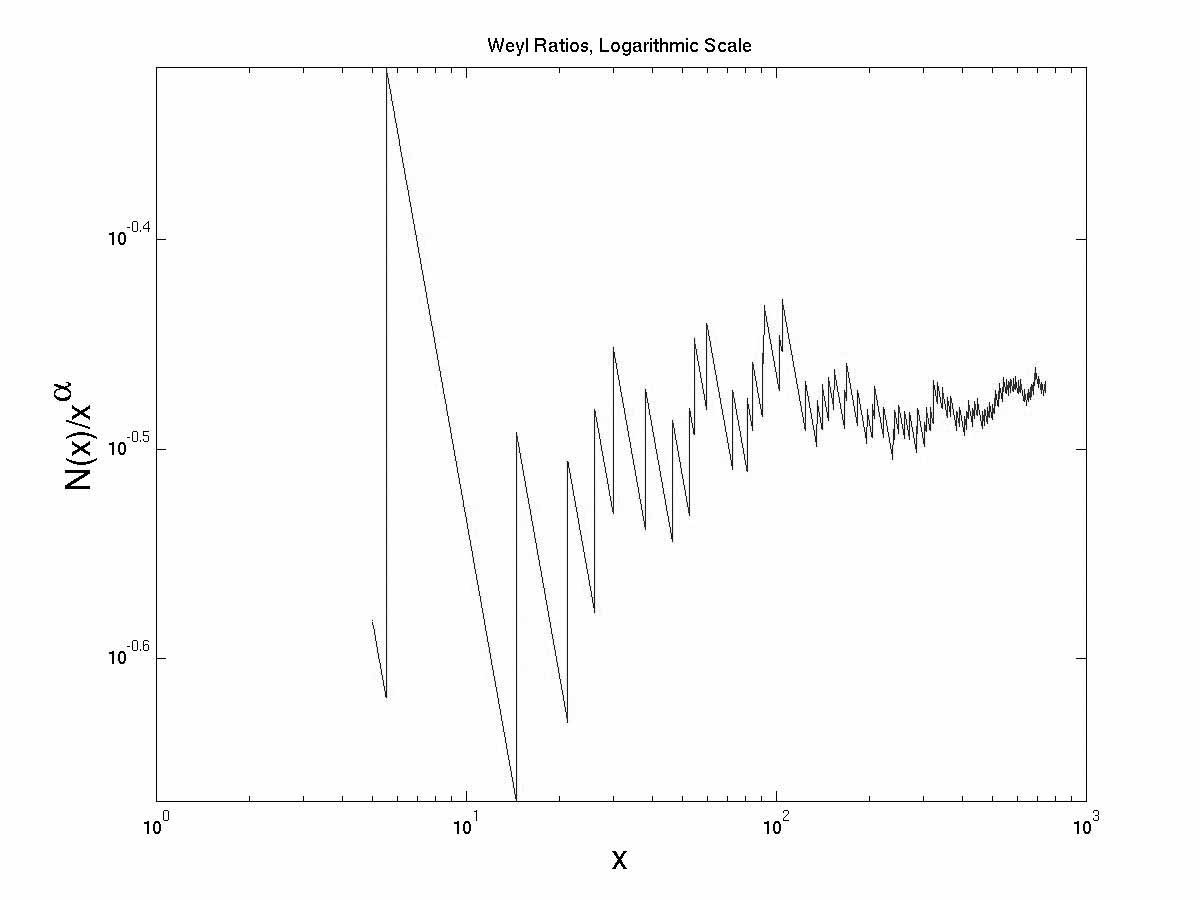} &
\includegraphics[width=\foursqwidtha\textwidth]{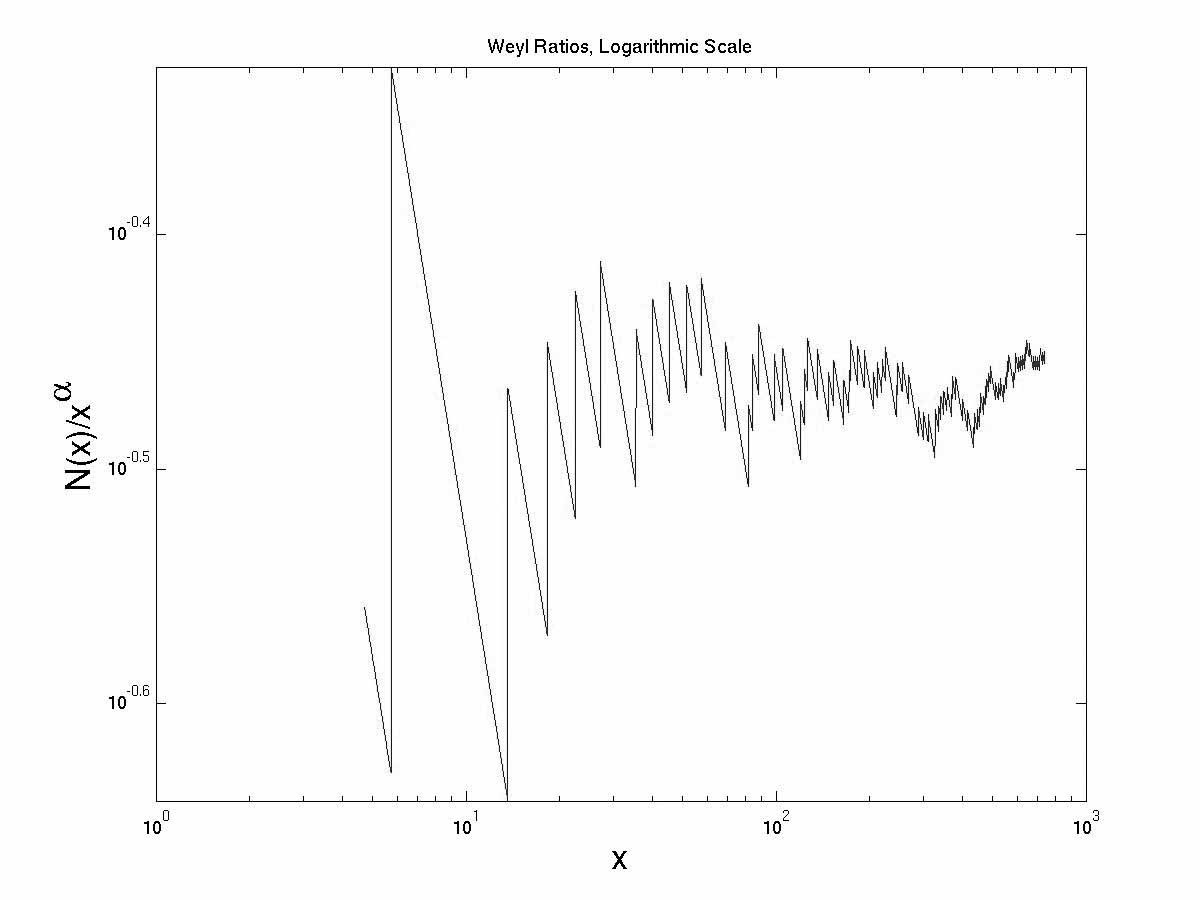}\\
\mbox{Bifurcation at Level $3$, $\alpha=.83368$} & \mbox{Bifurcation at Level $2$, $\alpha=.83019$}
\end{array}
$$
\caption{Weyl Ratios for $j=4$, $k=2$}
\label{figsix9}
\end{figure}

\begin{figure}
$$
\begin{array}{cc}
\includegraphics[width=\foursqwidth\textwidth]{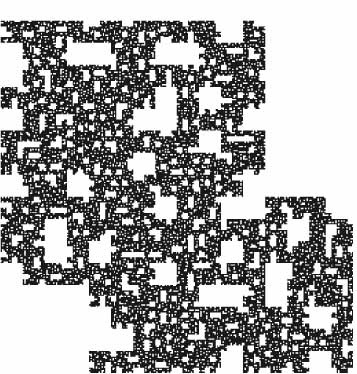} &
\includegraphics[width=\foursqwidth\textwidth]{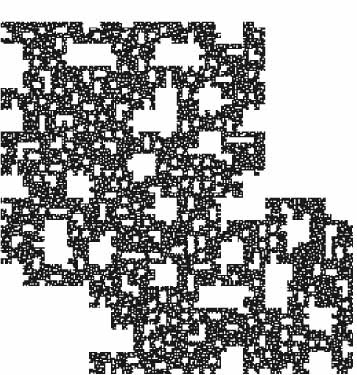}\\
\mbox{Original Carpet} & \mbox{Bifurcation at Level $4$}\\
\includegraphics[width=\foursqwidth\textwidth]{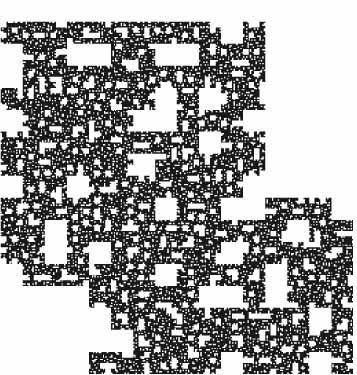} &
\includegraphics[width=\foursqwidth\textwidth]{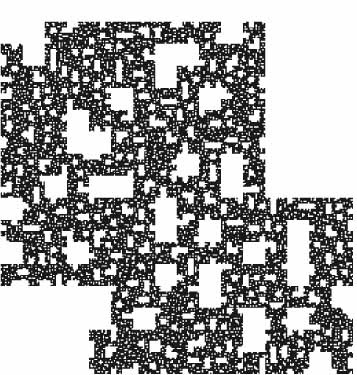}\\
\mbox{Bifurcation at Level $3$} & \mbox{Bifurcation at Level $2$}
\end{array}
$$
\caption{Carpet Bifurcations $\Omega_{4}$ for j=4, k=3}
\label{figsix10}
\end{figure}

\begin{figure}
$$
\begin{array}{cc}
\includegraphics[width=\foursqwidtha\textwidth]{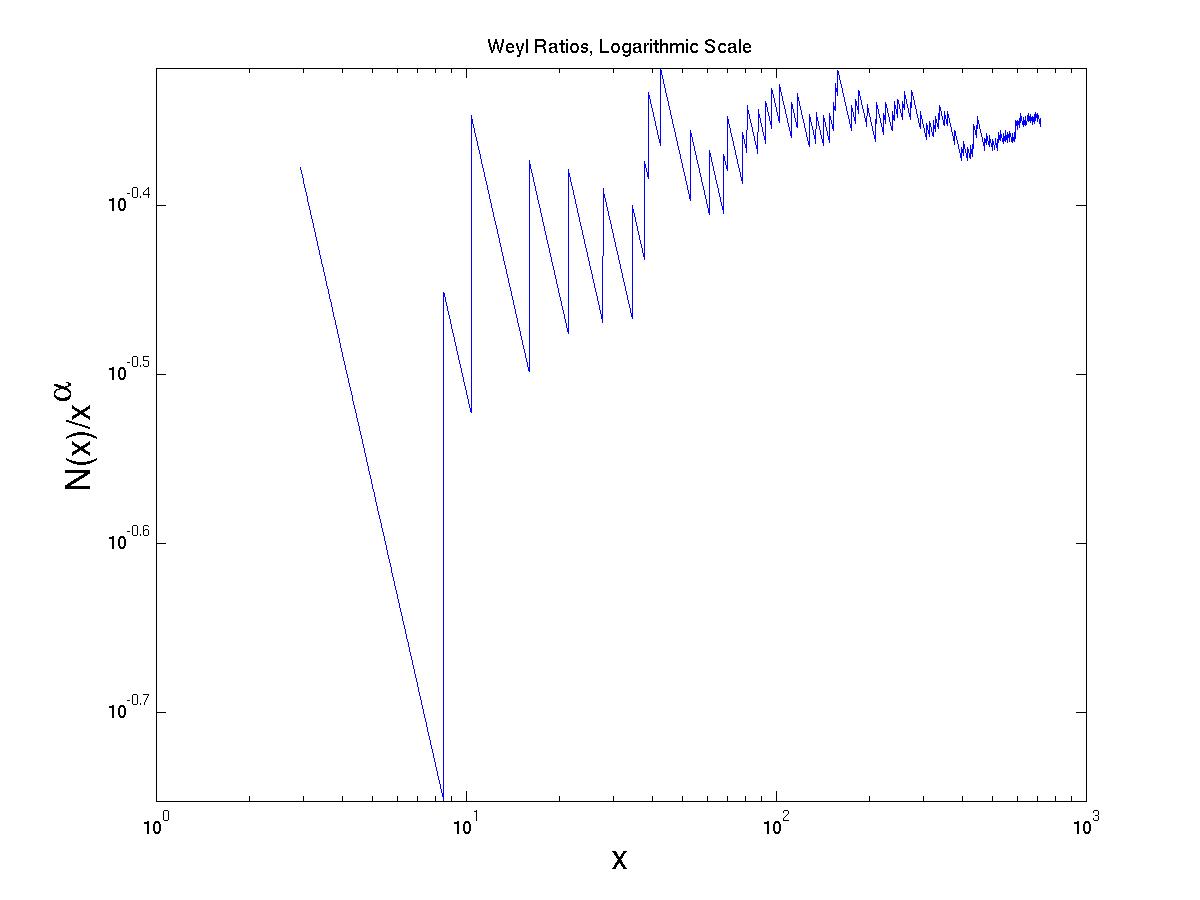} &
\includegraphics[width=\foursqwidtha\textwidth]{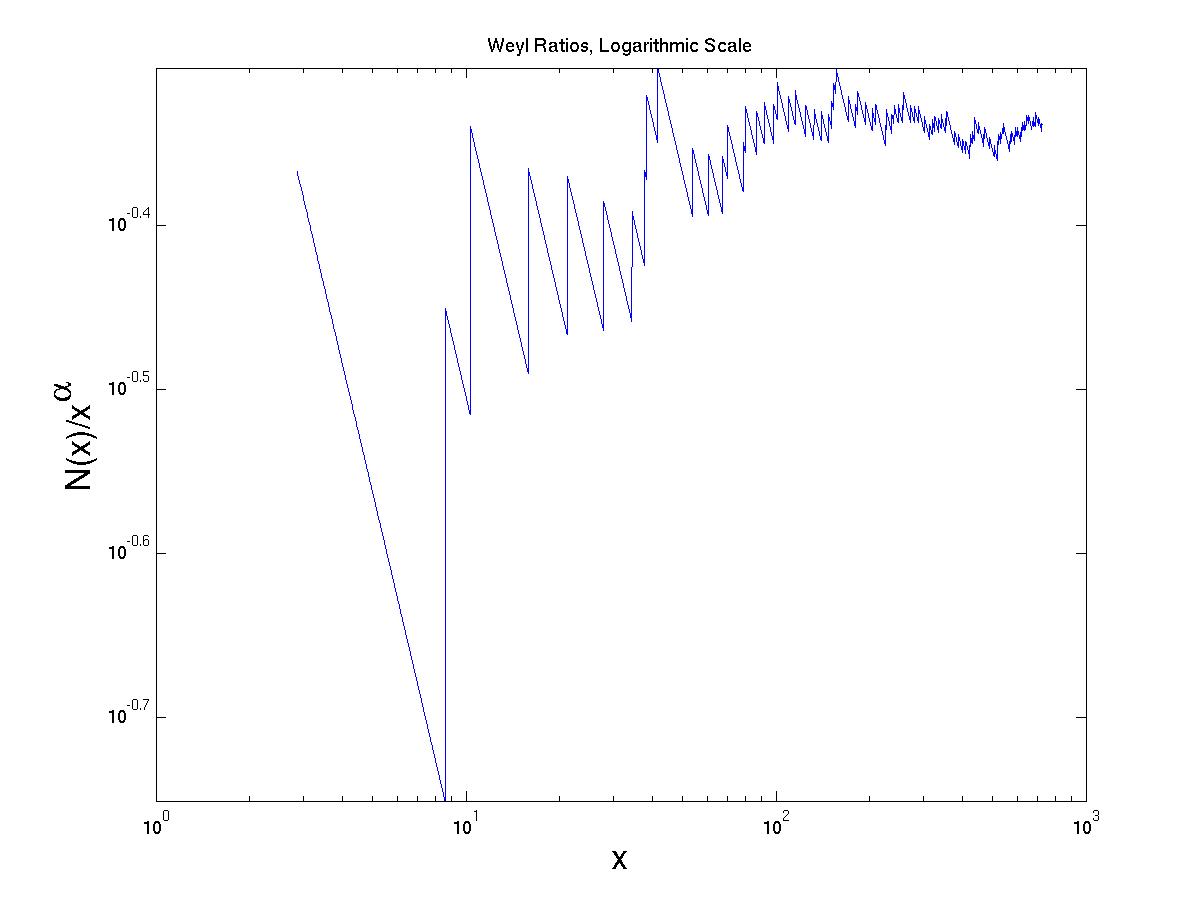}\\
\mbox{Original Carpet, $\alpha=.81013$} & \mbox{Bifurcation at Level $4$, $\alpha=.80544$}\\
\includegraphics[width=\foursqwidtha\textwidth]{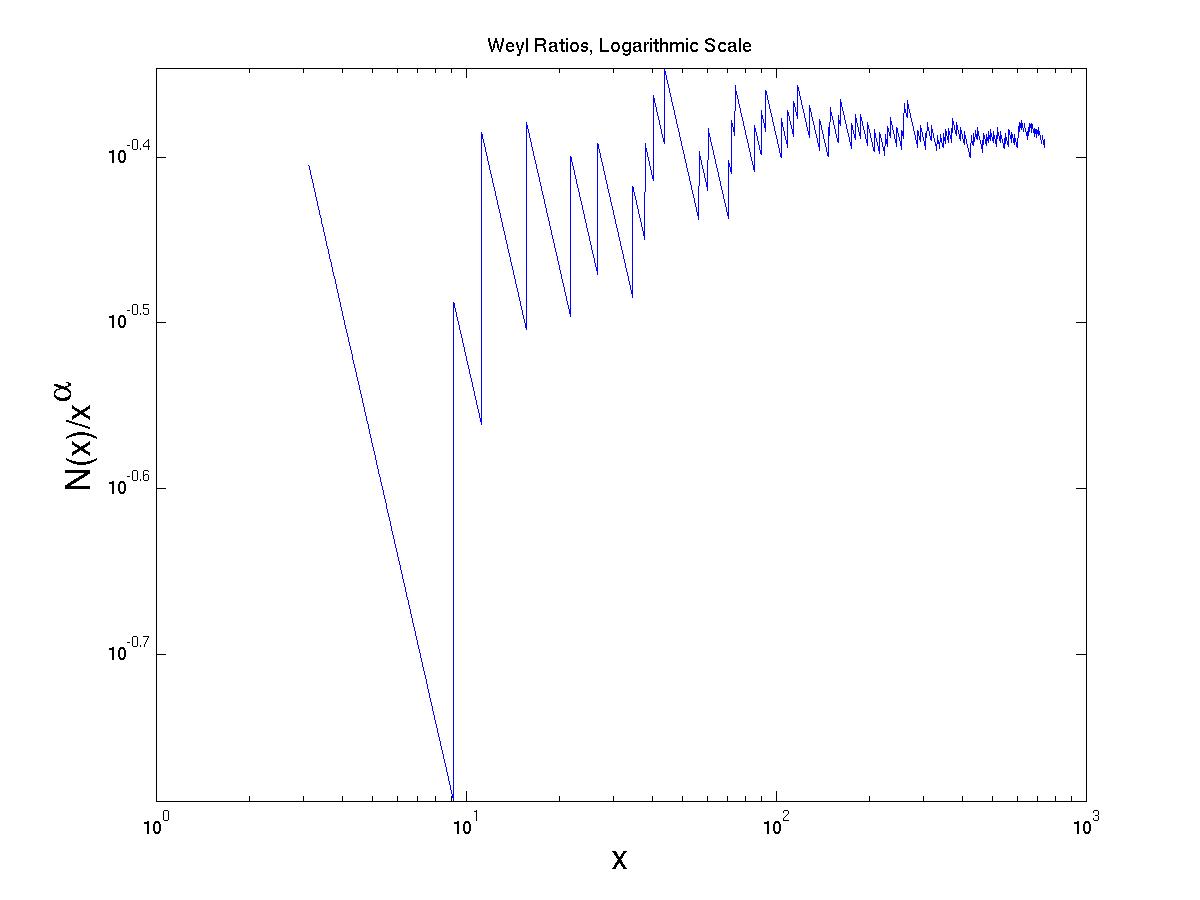} &
\includegraphics[width=\foursqwidtha\textwidth]{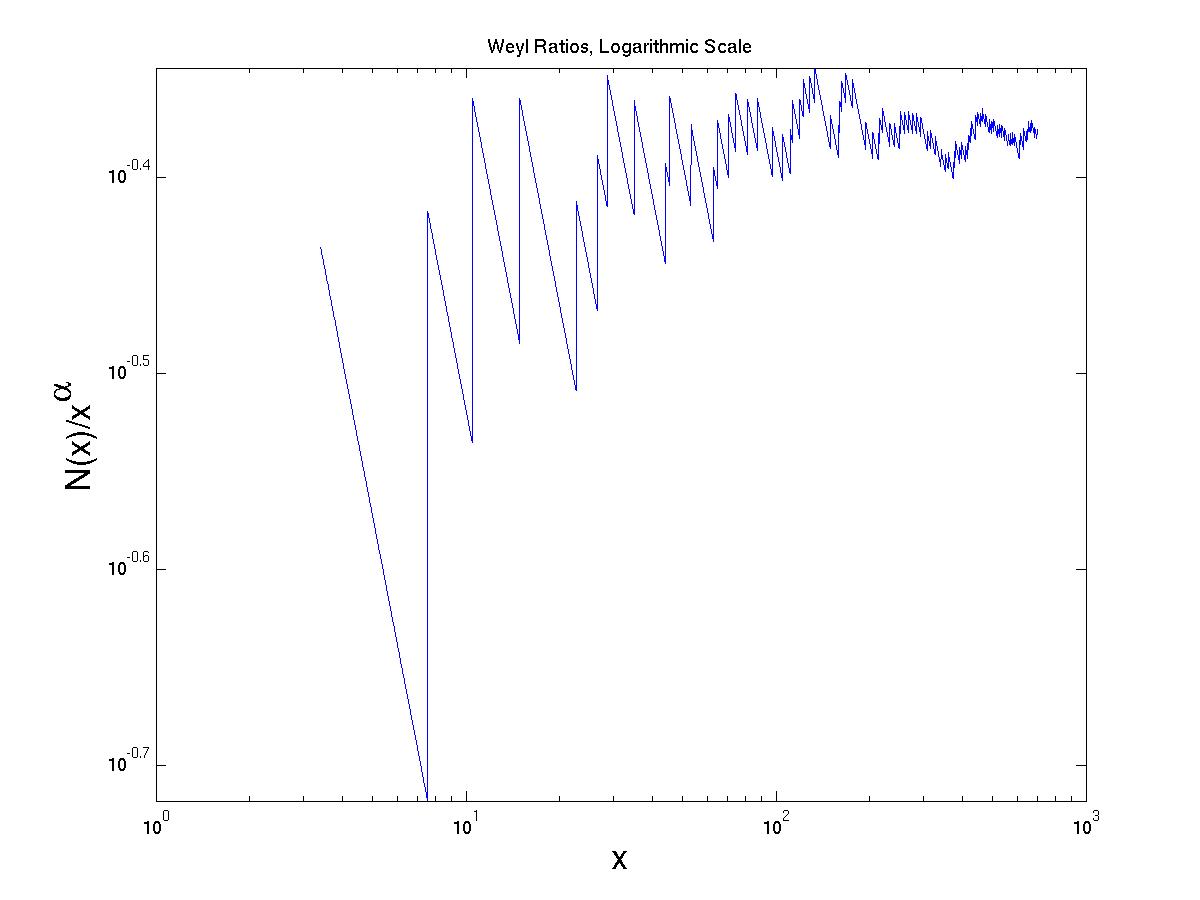}\\
\mbox{Bifurcation at Level $3$, $\alpha=.82086$} & \mbox{Bifurcation at Level $2$, $\alpha=.81975$}
\end{array}
$$
\caption{Weyl Ratios for $j=4$, $k=3$}
\label{figsix11}
\end{figure}

For further comparison of the Weyl ratios, we show those from another trial with $j=4$ and $k=2$, and those from another trial where $j=4$ and $k=3$.  The carpets for the new $j=4$, $k=2$ trial appear in Figure \ref{figsix8} with Weyl ratios in Figure \ref{figsix9}, while the carpets for the new $j=4$, $k=3$ trial appear in Figure \ref{figsix10} with Weyl ratios in Figure \ref{figsix11}.  It is clear that different random choices in the construction make a big difference in the spectrum.  We leave to the future the problem of formulating precise conjectures concerning the spectra of different random carpets.

Acknowledgments: We are grateful to Stacey Goff who contributed to the numerical experiments.

\begin{table}[htbp!]
  \resizebox{\tabfourwda}{!}{
  \begin{tabular}{r|*{4}r|r|rr|rrr||rrr}
    & \multicolumn{4}{p{5cm}}{Original Carpet}
    & \multicolumn{1}{p{1.5cm}}{Bifurc. at Level $4$}
    & \multicolumn{2}{p{2.3cm}}{Bifurcation at Level $3$}
    & \multicolumn{3}{p{3.5cm}}{Bifurcation at Level $2$}
    & \multicolumn{3}{p{3cm}}{Original Carpet Ratios $\lambda_{n}^{j+1}/\lambda_{n}^{j}$}\\
    \hline
    Level: & 1 & 2 & 3 & 4 & 4 & 3 & 4 & 2 & 3 & 4\\
    Refinement: & 2 & 1 & 0 & 0 & 0 & 0 & 0 & 1 & 0 & 0 \\
    \hline
    n\\
  1 &   5.580 &   4.524 &   3.961 &   3.331 &   3.349 &   3.885 &   3.248 &   5.011 &   4.393 &   3.689 &   0.811 &   0.875 &   0.841\\
  2 &   7.666 &   6.734 &   5.914 &   5.008 &   5.009 &   5.963 &   4.990 &   6.384 &   5.528 &   4.639 &   0.878 &   0.878 &   0.847\\
  3 &  18.031 &  15.575 &  13.671 &  11.556 &  11.543 &  13.514 &  11.311 &  15.597 &  13.664 &  11.478 &   0.864 &   0.878 &   0.845\\
  4 &  31.079 &  24.699 &  21.630 &  18.234 &  18.269 &  21.733 &  18.259 &  27.549 &  24.099 &  20.091 &   0.795 &   0.876 &   0.843\\
  5 &  40.933 &  35.373 &  31.097 &  26.152 &  25.983 &  31.230 &  26.123 &  35.484 &  31.262 &  26.392 &   0.864 &   0.879 &   0.841\\
  6 &  46.442 &  37.463 &  32.776 &  27.549 &  27.640 &  32.326 &  27.255 &  38.315 &  33.722 &  28.249 &   0.807 &   0.875 &   0.841\\
  7 &  49.757 &  41.840 &  36.519 &  30.975 &  30.850 &  36.427 &  30.614 &  45.424 &  39.840 &  33.521 &   0.841 &   0.873 &   0.848\\
  8 &  72.354 &  62.389 &  54.211 &  45.450 &  45.767 &  54.977 &  45.924 &  56.607 &  49.171 &  41.245 &   0.862 &   0.869 &   0.838\\
  9 &  88.309 &  74.938 &  65.259 &  54.360 &  54.880 &  65.444 &  54.860 &  70.649 &  61.489 &  51.576 &   0.849 &   0.871 &   0.833\\
 10 &  96.790 &  77.384 &  65.694 &  55.376 &  55.582 &  66.288 &  54.919 &  74.244 &  63.358 &  53.123 &   0.799 &   0.849 &   0.843\\
 11 & 103.355 &  80.835 &  70.391 &  59.301 &  59.323 &  70.224 &  59.343 &  83.948 &  72.797 &  60.697 &   0.782 &   0.871 &   0.842\\
 12 & 106.680 &  96.799 &  83.310 &  70.018 &  69.713 &  84.547 &  70.628 &  91.564 &  78.693 &  65.923 &   0.907 &   0.861 &   0.840\\
 13 & 138.947 & 114.330 &  90.751 &  74.898 &  75.720 &  93.994 &  78.280 & 100.291 &  88.203 &  73.751 &   0.823 &   0.794 &   0.825\\
 14 & 162.163 & 123.498 & 107.533 &  91.169 &  90.934 & 105.799 &  88.143 & 118.054 & 102.261 &  85.535 &   0.762 &   0.871 &   0.848\\
 15 & 162.163 & 124.726 & 109.598 &  92.432 &  91.856 & 109.995 &  92.539 & 133.591 & 114.904 &  95.879 &   0.769 &   0.879 &   0.843\\
 16 & 168.598 & 138.139 & 117.904 &  99.056 &  99.110 & 118.118 &  99.034 & 140.919 & 121.513 & 101.298 &   0.819 &   0.854 &   0.840\\
 17 & 170.390 & 144.519 & 125.422 & 104.987 & 105.280 & 128.532 & 108.234 & 149.625 & 134.221 & 112.025 &   0.848 &   0.868 &   0.837\\
 18 & 183.444 & 155.084 & 131.723 & 111.310 & 111.070 & 133.795 & 111.329 & 159.175 & 141.125 & 117.324 &   0.845 &   0.849 &   0.845\\
 19 & 198.436 & 172.548 & 152.591 & 128.765 & 129.400 & 150.123 & 125.782 & 164.504 & 143.437 & 119.346 &   0.870 &   0.884 &   0.844\\
 20 & 206.000 & 179.213 & 156.432 & 132.112 & 132.276 & 158.088 & 133.030 & 177.161 & 154.685 & 129.358 &   0.870 &   0.873 &   0.845\\
 21 & 214.522 & 188.011 & 162.499 & 137.582 & 137.272 & 164.643 & 138.022 & 183.438 & 158.717 & 132.738 &   0.876 &   0.864 &   0.847\\
 22 & 240.719 & 196.395 & 170.702 & 143.091 & 143.624 & 169.629 & 141.827 & 201.824 & 173.072 & 144.509 &   0.816 &   0.869 &   0.838\\
 23 & 250.543 & 206.038 & 173.663 & 146.420 & 146.806 & 173.727 & 146.243 & 209.485 & 180.972 & 151.168 &   0.822 &   0.843 &   0.843\\
 24 & 257.351 & 218.075 & 188.898 & 159.001 & 159.718 & 186.572 & 154.557 & 218.375 & 189.522 & 158.937 &   0.847 &   0.866 &   0.842\\
 25 & 266.899 & 224.684 & 195.801 & 164.282 & 164.779 & 195.558 & 163.810 & 230.264 & 195.063 & 162.503 &   0.842 &   0.871 &   0.839\\
 26 & 276.814 & 231.084 & 200.575 & 168.720 & 168.320 & 201.320 & 170.096 & 235.294 & 204.404 & 171.042 &   0.835 &   0.868 &   0.841\\
 27 & 279.112 & 253.950 & 209.948 & 176.419 & 174.688 & 220.595 & 183.338 & 241.413 & 207.898 & 173.358 &   0.910 &   0.827 &   0.840\\
 28 & 302.034 & 265.227 & 223.157 & 187.434 & 187.788 & 233.403 & 193.925 & 256.896 & 223.343 & 186.535 &   0.878 &   0.841 &   0.840\\
 29 & 331.424 & 277.522 & 232.146 & 194.488 & 193.716 & 238.930 & 199.095 & 278.109 & 245.174 & 203.467 &   0.837 &   0.836 &   0.838\\
 30 & 339.761 & 284.208 & 242.094 & 204.379 & 203.672 & 246.169 & 205.403 & 292.874 & 257.954 & 215.964 &   0.836 &   0.852 &   0.844\\
 31 & 372.054 & 299.978 & 251.492 & 212.246 & 211.069 & 257.519 & 215.248 & 317.191 & 271.234 & 226.783 &   0.806 &   0.838 &   0.844\\
 32 & 385.306 & 316.612 & 265.554 & 223.592 & 224.233 & 266.844 & 222.676 & 331.605 & 287.096 & 238.995 &   0.822 &   0.839 &   0.842\\
 33 & 393.257 & 324.027 & 278.512 & 233.507 & 232.679 & 278.277 & 231.720 & 342.980 & 301.702 & 254.300 &   0.824 &   0.860 &   0.838\\
 34 & 395.854 & 339.584 & 294.971 & 248.618 & 248.375 & 286.407 & 241.061 & 355.614 & 311.708 & 261.206 &   0.858 &   0.869 &   0.843\\
 35 & 405.837 & 361.323 & 303.434 & 254.390 & 253.273 & 301.936 & 253.741 & 371.636 & 324.228 & 269.291 &   0.890 &   0.840 &   0.838\\
 36 & 420.811 & 364.392 & 315.821 & 263.606 & 264.886 & 313.507 & 260.755 & 380.539 & 327.137 & 273.775 &   0.866 &   0.867 &   0.835\\
 37 & 424.847 & 375.125 & 320.413 & 268.424 & 268.127 & 317.266 & 261.999 & 386.150 & 333.468 & 277.348 &   0.883 &   0.854 &   0.838\\
 38 & 451.376 & 386.611 & 330.029 & 278.419 & 277.801 & 325.821 & 270.028 & 397.166 & 346.083 & 287.709 &   0.857 &   0.854 &   0.844\\
 39 & 462.812 & 405.467 & 338.324 & 282.917 & 283.186 & 339.347 & 282.063 & 409.025 & 359.088 & 301.704 &   0.876 &   0.834 &   0.836\\
 40 & 502.632 & 408.831 & 347.163 & 292.032 & 290.267 & 361.005 & 302.953 & 429.651 & 366.719 & 306.863 &   0.813 &   0.849 &   0.841\\
 41 & 528.577 & 447.618 & 366.670 & 305.564 & 307.685 & 387.191 & 316.739 & 448.759 & 389.395 & 325.524 &   0.847 &   0.819 &   0.833\\
 42 & 545.770 & 452.183 & 388.507 & 323.312 & 324.962 & 388.234 & 323.430 & 456.622 & 395.283 & 329.216 &   0.829 &   0.859 &   0.832\\
 43 & 551.007 & 464.111 & 394.833 & 329.113 & 330.967 & 389.063 & 324.844 & 462.266 & 399.880 & 333.236 &   0.842 &   0.851 &   0.834\\
 44 & 552.842 & 477.614 & 405.344 & 338.488 & 338.553 & 405.026 & 335.237 & 486.916 & 414.255 & 346.218 &   0.864 &   0.849 &   0.835\\
 45 & 564.216 & 494.116 & 413.127 & 345.354 & 345.961 & 416.434 & 344.196 & 495.591 & 429.603 & 356.793 &   0.876 &   0.836 &   0.836\\
 46 & 578.850 & 501.300 & 431.951 & 362.010 & 361.376 & 432.872 & 359.480 & 500.567 & 436.894 & 364.233 &   0.866 &   0.862 &   0.838\\
 47 & 598.075 & 523.765 & 440.593 & 367.808 & 366.624 & 444.307 & 363.054 & 512.076 & 438.755 & 367.034 &   0.876 &   0.841 &   0.835\\
 48 & 613.847 & 551.828 & 449.590 & 374.306 & 377.290 & 449.181 & 369.326 & 547.600 & 457.011 & 378.966 &   0.899 &   0.815 &   0.833\\
 49 & 631.941 & 558.302 & 477.459 & 399.598 & 402.493 & 476.826 & 396.391 & 554.463 & 465.816 & 386.917 &   0.883 &   0.855 &   0.837\\
 50 & 715.621 & 567.952 & 488.904 & 409.795 & 410.350 & 480.033 & 397.220 & 576.594 & 477.564 & 397.180 &   0.794 &   0.861 &   0.838\\
 51 & 725.020 & 592.871 & 499.609 & 416.054 & 418.770 & 488.846 & 401.513 & 587.455 & 492.776 & 410.999 &   0.818 &   0.843 &   0.833\\
 52 & 727.913 & 600.091 & 503.672 & 419.499 & 422.080 & 493.842 & 415.384 & 609.850 & 512.419 & 423.682 &   0.824 &   0.839 &   0.833\\
 53 & 731.705 & 616.830 & 508.602 & 421.296 & 423.380 & 505.115 & 419.245 & 612.367 & 521.155 & 428.834 &   0.843 &   0.825 &   0.828\\
 54 & 735.932 & 631.051 & 524.965 & 431.091 & 444.056 & 522.488 & 436.796 & 631.828 & 538.115 & 449.802 &   0.857 &   0.832 &   0.821\\
 55 & 736.330 & 637.768 & 532.002 & 438.897 & 447.894 & 538.642 & 448.956 & 640.834 & 543.298 & 455.736 &   0.866 &   0.834 &   0.825\\
 56 & 737.370 & 656.446 & 548.081 & 455.651 & 457.588 & 555.001 & 459.264 & 648.905 & 549.244 & 458.336 &   0.890 &   0.835 &   0.831\\
 57 & 763.993 & 671.049 & 563.190 & 467.554 & 468.467 & 567.954 & 471.301 & 668.772 & 568.342 & 472.118 &   0.878 &   0.839 &   0.830\\
 58 & 770.748 & 681.736 & 573.285 & 472.499 & 493.796 & 576.961 & 480.894 & 688.591 & 576.277 & 480.247 &   0.885 &   0.841 &   0.824\\
 59 & 771.833 & 688.547 & 591.230 & 495.228 & 497.748 & 581.500 & 488.252 & 705.282 & 604.405 & 502.457 &   0.892 &   0.859 &   0.838\\
 60 & 772.741 & 699.743 & 601.577 & 508.440 & 502.614 & 596.141 & 494.604 & 712.415 & 609.218 & 509.436 &   0.906 &   0.860 &   0.845\\
    \end{tabular}
    }
    \caption{Carpet Bifurcation Unnormalized Eigenvalues for $j=4$, $k=2$}
    \label{tablesix1}
\end{table}

\begin{table}[htbp!]
  \resizebox{\tabfivewda}{!}{
  \begin{tabular}{r|*{4}r|r|rr|rrr}
    & \multicolumn{4}{p{5cm}}{Original Carpet}
    & \multicolumn{1}{p{1.5cm}}{Bifurc. at Level $4$}
    & \multicolumn{2}{p{2.3cm}}{Bifurcation at Level $3$}
    & \multicolumn{3}{p{3.5cm}}{Bifurcation at Level $2$}\\
    \hline
    Level: & 1 & 2 & 3 & 4 & 4 & 3 & 4 & 2 & 3 & 4\\
    Refinement: & 2 & 1 & 0 & 0 & 0 & 0 & 0 & 1 & 0 & 0 \\
    \hline
    n\\
  1 &   1.000 &   1.000 &   1.000 &   1.000 &   1.000 &   1.000 &   1.000 &   1.000 &   1.000 &   1.000\\
  2 &   1.374 &   1.489 &   1.493 &   1.503 &   1.496 &   1.535 &   1.536 &   1.274 &   1.258 &   1.257\\
  3 &   3.231 &   3.443 &   3.452 &   3.469 &   3.447 &   3.479 &   3.482 &   3.113 &   3.110 &   3.111\\
  4 &   5.570 &   5.460 &   5.461 &   5.473 &   5.455 &   5.595 &   5.621 &   5.498 &   5.485 &   5.445\\
  5 &   7.336 &   7.819 &   7.852 &   7.850 &   7.758 &   8.039 &   8.042 &   7.082 &   7.116 &   7.153\\
  6 &   8.323 &   8.281 &   8.276 &   8.270 &   8.253 &   8.322 &   8.391 &   7.647 &   7.676 &   7.657\\
  7 &   8.917 &   9.249 &   9.221 &   9.298 &   9.211 &   9.377 &   9.425 &   9.066 &   9.068 &   9.086\\
  8 &  12.967 &  13.791 &  13.688 &  13.643 &  13.665 &  14.152 &  14.138 &  11.297 &  11.192 &  11.179\\
  9 &  15.826 &  16.565 &  16.477 &  16.318 &  16.387 &  16.847 &  16.889 &  14.100 &  13.996 &  13.979\\
 10 &  17.346 &  17.106 &  16.587 &  16.623 &  16.596 &  17.064 &  16.907 &  14.817 &  14.421 &  14.399\\
 11 &  18.522 &  17.868 &  17.773 &  17.801 &  17.713 &  18.077 &  18.269 &  16.754 &  16.570 &  16.451\\
 12 &  19.118 &  21.397 &  21.035 &  21.018 &  20.816 &  21.764 &  21.743 &  18.274 &  17.912 &  17.868\\
 13 &  24.901 &  25.272 &  22.914 &  22.483 &  22.609 &  24.196 &  24.099 &  20.016 &  20.076 &  19.989\\
 14 &  29.062 &  27.299 &  27.151 &  27.367 &  27.152 &  27.235 &  27.135 &  23.561 &  23.276 &  23.183\\
 15 &  29.062 &  27.570 &  27.673 &  27.746 &  27.427 &  28.315 &  28.489 &  26.662 &  26.154 &  25.987\\
 16 &  30.215 &  30.535 &  29.770 &  29.734 &  29.593 &  30.406 &  30.488 &  28.124 &  27.658 &  27.456\\
 17 &  30.536 &  31.946 &  31.668 &  31.515 &  31.435 &  33.087 &  33.320 &  29.862 &  30.551 &  30.363\\
 18 &  32.875 &  34.281 &  33.259 &  33.413 &  33.164 &  34.442 &  34.273 &  31.768 &  32.122 &  31.800\\
 19 &  35.562 &  38.141 &  38.528 &  38.652 &  38.637 &  38.645 &  38.723 &  32.831 &  32.648 &  32.348\\
 20 &  36.918 &  39.615 &  39.498 &  39.657 &  39.496 &  40.696 &  40.954 &  35.357 &  35.208 &  35.061\\
 21 &  38.445 &  41.559 &  41.029 &  41.299 &  40.988 &  42.383 &  42.491 &  36.610 &  36.126 &  35.977\\
 22 &  43.140 &  43.413 &  43.101 &  42.952 &  42.885 &  43.666 &  43.662 &  40.279 &  39.394 &  39.168\\
 23 &  44.900 &  45.544 &  43.848 &  43.952 &  43.835 &  44.721 &  45.022 &  41.808 &  41.192 &  40.973\\
 24 &  46.120 &  48.205 &  47.695 &  47.728 &  47.690 &  48.028 &  47.581 &  43.583 &  43.138 &  43.078\\
 25 &  47.831 &  49.666 &  49.438 &  49.313 &  49.201 &  50.341 &  50.430 &  45.955 &  44.399 &  44.045\\
 26 &  49.608 &  51.081 &  50.643 &  50.646 &  50.259 &  51.825 &  52.365 &  46.959 &  46.525 &  46.359\\
 27 &  50.020 &  56.135 &  53.010 &  52.957 &  52.160 &  56.786 &  56.442 &  48.180 &  47.321 &  46.987\\
 28 &  54.128 &  58.628 &  56.345 &  56.263 &  56.071 &  60.083 &  59.701 &  51.271 &  50.836 &  50.559\\
 29 &  59.395 &  61.346 &  58.615 &  58.381 &  57.841 &  61.506 &  61.293 &  55.504 &  55.805 &  55.148\\
 30 &  60.889 &  62.824 &  61.126 &  61.350 &  60.814 &  63.370 &  63.235 &  58.451 &  58.714 &  58.535\\
 31 &  66.676 &  66.310 &  63.499 &  63.711 &  63.023 &  66.291 &  66.265 &  63.304 &  61.737 &  61.467\\
 32 &  69.051 &  69.986 &  67.050 &  67.117 &  66.954 &  68.692 &  68.552 &  66.181 &  65.347 &  64.777\\
 33 &  70.476 &  71.626 &  70.322 &  70.093 &  69.475 &  71.635 &  71.337 &  68.451 &  68.672 &  68.926\\
 34 &  70.942 &  75.064 &  74.477 &  74.629 &  74.162 &  73.728 &  74.212 &  70.972 &  70.949 &  70.797\\
 35 &  72.731 &  79.870 &  76.614 &  76.362 &  75.624 &  77.725 &  78.116 &  74.170 &  73.799 &  72.989\\
 36 &  75.414 &  80.548 &  79.742 &  79.128 &  79.092 &  80.704 &  80.275 &  75.947 &  74.461 &  74.204\\
 37 &  76.138 &  82.921 &  80.901 &  80.574 &  80.060 &  81.672 &  80.658 &  77.067 &  75.902 &  75.173\\
 38 &  80.892 &  85.460 &  83.329 &  83.574 &  82.948 &  83.874 &  83.130 &  79.265 &  78.773 &  77.981\\
 39 &  82.941 &  89.628 &  85.424 &  84.925 &  84.556 &  87.356 &  86.835 &  81.632 &  81.734 &  81.774\\
 40 &  90.078 &  90.371 &  87.655 &  87.661 &  86.670 &  92.931 &  93.266 &  85.748 &  83.471 &  83.172\\
 41 &  94.727 &  98.945 &  92.581 &  91.723 &  91.871 &  99.672 &  97.510 &  89.562 &  88.632 &  88.230\\
 42 &  97.808 &  99.954 &  98.094 &  97.050 &  97.030 &  99.941 &  99.570 &  91.131 &  89.972 &  89.231\\
 43 &  98.747 & 102.591 &  99.691 &  98.792 &  98.823 & 100.154 & 100.005 &  92.258 &  91.019 &  90.320\\
 44 &  99.076 & 105.576 & 102.345 & 101.606 & 101.088 & 104.263 & 103.205 &  97.177 &  94.290 &  93.839\\
 45 & 101.114 & 109.223 & 104.311 & 103.667 & 103.300 & 107.200 & 105.963 &  98.908 &  97.784 &  96.705\\
 46 & 103.737 & 110.811 & 109.063 & 108.667 & 107.903 & 111.431 & 110.668 &  99.902 &  99.444 &  98.722\\
 47 & 107.182 & 115.777 & 111.246 & 110.407 & 109.470 & 114.375 & 111.768 & 102.198 &  99.867 &  99.481\\
 48 & 110.008 & 121.981 & 113.517 & 112.357 & 112.654 & 115.630 & 113.699 & 109.288 & 104.022 & 102.715\\
 49 & 113.251 & 123.412 & 120.554 & 119.949 & 120.180 & 122.746 & 122.031 & 110.658 & 106.026 & 104.870\\
 50 & 128.248 & 125.545 & 123.443 & 123.010 & 122.526 & 123.572 & 122.287 & 115.075 & 108.701 & 107.652\\
 51 & 129.932 & 131.053 & 126.147 & 124.889 & 125.040 & 125.840 & 123.608 & 117.242 & 112.163 & 111.397\\
 52 & 130.451 & 132.649 & 127.172 & 125.923 & 126.028 & 127.126 & 127.878 & 121.712 & 116.634 & 114.835\\
 53 & 131.130 & 136.349 & 128.417 & 126.463 & 126.416 & 130.029 & 129.067 & 122.214 & 118.622 & 116.231\\
 54 & 131.888 & 139.493 & 132.549 & 129.403 & 132.590 & 134.501 & 134.470 & 126.098 & 122.483 & 121.915\\
 55 & 131.959 & 140.977 & 134.325 & 131.746 & 133.736 & 138.659 & 138.214 & 127.896 & 123.663 & 123.523\\
 56 & 132.145 & 145.106 & 138.385 & 136.775 & 136.630 & 142.870 & 141.387 & 129.506 & 125.016 & 124.228\\
 57 & 136.916 & 148.334 & 142.200 & 140.348 & 139.879 & 146.205 & 145.093 & 133.471 & 129.363 & 127.963\\
 58 & 138.127 & 150.696 & 144.749 & 141.833 & 147.442 & 148.523 & 148.046 & 137.427 & 131.169 & 130.166\\
 59 & 138.322 & 152.202 & 149.280 & 148.655 & 148.622 & 149.692 & 150.311 & 140.758 & 137.571 & 136.186\\
 60 & 138.484 & 154.677 & 151.893 & 152.621 & 150.075 & 153.461 & 152.267 & 142.182 & 138.667 & 138.078\\
    \end{tabular}
    }
    \caption{Carpet Bifurcation Normalized Eigenvalues for $j=4$, $k=2$}
    \label{tablesix2}
\end{table}

\begin{table}[htbp!]
  \resizebox{\tabfourwda}{!}{
  \begin{tabular}{r|*{4}r|r|rr|rrr||rrr}
    & \multicolumn{4}{p{5cm}}{Original Carpet}
    & \multicolumn{1}{p{1.5cm}}{Bifurc. at Level $4$}
    & \multicolumn{2}{p{2.3cm}}{Bifurcation at Level $3$}
    & \multicolumn{3}{p{3.5cm}}{Bifurcation at Level $2$}
    & \multicolumn{3}{p{3cm}}{Original Carpet Ratios $\lambda_{n}^{j+1}/\lambda_{n}^{j}$}\\
    \hline
    Level: & 1 & 2 & 3 & 4 & 4 & 3 & 4 & 2 & 3 & 4\\
    Refinement: & 2 & 1 & 0 & 0 & 0 & 0 & 0 & 1 & 0 & 0 \\
    \hline
    n\\
  1 &   7.092 &   4.504 &   3.375 &   2.426 &   2.445 &   3.493 &   2.631 &   6.127 &   4.897 &   3.695 &   0.635 &   0.749 &   0.719\\
  2 &  11.728 &   8.197 &   6.565 &   4.765 &   4.809 &   6.560 &   4.906 &   9.482 &   7.315 &   5.436 &   0.699 &   0.801 &   0.726\\
  3 &  24.546 &  15.759 &  12.132 &   9.007 &   8.969 &  12.351 &   9.397 &  21.222 &  16.397 &  12.210 &   0.642 &   0.770 &   0.742\\
  4 &  30.185 &  18.800 &  15.634 &  11.511 &  11.431 &  15.081 &  11.247 &  23.081 &  17.654 &  13.228 &   0.623 &   0.832 &   0.736\\
  5 &  42.518 &  27.342 &  22.736 &  16.597 &  17.156 &  21.355 &  15.651 &  35.771 &  28.705 &  21.604 &   0.643 &   0.832 &   0.730\\
  6 &  58.544 &  39.332 &  31.633 &  23.474 &  23.415 &  29.450 &  21.692 &  42.024 &  34.860 &  25.822 &   0.672 &   0.804 &   0.742\\
  7 &  61.533 &  48.008 &  39.592 &  29.343 &  28.460 &  37.253 &  26.997 &  51.343 &  38.677 &  28.823 &   0.780 &   0.825 &   0.741\\
  8 &  77.637 &  55.316 &  44.954 &  33.046 &  33.484 &  43.576 &  33.038 &  62.867 &  47.085 &  35.453 &   0.712 &   0.813 &   0.735\\
  9 &  83.257 &  65.486 &  50.292 &  37.106 &  37.063 &  50.448 &  38.068 &  69.055 &  53.226 &  39.577 &   0.787 &   0.768 &   0.738\\
 10 & 104.768 &  73.075 &  55.186 &  39.936 &  40.286 &  56.936 &  42.856 &  82.013 &  62.297 &  45.968 &   0.697 &   0.755 &   0.724\\
 11 & 113.834 &  80.749 &  63.592 &  47.756 &  46.867 &  61.053 &  45.682 &  89.184 &  71.478 &  53.499 &   0.709 &   0.788 &   0.751\\
 12 & 150.330 & 107.208 &  87.607 &  64.387 &  64.396 &  83.453 &  61.747 & 104.399 &  78.918 &  58.162 &   0.713 &   0.817 &   0.735\\
 13 & 162.163 & 114.553 &  93.828 &  68.725 &  68.441 &  87.807 &  65.440 & 121.120 &  97.702 &  73.003 &   0.706 &   0.819 &   0.732\\
 14 & 162.163 & 123.236 &  97.285 &  71.388 &  72.414 &  96.939 &  73.436 & 128.037 & 102.533 &  76.207 &   0.760 &   0.789 &   0.734\\
 15 & 172.077 & 134.430 & 108.788 &  80.883 &  81.781 & 105.258 &  78.237 & 141.405 & 111.031 &  82.675 &   0.781 &   0.809 &   0.743\\
 16 & 175.214 & 142.039 & 112.803 &  85.010 &  83.188 & 110.001 &  79.855 & 148.384 & 120.540 &  89.074 &   0.811 &   0.794 &   0.754\\
 17 & 188.670 & 151.413 & 120.069 &  88.113 &  88.222 & 116.098 &  84.337 & 157.156 & 125.459 &  92.636 &   0.803 &   0.793 &   0.734\\
 18 & 195.969 & 155.785 & 124.187 &  92.302 &  93.366 & 125.759 &  91.453 & 166.793 & 129.710 &  96.150 &   0.795 &   0.797 &   0.743\\
 19 & 202.684 & 172.478 & 133.879 &  97.392 &  98.402 & 139.613 & 103.979 & 175.989 & 137.547 & 102.029 &   0.851 &   0.776 &   0.727\\
 20 & 240.609 & 181.954 & 145.005 & 105.777 & 106.864 & 147.803 & 108.886 & 181.217 & 144.039 & 106.966 &   0.756 &   0.797 &   0.729\\
 21 & 244.457 & 187.959 & 146.514 & 109.927 & 109.130 & 153.475 & 113.767 & 189.592 & 153.756 & 113.071 &   0.769 &   0.780 &   0.750\\
 22 & 274.309 & 195.733 & 161.327 & 116.845 & 117.350 & 158.670 & 118.341 & 215.292 & 165.384 & 119.406 &   0.714 &   0.824 &   0.724\\
 23 & 280.954 & 213.619 & 164.799 & 120.048 & 120.767 & 168.550 & 127.045 & 218.276 & 174.245 & 129.845 &   0.760 &   0.771 &   0.728\\
 24 & 310.772 & 221.652 & 171.028 & 125.428 & 126.050 & 176.933 & 132.774 & 240.206 & 186.285 & 134.818 &   0.713 &   0.772 &   0.733\\
 25 & 315.838 & 234.301 & 178.988 & 132.647 & 130.963 & 181.365 & 135.580 & 264.516 & 204.123 & 154.817 &   0.742 &   0.764 &   0.741\\
 26 & 325.052 & 249.165 & 191.045 & 141.451 & 139.979 & 191.300 & 140.370 & 271.605 & 208.368 & 156.500 &   0.767 &   0.767 &   0.740\\
 27 & 331.424 & 256.361 & 197.974 & 145.517 & 144.028 & 202.974 & 148.117 & 277.958 & 218.970 & 162.376 &   0.774 &   0.772 &   0.735\\
 28 & 360.177 & 265.361 & 201.315 & 149.385 & 146.794 & 211.209 & 159.046 & 285.512 & 224.683 & 167.545 &   0.737 &   0.759 &   0.742\\
 29 & 385.417 & 283.819 & 227.801 & 162.546 & 167.689 & 217.654 & 162.927 & 309.537 & 242.878 & 180.433 &   0.736 &   0.803 &   0.714\\
 30 & 389.957 & 302.005 & 238.307 & 174.063 & 174.823 & 234.133 & 173.773 & 315.867 & 251.600 & 187.139 &   0.774 &   0.789 &   0.730\\
 31 & 404.260 & 313.891 & 253.256 & 181.953 & 186.426 & 234.618 & 175.749 & 322.848 & 265.116 & 192.433 &   0.776 &   0.807 &   0.718\\
 32 & 415.828 & 319.907 & 261.191 & 192.491 & 192.163 & 254.854 & 184.953 & 340.518 & 267.937 & 197.979 &   0.769 &   0.816 &   0.737\\
 33 & 438.883 & 339.911 & 275.572 & 197.205 & 202.084 & 263.163 & 194.662 & 357.545 & 271.540 & 198.556 &   0.774 &   0.811 &   0.716\\
 34 & 446.211 & 351.282 & 280.260 & 208.609 & 207.833 & 268.417 & 200.150 & 360.353 & 289.120 & 213.750 &   0.787 &   0.798 &   0.744\\
 35 & 459.422 & 358.047 & 294.617 & 217.058 & 218.615 & 289.697 & 212.752 & 375.481 & 300.539 & 223.974 &   0.779 &   0.823 &   0.737\\
 36 & 496.933 & 372.237 & 303.439 & 221.374 & 219.487 & 294.416 & 216.619 & 409.773 & 314.193 & 228.308 &   0.749 &   0.815 &   0.730\\
 37 & 497.815 & 399.070 & 310.002 & 227.464 & 231.842 & 305.155 & 225.425 & 420.361 & 322.296 & 239.171 &   0.802 &   0.777 &   0.734\\
 38 & 508.950 & 411.214 & 318.971 & 233.374 & 233.036 & 312.163 & 229.676 & 430.050 & 326.045 & 240.535 &   0.808 &   0.776 &   0.732\\
 39 & 527.137 & 420.241 & 330.335 & 243.651 & 240.021 & 323.612 & 237.958 & 442.270 & 336.285 & 246.673 &   0.797 &   0.786 &   0.738\\
 40 & 551.265 & 426.786 & 333.384 & 249.500 & 243.044 & 332.791 & 248.954 & 446.720 & 343.015 & 251.129 &   0.774 &   0.781 &   0.748\\
 41 & 567.597 & 450.642 & 355.913 & 265.783 & 258.013 & 345.979 & 254.543 & 457.895 & 357.701 & 264.114 &   0.794 &   0.790 &   0.747\\
 42 & 583.940 & 472.465 & 361.077 & 267.320 & 261.429 & 352.508 & 260.712 & 464.320 & 367.546 & 271.991 &   0.809 &   0.764 &   0.740\\
 43 & 597.055 & 484.208 & 368.120 & 272.725 & 269.302 & 356.974 & 266.488 & 480.125 & 380.877 & 280.646 &   0.811 &   0.760 &   0.741\\
 44 & 604.989 & 488.128 & 377.891 & 278.169 & 275.595 & 379.520 & 279.242 & 519.333 & 390.951 & 285.829 &   0.807 &   0.774 &   0.736\\
 45 & 666.885 & 501.510 & 396.592 & 291.488 & 284.793 & 390.723 & 287.972 & 534.943 & 417.099 & 303.615 &   0.752 &   0.791 &   0.735\\
 46 & 696.767 & 528.048 & 410.599 & 301.706 & 303.345 & 402.895 & 293.498 & 542.283 & 422.342 & 309.099 &   0.758 &   0.778 &   0.735\\
 47 & 724.838 & 544.995 & 416.860 & 305.621 & 305.132 & 411.603 & 304.224 & 559.085 & 432.436 & 320.626 &   0.752 &   0.765 &   0.733\\
 48 & 725.020 & 560.630 & 430.081 & 311.371 & 318.362 & 418.850 & 306.840 & 571.554 & 444.148 & 325.198 &   0.773 &   0.767 &   0.724\\
 49 & 731.705 & 564.489 & 442.746 & 321.252 & 325.115 & 437.407 & 323.790 & 593.800 & 457.423 & 335.551 &   0.771 &   0.784 &   0.726\\
 50 & 733.652 & 579.186 & 445.729 & 328.599 & 329.796 & 442.425 & 324.682 & 619.347 & 472.066 & 348.351 &   0.789 &   0.770 &   0.737\\
 51 & 740.114 & 612.913 & 452.752 & 331.363 & 332.623 & 460.158 & 336.138 & 624.425 & 485.492 & 356.737 &   0.828 &   0.739 &   0.732\\
 52 & 742.828 & 630.679 & 475.019 & 341.589 & 343.604 & 468.450 & 344.312 & 653.206 & 500.482 & 372.797 &   0.849 &   0.753 &   0.719\\
 53 & 760.150 & 646.953 & 497.418 & 364.433 & 360.065 & 479.135 & 353.972 & 658.561 & 506.439 & 378.011 &   0.851 &   0.769 &   0.733\\
 54 & 770.961 & 655.027 & 523.642 & 372.127 & 371.947 & 500.853 & 366.851 & 663.004 & 520.723 & 380.865 &   0.850 &   0.799 &   0.711\\
 55 & 772.604 & 686.275 & 532.901 & 385.695 & 377.962 & 513.718 & 377.732 & 689.445 & 521.847 & 388.528 &   0.888 &   0.777 &   0.724\\
 56 & 809.840 & 696.798 & 544.755 & 393.847 & 396.853 & 517.891 & 381.072 & 717.226 & 545.963 & 404.472 &   0.860 &   0.782 &   0.723\\
 57 & 815.252 & 711.641 & 553.161 & 399.410 & 403.815 & 544.312 & 395.388 & 723.980 & 575.366 & 428.062 &   0.873 &   0.777 &   0.722\\
 58 & 840.638 & 722.733 & 572.093 & 408.348 & 412.322 & 549.102 & 400.724 & 751.200 & 582.018 & 429.608 &   0.860 &   0.792 &   0.714\\
 59 & 884.506 & 746.735 & 586.715 & 413.175 & 425.816 & 566.350 & 409.882 & 771.766 & 595.182 & 434.700 &   0.844 &   0.786 &   0.704\\
 60 & 925.356 & 769.117 & 598.816 & 436.475 & 434.588 & 576.353 & 421.201 & 801.163 & 615.474 & 448.107 &   0.831 &   0.779 &   0.729\\
    \end{tabular}
    }
    \caption{Carpet Bifurcation Unnormalized Eigenvalues for $j=4$, $k=3$}
    \label{tablesix3}
\end{table}

\begin{table}[htbp!]
  \resizebox{\tabfivewda}{!}{
  \begin{tabular}{r|*{4}r|r|rr|rrr}
    & \multicolumn{4}{p{5cm}}{Original Carpet}
    & \multicolumn{1}{p{1.5cm}}{Bifurc. at Level $4$}
    & \multicolumn{2}{p{2.3cm}}{Bifurcation at Level $3$}
    & \multicolumn{3}{p{3.5cm}}{Bifurcation at Level $2$}\\
    \hline
    Level: & 1 & 2 & 3 & 4 & 4 & 3 & 4 & 2 & 3 & 4\\
    Refinement: & 2 & 1 & 0 & 0 & 0 & 0 & 0 & 1 & 0 & 0 \\
    \hline
    n\\
  1 &   1.000 &   1.000 &   1.000 &   1.000 &   1.000 &   1.000 &   1.000 &   1.000 &   1.000 &   1.000\\
  2 &   1.654 &   1.820 &   1.945 &   1.964 &   1.967 &   1.878 &   1.865 &   1.548 &   1.494 &   1.471\\
  3 &   3.461 &   3.499 &   3.595 &   3.712 &   3.668 &   3.536 &   3.572 &   3.464 &   3.349 &   3.305\\
  4 &   4.256 &   4.174 &   4.633 &   4.744 &   4.675 &   4.317 &   4.275 &   3.767 &   3.605 &   3.580\\
  5 &   5.995 &   6.071 &   6.737 &   6.840 &   7.016 &   6.113 &   5.950 &   5.839 &   5.862 &   5.847\\
  6 &   8.255 &   8.732 &   9.374 &   9.674 &   9.576 &   8.431 &   8.246 &   6.859 &   7.119 &   6.989\\
  7 &   8.676 &  10.659 &  11.732 &  12.093 &  11.640 &  10.665 &  10.263 &   8.380 &   7.899 &   7.801\\
  8 &  10.947 &  12.281 &  13.321 &  13.619 &  13.694 &  12.475 &  12.559 &  10.262 &   9.616 &   9.595\\
  9 &  11.740 &  14.539 &  14.903 &  15.292 &  15.158 &  14.442 &  14.471 &  11.272 &  10.870 &  10.712\\
 10 &  14.773 &  16.224 &  16.353 &  16.459 &  16.476 &  16.300 &  16.291 &  13.387 &  12.723 &  12.441\\
 11 &  16.051 &  17.928 &  18.844 &  19.681 &  19.167 &  17.478 &  17.366 &  14.557 &  14.597 &  14.479\\
 12 &  21.197 &  23.802 &  25.960 &  26.536 &  26.336 &  23.891 &  23.473 &  17.040 &  16.117 &  15.742\\
 13 &  22.866 &  25.433 &  27.804 &  28.324 &  27.991 &  25.137 &  24.876 &  19.770 &  19.953 &  19.758\\
 14 &  22.866 &  27.361 &  28.828 &  29.421 &  29.616 &  27.752 &  27.916 &  20.899 &  20.940 &  20.625\\
 15 &  24.264 &  29.846 &  32.237 &  33.334 &  33.446 &  30.133 &  29.741 &  23.081 &  22.675 &  22.376\\
 16 &  24.706 &  31.535 &  33.427 &  35.035 &  34.022 &  31.491 &  30.356 &  24.220 &  24.617 &  24.108\\
 17 &  26.604 &  33.617 &  35.580 &  36.314 &  36.081 &  33.236 &  32.060 &  25.652 &  25.622 &  25.072\\
 18 &  27.633 &  34.587 &  36.800 &  38.040 &  38.184 &  36.002 &  34.765 &  27.225 &  26.490 &  26.023\\
 19 &  28.580 &  38.294 &  39.672 &  40.138 &  40.244 &  39.968 &  39.527 &  28.726 &  28.090 &  27.614\\
 20 &  33.927 &  40.398 &  42.969 &  43.593 &  43.705 &  42.313 &  41.392 &  29.579 &  29.416 &  28.950\\
 21 &  34.470 &  41.731 &  43.416 &  45.304 &  44.632 &  43.937 &  43.248 &  30.946 &  31.401 &  30.603\\
 22 &  38.679 &  43.457 &  47.805 &  48.155 &  47.993 &  45.424 &  44.987 &  35.141 &  33.775 &  32.317\\
 23 &  39.616 &  47.428 &  48.834 &  49.475 &  49.391 &  48.252 &  48.295 &  35.628 &  35.585 &  35.143\\
 24 &  43.821 &  49.211 &  50.680 &  51.692 &  51.551 &  50.652 &  50.473 &  39.208 &  38.044 &  36.488\\
 25 &  44.535 &  52.020 &  53.039 &  54.667 &  53.561 &  51.921 &  51.540 &  43.176 &  41.687 &  41.901\\
 26 &  45.834 &  55.320 &  56.612 &  58.296 &  57.248 &  54.765 &  53.361 &  44.333 &  42.554 &  42.357\\
 27 &  46.733 &  56.917 &  58.665 &  59.972 &  58.904 &  58.107 &  56.305 &  45.370 &  44.719 &  43.947\\
 28 &  50.787 &  58.916 &  59.655 &  61.566 &  60.035 &  60.465 &  60.460 &  46.603 &  45.886 &  45.346\\
 29 &  54.346 &  63.013 &  67.503 &  66.990 &  68.581 &  62.310 &  61.936 &  50.524 &  49.602 &  48.834\\
 30 &  54.986 &  67.051 &  70.616 &  71.736 &  71.499 &  67.027 &  66.059 &  51.557 &  51.383 &  50.649\\
 31 &  57.003 &  69.690 &  75.046 &  74.988 &  76.244 &  67.166 &  66.810 &  52.697 &  54.143 &  52.082\\
 32 &  58.634 &  71.026 &  77.398 &  79.331 &  78.590 &  72.959 &  70.309 &  55.581 &  54.719 &  53.583\\
 33 &  61.885 &  75.467 &  81.659 &  81.273 &  82.647 &  75.338 &  73.999 &  58.360 &  55.455 &  53.739\\
 34 &  62.918 &  77.992 &  83.048 &  85.974 &  84.999 &  76.842 &  76.086 &  58.819 &  59.045 &  57.852\\
 35 &  64.781 &  79.494 &  87.303 &  89.455 &  89.409 &  82.934 &  80.876 &  61.288 &  61.377 &  60.619\\
 36 &  70.071 &  82.644 &  89.917 &  91.234 &  89.765 &  84.285 &  82.346 &  66.885 &  64.166 &  61.792\\
 37 &  70.195 &  88.602 &  91.862 &  93.744 &  94.818 &  87.360 &  85.694 &  68.614 &  65.821 &  64.732\\
 38 &  71.765 &  91.298 &  94.519 &  96.180 &  95.306 &  89.366 &  87.310 &  70.195 &  66.586 &  65.101\\
 39 &  74.329 &  93.302 &  97.887 & 100.415 &  98.163 &  92.643 &  90.458 &  72.190 &  68.678 &  66.762\\
 40 &  77.732 &  94.755 &  98.790 & 102.826 &  99.399 &  95.271 &  94.638 &  72.916 &  70.052 &  67.968\\
 41 &  80.035 & 100.052 & 105.466 & 109.536 & 105.521 &  99.047 &  96.763 &  74.740 &  73.051 &  71.483\\
 42 &  82.339 & 104.897 & 106.996 & 110.170 & 106.918 & 100.916 &  99.108 &  75.789 &  75.062 &  73.614\\
 43 &  84.188 & 107.504 & 109.084 & 112.397 & 110.138 & 102.194 & 101.304 &  78.369 &  77.784 &  75.957\\
 44 &  85.307 & 108.374 & 111.979 & 114.641 & 112.712 & 108.649 & 106.152 &  84.768 &  79.842 &  77.360\\
 45 &  94.035 & 111.345 & 117.521 & 120.130 & 116.474 & 111.856 & 109.471 &  87.316 &  85.182 &  82.174\\
 46 &  98.248 & 117.237 & 121.671 & 124.341 & 124.061 & 115.340 & 111.571 &  88.514 &  86.252 &  83.658\\
 47 & 102.206 & 121.000 & 123.526 & 125.955 & 124.792 & 117.833 & 115.648 &  91.257 &  88.314 &  86.778\\
 48 & 102.232 & 124.471 & 127.444 & 128.325 & 130.202 & 119.908 & 116.643 &  93.292 &  90.706 &  88.015\\
 49 & 103.175 & 125.328 & 131.197 & 132.397 & 132.964 & 125.221 & 123.087 &  96.923 &  93.417 &  90.817\\
 50 & 103.449 & 128.591 & 132.081 & 135.425 & 134.879 & 126.657 & 123.426 & 101.093 &  96.407 &  94.281\\
 51 & 104.360 & 136.079 & 134.162 & 136.564 & 136.035 & 131.734 & 127.781 & 101.922 &  99.149 &  96.551\\
 52 & 104.743 & 140.023 & 140.760 & 140.778 & 140.526 & 134.108 & 130.888 & 106.620 & 102.210 & 100.898\\
 53 & 107.186 & 143.637 & 147.398 & 150.193 & 147.258 & 137.166 & 134.560 & 107.494 & 103.427 & 102.309\\
 54 & 108.710 & 145.429 & 155.169 & 153.363 & 152.117 & 143.384 & 139.456 & 108.219 & 106.344 & 103.081\\
 55 & 108.942 & 152.367 & 157.912 & 158.955 & 154.578 & 147.067 & 143.592 & 112.535 & 106.574 & 105.155\\
 56 & 114.192 & 154.703 & 161.425 & 162.315 & 162.303 & 148.262 & 144.862 & 117.069 & 111.499 & 109.471\\
 57 & 114.955 & 157.999 & 163.916 & 164.608 & 165.151 & 155.825 & 150.304 & 118.172 & 117.504 & 115.855\\
 58 & 118.535 & 160.461 & 169.526 & 168.291 & 168.630 & 157.197 & 152.332 & 122.615 & 118.862 & 116.274\\
 59 & 124.721 & 165.790 & 173.859 & 170.281 & 174.149 & 162.134 & 155.814 & 125.972 & 121.551 & 117.652\\
 60 & 130.481 & 170.760 & 177.445 & 179.883 & 177.736 & 164.998 & 160.117 & 130.770 & 125.695 & 121.280\\
    \end{tabular}
    }
    \caption{Carpet Bifurcation Normalized Eigenvalues for $j=4$, $k=3$}
    \label{tablesix4}
\end{table}

\begin{table}[htbp!]
  \resizebox{\tabninewd}{!}{
  \begin{tabular}{r|*{5}r||*{4}r}
    & \multicolumn{5}{p{4cm}}{Eigenvalue Data}
    & \multicolumn{4}{p{4cm}}{Ratios $\lambda_{n}^{j+1}/\lambda_{n}^{j}$}\\
    \hline
    Level: & 1 & 2 & 3 & 4 & 5 & $j=1$ & $j=2$ & $j=3$ & $j=4$\\
    Refinement: & 2 & 1 & 0 & 0 & 0 \\
    \hline
    n\\	
  1 &   7.812 &   5.846 &   5.041 &   3.855 &   3.160 &   0.748 &   0.862 &   0.765 &   0.820\\
  2 &  11.846 &   8.843 &   7.739 &   5.918 &   4.865 &   0.746 &   0.875 &   0.765 &   0.822\\
  3 &  16.892 &  11.188 &   9.621 &   7.327 &   6.019 &   0.662 &   0.860 &   0.762 &   0.821\\
  4 &  31.783 &  21.564 &  18.685 &  14.339 &  11.813 &   0.678 &   0.867 &   0.767 &   0.824\\
  5 &  38.849 &  27.674 &  24.037 &  18.351 &  15.056 &   0.712 &   0.869 &   0.763 &   0.820\\
  6 &  44.579 &  33.739 &  28.815 &  21.870 &  17.953 &   0.757 &   0.854 &   0.759 &   0.821\\
  7 &  66.179 &  53.267 &  45.651 &  34.817 &  28.614 &   0.805 &   0.857 &   0.763 &   0.822\\
  8 &  79.132 &  57.301 &  50.164 &  38.328 &  31.538 &   0.724 &   0.875 &   0.764 &   0.823\\
  9 &  91.235 &  64.613 &  55.715 &  42.305 &  34.727 &   0.708 &   0.862 &   0.759 &   0.821\\
 10 &  93.671 &  71.615 &  60.579 &  45.932 &  37.773 &   0.765 &   0.846 &   0.758 &   0.822\\
 11 & 115.299 &  72.824 &  65.016 &  48.923 &  40.398 &   0.632 &   0.893 &   0.752 &   0.826\\
 12 & 118.662 &  89.806 &  79.457 &  60.983 &  50.215 &   0.757 &   0.885 &   0.768 &   0.823\\
 13 & 162.163 &  95.052 &  83.209 &  63.326 &  51.851 &   0.586 &   0.875 &   0.761 &   0.819\\
 14 & 162.163 & 106.641 &  90.290 &  63.836 &  52.366 &   0.658 &   0.847 &   0.707 &   0.820\\
 15 & 163.646 & 118.545 &  98.672 &  74.716 &  61.516 &   0.724 &   0.832 &   0.757 &   0.823\\
 16 & 170.709 & 125.049 & 102.625 &  78.056 &  63.851 &   0.733 &   0.821 &   0.761 &   0.818\\
 17 & 174.715 & 135.278 & 113.855 &  86.167 &  70.686 &   0.774 &   0.842 &   0.757 &   0.820\\
 18 & 188.946 & 150.973 & 120.531 &  91.075 &  74.673 &   0.799 &   0.798 &   0.756 &   0.820\\
 19 & 201.351 & 161.952 & 138.235 & 105.178 &  86.257 &   0.804 &   0.854 &   0.761 &   0.820\\
 20 & 204.687 & 170.606 & 144.171 & 108.697 &  89.387 &   0.833 &   0.845 &   0.754 &   0.822\\
 21 & 244.350 & 182.087 & 153.859 & 110.986 &  90.735 &   0.745 &   0.845 &   0.721 &   0.818\\
 22 & 247.462 & 182.669 & 155.701 & 115.922 &  95.112 &   0.738 &   0.852 &   0.745 &   0.820\\
 23 & 264.337 & 202.217 & 175.672 & 132.411 & 108.679 &   0.765 &   0.869 &   0.754 &   0.821\\
 24 & 274.203 & 210.504 & 181.596 & 137.161 & 112.420 &   0.768 &   0.863 &   0.755 &   0.820\\
 25 & 281.799 & 218.830 & 186.415 & 141.338 & 115.948 &   0.777 &   0.852 &   0.758 &   0.820\\
 26 & 290.260 & 227.945 & 196.993 & 150.327 & 123.167 &   0.785 &   0.864 &   0.763 &   0.819\\
 27 & 325.832 & 234.838 & 202.526 & 153.692 & 126.141 &   0.721 &   0.862 &   0.759 &   0.821\\
 28 & 331.424 & 266.942 & 224.635 & 172.469 & 141.217 &   0.805 &   0.842 &   0.768 &   0.819\\
 29 & 336.355 & 289.198 & 236.851 & 179.569 & 147.123 &   0.860 &   0.819 &   0.758 &   0.819\\
 30 & 361.265 & 308.120 & 250.806 & 187.350 & 153.399 &   0.853 &   0.814 &   0.747 &   0.819\\
 31 & 384.370 & 312.120 & 257.020 & 194.251 & 160.344 &   0.812 &   0.823 &   0.756 &   0.825\\
 32 & 394.115 & 316.748 & 264.292 & 200.424 & 164.390 &   0.804 &   0.834 &   0.758 &   0.820\\
 33 & 409.812 & 335.456 & 274.591 & 205.929 & 168.413 &   0.819 &   0.819 &   0.750 &   0.818\\
 34 & 411.970 & 343.218 & 289.720 & 216.514 & 177.282 &   0.833 &   0.844 &   0.747 &   0.819\\
 35 & 439.508 & 359.087 & 304.720 & 232.617 & 190.896 &   0.817 &   0.849 &   0.763 &   0.821\\
 36 & 440.834 & 363.915 & 310.972 & 233.960 & 191.294 &   0.826 &   0.855 &   0.752 &   0.818\\
 37 & 460.863 & 375.298 & 321.912 & 241.177 & 196.161 &   0.814 &   0.858 &   0.749 &   0.813\\
 38 & 498.602 & 386.461 & 336.603 & 245.761 & 201.687 &   0.775 &   0.871 &   0.730 &   0.821\\
 39 & 510.810 & 405.711 & 342.497 & 255.459 & 209.461 &   0.794 &   0.844 &   0.746 &   0.820\\
 40 & 523.667 & 413.699 & 345.709 & 264.311 & 217.396 &   0.790 &   0.836 &   0.765 &   0.822\\
 41 & 531.231 & 430.848 & 352.572 & 265.834 & 218.537 &   0.811 &   0.818 &   0.754 &   0.822\\
 42 & 546.687 & 437.950 & 369.987 & 281.332 & 230.437 &   0.801 &   0.845 &   0.760 &   0.819\\
 43 & 559.577 & 440.856 & 377.907 & 286.941 & 235.441 &   0.788 &   0.857 &   0.759 &   0.821\\
 44 & 577.182 & 471.699 & 383.805 & 291.534 & 239.439 &   0.817 &   0.814 &   0.760 &   0.821\\
 45 & 581.677 & 479.787 & 394.602 & 298.709 & 245.557 &   0.825 &   0.822 &   0.757 &   0.822\\
 46 & 598.344 & 498.649 & 409.685 & 313.466 & 255.821 &   0.833 &   0.822 &   0.765 &   0.816\\
 47 & 626.841 & 504.051 & 421.465 & 321.144 & 263.599 &   0.804 &   0.836 &   0.762 &   0.821\\
 48 & 633.797 & 518.365 & 430.774 & 326.844 & 268.188 &   0.818 &   0.831 &   0.759 &   0.821\\
 49 & 703.078 & 527.983 & 439.453 & 333.442 & 272.346 &   0.751 &   0.832 &   0.759 &   0.817\\
 50 & 717.798 & 547.394 & 449.694 & 334.683 & 273.914 &   0.763 &   0.822 &   0.744 &   0.818\\
 51 & 725.020 & 563.709 & 453.597 & 344.470 & 283.262 &   0.778 &   0.805 &   0.759 &   0.822\\
 52 & 731.705 & 579.899 & 471.292 & 358.560 & 294.274 &   0.793 &   0.813 &   0.761 &   0.821\\
 53 & 737.699 & 592.267 & 491.538 & 371.267 & 303.822 &   0.803 &   0.830 &   0.755 &   0.818\\
 54 & 741.151 & 605.924 & 508.950 & 382.992 & 314.146 &   0.818 &   0.840 &   0.753 &   0.820\\
 55 & 745.273 & 622.510 & 512.829 & 387.189 & 317.207 &   0.835 &   0.824 &   0.755 &   0.819\\
 56 & 765.924 & 648.652 & 543.566 & 410.317 & 336.800 &   0.847 &   0.838 &   0.755 &   0.821\\
 57 & 769.461 & 671.205 & 547.092 & 413.608 & 339.283 &   0.872 &   0.815 &   0.756 &   0.820\\
 58 & 771.967 & 672.947 & 553.461 & 418.171 & 343.727 &   0.872 &   0.822 &   0.756 &   0.822\\
 59 & 784.362 & 697.036 & 579.292 & 434.757 & 356.559 &   0.889 &   0.831 &   0.750 &   0.820\\
 60 & 819.314 & 722.128 & 598.492 & 448.196 & 367.600 &   0.881 &   0.829 &   0.749 &   0.820\\
    \end{tabular}
    }
    \caption{Carpet Mixed Unnormalized Eigenvalues and Ratios, 23232}
    \label{tablesix5}
\end{table}

\begin{table}[htbp!]
  \resizebox{\tabtenwd}{!}{
  \begin{tabular}{r|*{5}r}
    \multicolumn{6}{p{4cm}}{blah blah}\\
    \hline
    Level: & 1 & 2 & 3 & 4 & 5 \\
    Refinement: & 2 & 1 & 0 & 0 & 0 \\
    \hline
    n\\
1 &   1.000 &   1.000 &   1.000 &   1.000 &   1.000\\
2 &   1.516 &   1.513 &   1.535 &   1.535 &   1.540\\
3 &   2.162 &   1.914 &   1.908 &   1.901 &   1.905\\
4 &   4.068 &   3.689 &   3.707 &   3.720 &   3.739\\
5 &   4.973 &   4.734 &   4.768 &   4.761 &   4.765\\
6 &   5.706 &   5.771 &   5.716 &   5.673 &   5.681\\
7 &   8.471 &   9.112 &   9.056 &   9.032 &   9.055\\
8 &  10.129 &   9.802 &   9.951 &   9.943 &   9.981\\
9 &  11.679 &  11.052 &  11.052 &  10.975 &  10.990\\
10 &  11.991 &  12.250 &  12.017 &  11.916 &  11.954\\
11 &  14.759 &  12.457 &  12.897 &  12.692 &  12.785\\
12 &  15.190 &  15.362 &  15.762 &  15.820 &  15.892\\
13 &  20.758 &  16.259 &  16.506 &  16.428 &  16.409\\
14 &  20.758 &  18.242 &  17.911 &  16.561 &  16.572\\
15 &  20.948 &  20.278 &  19.574 &  19.383 &  19.468\\
16 &  21.852 &  21.390 &  20.358 &  20.249 &  20.207\\
17 &  22.365 &  23.140 &  22.585 &  22.354 &  22.370\\
18 &  24.186 &  25.825 &  23.910 &  23.627 &  23.632\\
19 &  25.774 &  27.703 &  27.422 &  27.286 &  27.298\\
20 &  26.201 &  29.183 &  28.599 &  28.198 &  28.288\\
21 &  31.278 &  31.147 &  30.521 &  28.792 &  28.715\\
22 &  31.677 &  31.247 &  30.887 &  30.073 &  30.100\\
23 &  33.837 &  34.591 &  34.848 &  34.351 &  34.394\\
24 &  35.100 &  36.008 &  36.023 &  35.583 &  35.577\\
25 &  36.072 &  37.432 &  36.979 &  36.666 &  36.694\\
26 &  37.155 &  38.991 &  39.078 &  38.998 &  38.979\\
27 &  41.709 &  40.171 &  40.175 &  39.871 &  39.920\\
28 &  42.425 &  45.662 &  44.561 &  44.743 &  44.691\\
29 &  43.056 &  49.469 &  46.984 &  46.584 &  46.560\\
30 &  46.244 &  52.706 &  49.752 &  48.603 &  48.546\\
31 &  49.202 &  53.390 &  50.985 &  50.393 &  50.744\\
32 &  50.449 &  54.182 &  52.428 &  51.995 &  52.024\\
33 &  52.459 &  57.382 &  54.471 &  53.423 &  53.298\\
34 &  52.735 &  58.710 &  57.472 &  56.169 &  56.104\\
35 &  56.260 &  61.424 &  60.447 &  60.346 &  60.413\\
36 &  56.430 &  62.250 &  61.688 &  60.695 &  60.539\\
37 &  58.994 &  64.197 &  63.858 &  62.567 &  62.079\\
38 &  63.824 &  66.107 &  66.772 &  63.756 &  63.828\\
39 &  65.387 &  69.399 &  67.941 &  66.272 &  66.288\\
40 &  67.033 &  70.766 &  68.578 &  68.569 &  68.799\\
41 &  68.001 &  73.699 &  69.940 &  68.964 &  69.160\\
42 &  69.980 &  74.914 &  73.394 &  72.984 &  72.926\\
43 &  71.630 &  75.411 &  74.966 &  74.439 &  74.510\\
44 &  73.883 &  80.687 &  76.136 &  75.631 &  75.775\\
45 &  74.459 &  82.071 &  78.277 &  77.492 &  77.711\\
46 &  76.592 &  85.297 &  81.269 &  81.320 &  80.960\\
47 &  80.240 &  86.221 &  83.606 &  83.312 &  83.421\\
48 &  81.130 &  88.670 &  85.453 &  84.791 &  84.873\\
49 &  89.999 &  90.315 &  87.174 &  86.503 &  86.189\\
50 &  91.883 &  93.635 &  89.206 &  86.825 &  86.686\\
51 &  92.807 &  96.426 &  89.980 &  89.364 &  89.644\\
52 &  93.663 &  99.195 &  93.490 &  93.019 &  93.129\\
53 &  94.431 & 101.311 &  97.507 &  96.315 &  96.150\\
54 &  94.872 & 103.647 & 100.961 &  99.357 &  99.418\\
55 &  95.400 & 106.484 & 101.730 & 100.446 & 100.387\\
56 &  98.044 & 110.956 & 107.827 & 106.446 & 106.587\\
57 &  98.496 & 114.814 & 108.527 & 107.300 & 107.373\\
58 &  98.817 & 115.112 & 109.790 & 108.483 & 108.779\\
59 & 100.404 & 119.232 & 114.914 & 112.786 & 112.840\\
60 & 104.878 & 123.525 & 118.723 & 116.273 & 116.334\\
    \end{tabular}
    }
    \caption{Carpet Mixed Normalized Eigenvalues 23232}
    \label{tablesix6}
\end{table}

\clearpage


Department of Mathematical Sciences, George Mason University, Fairfax, VA  22030

\textit{E-mail address}: \verb!tberry@gmu.edu!

Department of Mathematics, Cornell University, Ithaca, NY 14850-4201

\textit{E-mail address}: \verb!smh82@cornell.edu!

Department of Mathematics, Cornell University, Ithaca, NY 14850-4201

\textit{E-mail address}: \verb!str@math.cornell.edu!

\end{document}